%% file: MechOpInf.tex
\documentclass[%
  onecolumn
]{mpi2015-cscpreprint}

\usepackage[american]{babel}
\usepackage{tikz}
\usepackage{pgfplots}

\usepackage{graphicx}

\usepackage{amssymb}
\usepackage{amsthm}
\usepackage{amsmath}

\usepackage{caption}
\usepackage{subcaption}

\newtheorem{theorem}{Theorem}
\newtheorem{assumption}{Assumption}
\newtheorem*{remark}{Remark}

\usepackage[normalem]{ulem}

\usepackage{todonotes}
\usepackage{cleveref}
\usepackage[software, hardware]{mymacros}

\begin{document}
  
%%%%%%%%%%%%%%%%%%%%%%%%%%%%%%%%%%%%%%%%%%%%%%%%%%%%%%%%%%%%%%%%%%%%%%%%%%%%%%%%
% PAPER INFORMATION.                                                           %
%%%%%%%%%%%%%%%%%%%%%%%%%%%%%%%%%%%%%%%%%%%%%%%%%%%%%%%%%%%%%%%%%%%%%%%%%%%%%%%%

\title{An Operator Inference Oriented Approach for Mechanical Systems}
  
\author[$\ast$]{Yevgeniya Filanova}
\affil[$\ast$]{Otto-von-Guericke University Magdeburg, Max Planck Institute for Dynamics of Complex Technical Systems Magdeburg\authorcr
  \email{filanova@mpi-magdeburg.mpg.de}, \orcid{0000-0002-8599-3747}}
  
\author[$\dagger$]{Igor Pontes Duff}
\affil[$\dagger$]{Max Planck Institute for Dynamics of Complex Technical Systems\authorcr
  \email{pontes@mpi-magdeburg.mpg.de}, \orcid{0000-0001-6433-6142}}

\author[$\ddagger$]{Pawan Goyal}
\affil[$\ddagger$]{Max Planck Institute for Dynamics of Complex Technical Systems\authorcr
	\email{goyalp@mpi-magdeburg.mpg.de}, \orcid{0000-0003-3072-7780}}

\author[$\dagger\dagger$]{\\ Peter Benner}
\affil[$\dagger\dagger$]{Max Planck Institute for Dynamics of Complex Technical Systems, Otto-von-Guericke University Magdeburg\authorcr
	\email{benner@mpi-magdeburg.mpg.de}, \orcid{0000-0003-3362-4103}}
  
\shorttitle{An Operator Inference Oriented Approach for Mechanical Systems}
\shortauthor{Y. Filanova et al.}
\shortdate{}
  

%\msc{MSC1, MSC2, MSC3}
\abstract{%
Model-order reduction techniques allow the construction of low-dimensional surrogate models that can accelerate engineering design processes. Often, these techniques are intrusive, meaning that they require direct access to underlying high-fidelity models. Accessing these models is laborious or may not even be possible in some cases. Therefore, there is an interest in developing non-intrusive model reduction techniques to construct low-dimensional models directly from simulated or experimental data. In this work, we focus on a recent data-driven methodology, namely operator inference, that aims at inferring the reduced operators using only trajectories of high-fidelity models. We present an extension of operator inference for mechanical systems, preserving the second-order structure. We also study a particular case in which complete information about the external forces is available. In this formulation, the reduced operators having certain properties inspired by the original system matrices are enforced by adding constraints to the optimization problem. We illustrate the presented methodology using three numerical examples.}

\novelty{In this work, we extend the operator inference methodology to the second-order system structure, considering two cases. In the first case, we obtain the system operators from the input, state, and derivative information and prove their asymptotic closeness to the intrusive reduced operators. In the second case, we assume the external force information to be available and obtain the system operators, preserving their symmetric positive definite properties.}

\maketitle
%%%%%%%%%%%%%%%%%%%%%%%%%%%%%%%%%%%%%%%%%%%%%%%%%%%%%%%%%%%%%%%%%%%%%%%%%%%%%%%%
% PAPER CONTENT.                                                               %
%%%%%%%%%%%%%%%%%%%%%%%%%%%%%%%%%%%%%%%%%%%%%%%%%%%%%%%%%%%%%%%%%%%%%%%%%%%%%%%%
\section{Introduction}%
\label{sec:intro}
Mathematical models of mechanical systems describe their dynamic behaviors and robustness, allowing to anticipate the state of the system under the influence of certain external factors.
Mechanical models can be designed in various ways, depending on the goals pursued and the system type. Dynamic behavior of interconnected rigid or flexible bodies can be analyzed using multibody system formalism \cite{EicF98}. It is widely used in robotics, vehicle dynamics, and for different types of mechanisms to characterize the motion, e.g., to obtain trajectories, critical speeds, etc. The modeling is based on representing a given system as a number of solid bodies which are connected with joints or force elements. The governing system of ordinary differential (-algebraic) equations is derived using Lagrange's equation followed via the D'Alembert principle.

On the other hand, if the dynamic behavior of a continuous object is of interest, methods from solid mechanics can be utilized. They allow one to identify the displacements, inner stresses, and strains of the structure \cite{Alt15}. 
Considering the general physical principles common to all media, such as the balance of energy, the conservation of mass and momentum, etc., the governing equations are often derived either in integral or differential form. The latter form is essential for most structural analysis problems. It comes as a partial differential balance equation, which is assumed to be satisfied at every point of the field of interest. 
The central part of continuum mechanics consists of the additional constitutive equations, which define the material law. Together with the local balance equation, they allow to completely describe the inner stress-strain state of an object \cite{LaiRetal10}. In practice,  numerical solutions of the governing partial differential equations (PDEs) are arguably most often computed by the finite element method (FEM) \cite{ZieTD14}, which provides a spatial discretization of the solution field and leads to a system of second-order ODEs with specific mechanical properties.

All these are accompanied by the development of new dynamic and material models, solution methods, simulation software, and at the same time---model-order reduction (MOR) methods.
Increasing simulation costs while carrying out engineering design gives rise to the necessity of having surrogate models with lower complexity yet acceptable accuracy. The construction of the lower dimensional models is typically done by projection-based MOR methods.
The main idea is to find a low-dimensional subspace of solution-trajectories and project the system operators onto these subspaces; see, e.g., \cite{morHaa17,morBenMS05,morKunV01} for the details. 

There exist many well-known reduction methods that can be efficiently applied for mechanical systems, such as modal truncation \cite{morDav66}, moment matching \cite{morFre03,morBeaG05,morBenGP19}, and balanced truncation \cite{morMoo79,morChaLVetal06}. These methods rely on constructing projection matrices with a particular focus. For instance, balanced truncation aims at determining the projection matrices containing the subspaces that are easy to reach as well as easy to observe. In \cite{morSaaSW19}, an overview and a comparison of many such methods for linear mechanical systems are provided with applications to a high-dimensional robotic fishtail model. 
%_______________________________________________________________________________________
%POD
It is worthwhile highlighting a snapshot-based approach, namely Proper Orthogonal Decomposition (POD), where a projection matrix or reduced basis is constructed from the state snapshots of the full model.
This method utilizes an orthogonal basis for representing the given data in the least squares optimal sense \cite{morKunV01,morKunV08,morBerHL93,morLuJCetal19}.

%This makes POD attractive to partly use in the non-intrusive MOR methods in order to learn low-dimensional models.

%_______________________________________________________________________________________
%Master-slave reduction techniques
In the engineering literature, common reduction methods are based on dividing the generalized system coordinates into \emph{master} and \emph{slave} coordinates. 
This interpretation reflects the intuitive background of all reduction methods --- that is, some parts of the system of equations may be unimportant for the system dynamics and thus can be omitted. 
Historically, Guyan reduction \cite{morGuy65} was the first important technique in this category. This method is also known as \textit{static condensation} because it does not take into account dynamical effects and provides exact results for static simulations. For dynamical simulation, meaningful results are possible only for the loading frequency range close to the lowest eigenfrequencies of the system; otherwise, the results are too stiff \cite{FloPS14}.
Guyan reduction forms the basis for other more advanced methods, such as Craig-Bampton reduction~\cite{morCraBM68}, Improved Reduction System (IRS)~\cite{morGor94}, and System Equivalent Reduction Expansion Process (SEREP)~\cite{morMarLetal20}. 
These methods have improved accuracy due to the consideration of the eigenmodes of the omitted system as in Craig-Bampton reduction, or due to the approximation of the inertia forces as in IRS and SEREP.

%_______________________________________________________________________________________
%non-intrusive way

All these mentioned methods require access to the system operators. Thus, these methods are referred to as \emph{intrusive} ones. However, in many scenarios, obtaining a full-order model in an explicit form can be very laborious or may not even be possible in many scenarios. Experimental measurements can also characterize a mechanical system, where the actual model behind the experiment may be unknown. Not only these, but very often, simulations of structural, dynamical processes are done via commercial software, and the governing equations are impossible to extract. Therefore, there is considerable interest in constructing potentially low-dimensional models in a \emph{non-intrusive way} using only data that are either obtained using simulations or experiments. In this process, we explicitly eliminate the need for the full-order model but leverage the model hypothesis, which can either be known empirically or given by experts. %_____________________________________________________________________________________
%Examples of non-intrusive methods for mechanical systems.

In recent times, many non-intrusive reduced-order methodologies have been developed. Often, linear dynamical models can be learned using data either obtained in the time domain or frequency domain. The construction of reduced-order models using frequency domain data has been originally developed for first-order and extended to second-order systems (that often arise in mechanical systems): the Loewner framework \cite{morMayA07,morPonGB20a}, the vector-fitting \cite{morGusS99, morWerGG21}, and the AAA algorithm \cite{NakST18,morGosG20} are instances of frequency-domain reduced-order modeling approaches.
There exist several methodologies to learn models from time-domain data. A widely used method for learning discrete-time systems is Dynamic Mode Decomposition (DMD) \cite{morSch10,morCheTR12,morTuRLetal14}, which is an attractive reduction technique related to Koopman operator approximation. 
%DMD constructs the best-fit (least-square) linear dynamical system to the nonlinear dynamical system generating the data (?). 
%The DMD method approximates the modes of the Koopman operator. 
%The main idea of the Koopman operator is lifting the nonlinear dynamics to the linear but infinite-dimensional regime.
The basis of this method is collecting data from a dynamical system and solving a minimization problem to find the linear system operator. 

Another method introduced for first-order parametric systems in \cite{morPehW16} is operator inference which uses hypotheses based on the structure of the PDE level. The essence of the method is learning the unknown operators using the data compressed to a low-dimensional subspace, followed by solving a least-squares problem. Several extensions of operator inference to parametric and nonlinear systems can be found in \cite{morBenGKetal20,morQiaKPetal20,morYilGBetal20}.
% incompressible flow \cite{morBenGHetal20}. 
Although the operator inference approach was developed for continuous-time systems, it shares an analogy to the DMD approaches. %However, operator inference distinguishes a different approach for constructing nonlinear reduced order systems and is developed for the application to continuous-time systems.

%A common issue, which is mentioned already in the original paper on operator inference, is that the optimization problem is usually not well-posed due to the ill-conditioning of the matrix, assembled from data. Therefore, the effective solution to the operator inference optimization problem is closely related to the development of regularization techniques. Various approaches are usually based on Tikhonov regularization (also known as rigid regression) \cite{TikGetal77}. The solution norm is restricted, using some penalty parameters, thus preventing the solution from "blowing up". To achieve a compromise between the solution accuracy and regularity, some authors propose physics-informed regularization for nonlinear dynamical systems \cite{morSawKP21,morMcqHW21}. Considering the different physical nature of the linear and nonlinear system operators, they are assigned to the different "regularization levels", characterized by different penalty parameters.

Most operator inference methods focus on learning first-order ODE systems with a prior hypothesis on the form of the model. However, mechanical systems are distinguished by the second-order specific ODE structure, where system matrices also have a physical meaning. Although second-order ODE systems can be transformed into their first-order companion form, it,  first of all, leads to the system being twice as large as the original one; secondly, a subsequent naive reduction of these systems not only violates its original structure but also can lead to non-physical behavior. Therefore, we focus on preserving second-order structures in the learning process to obtain better interpretability. 
We mention the recent attempt in the direction in \cite{morShaK22}, where the operator inference methodology is described for Lagrangian mechanical models. Therein, the Lagrangian approach to derive the governing equations is presented, together with the formulation of operator inference that preserves the second-order structure and symmetric positive definite (s.p.d.) properties of operators. The methodology in \cite{morShaK22} is presented for the particular case, when the reduced system mass matrix is equal to the identity. 
%_______________________________________________________________________________________
%In this paper we are going to present the work that we did in parallel, focusing on continuum mechanics. (Description: recover FEM formulation, classical POD intrusive reduction and explain methodology operator inference).
In this paper, we present the work in a similar direction. The operator inference procedure is tailored to the mechanical system structure, focusing on the data obtained from the FEM simulations. Firstly, we propose an extension of the operator inference approach to obtain second-order dynamics. We discuss the connection between the inferred operators and the matrices obtained via intrusive POD reduction. Then we tailor the learning process for the case when the external loads are completely known and develop the operator inference approach with additional constraints in order to enforce the reduced operators to be symmetric positive definite \footnote{We note that the work presented in \cite{morShaK22} paralleled our research, of which a first idea was contained in \cite{morBetal}. This development happened independently without both groups knowing of the work of the other.}. 
%The simulation data for the learning phase is assumed to be given
%from the simulation, including the derivative data, which can be provided from a second-order time integrator, typically used in FEM-Software. \todo{We can skip this line in the introduction in the paper at the appropriate place. }

The remainder of the paper is organized as follows. \Cref{sec:fem} briefly presents an overview of continuum mechanics and the derivation of FEM governing equations. \Cref{sec:setup} gives the information about the available data and the time-integration algorithm.  \Cref{sec:pod} briefly explain the intrusive data-driven POD method.  \Cref{sec:opinf} represents the operator inference for mechanical systems in the simplest form and its constrained version. Finally, numerical results are presented in \Cref{sec:num} to illustrate the proposed methodologies. We provide our outlook in  \Cref{sec:conc}.

%%%%%%%%%%%%%%%%%%%%%%%%%%%%%%%%%%%%%%%%%%%%%%%%%%%%%%%%%%%%%%%%%%%%%%%%%%%%%%%%
%%%%%%%%%%%%%%%%%%%%%%%%%%%%%%%%%%%%%%%%%%%%%%%%%%%%%%%%%%%%%%%%%%%%%%%%%%%%%%%%
%_____________________________________________________________________________________________
% MAIN BODY
%_____________________________________________________________________________________________
\section{Continuum mechanics and finite element formulation}
\label{sec:fem}

In this section, we shall briefly review the equation of motion from the continuum solid mechanics viewpoint. The primary interest of solid mechanics is the response of an object to the forces that are acting on it, namely, the identification of the displacement field and stress-strain state. All the characteristic quantities are connected through the kinematic and constitutive relations, and can be found in the solution of local impulse balance PDEs. These concepts are explained in the continuum mechanics literature \cite{LaiRetal10, BerG15,Alt15}.

Next, we present the spatial discretization of the solution domain using the finite element method (FEM). It is a widely used approach for solving structural mechanics problems and is implemented in many powerful simulation packages, which are predominantly used in engineering practice.
More detailed description of the FEM can be found, e.g., in  \cite{ZieTD14, BorCRetal12, LiuMBetal13}. In this paper, we focus on the small deformation theory and linear elastic material law, which cover a wide range of structural mechanics problems.
%Otherwise it would require a more detailed immersion in the finite strain theory.

%\subsection{Basic continuum mechanics relations}
%\label{subsec:continuum}

Consider an object that is exposed to external forces. The various displacements of the body are described at each point by values in the corresponding coordinate directions $x_1, x_2$, and $x_3$ gathered in the displacement vector
\begin{equation*}
\boldsymbol{x}^\top = \left[x_1, ~x_2,~x_3\right].
\end{equation*}
If the body is not rigid, displacements appear together with the changes in size and shape of an object, called deformations or strains. 
We assume that deformations are sufficiently small (less than 5\%) and connected with displacements via kinematic relations, forming a symmetric Cauchy-strain tensor $\boldsymbol{\varepsilon}$ as follows:
\begin{equation}\label{eq:kinem}
\boldsymbol{\varepsilon} = \frac{1}{2} \left( \nabla \boldsymbol{x} + \left(\nabla \boldsymbol{x}\right)^\top \right). 
\end{equation}
For more significant deformations, the symmetry of the strain tensor cannot be assured because of the different formulations in the Lagrangian and Eulerian coordinate systems \cite{Alt15}. Thus, in these scenarios, it is essential to use the finite strain theory, which is the out of the scope of this paper. Finally, an important part of structural analysis is the identification of stresses as a reaction to external and internal loads.
A stress vector $\boldsymbol{t}$ is defined as an inner force $\boldsymbol{f}$ acting in an imaginary cut of a body on an arbitrary small area $\Delta S$:
\begin{equation*}
 \boldsymbol{t} = \lim_{\Delta S \to 0} \frac{\Delta \boldsymbol{f}}{\Delta S} = \frac{\mathrm{d} \boldsymbol{f}}{\mathrm{d}S}.
\end{equation*}
By choosing the vectorial basis for each plane such that the first axis coincides with the normal vector $\boldsymbol{n}$ to the plane, and the second and third axes are two mutually orthogonal
vectors, the stress vector can be presented with three corresponding components. As a result, it forms a second-order Cauchy stress tensor $\sigma$, which describes a three-dimensional stress state in a point of the solid.
%\begin{equation*}
% \boldsymbol{\sigma} = \begin{bmatrix}
%           \sigma_{11} & \sigma_{12} & \sigma_{13} \\
%           \sigma_{21} & \sigma_{22} & \sigma_{23} \\
%           \sigma_{31} & \sigma_{32} & \sigma_{33} 
%          \end{bmatrix}.
%\end{equation*}
In fact, there are only six independent stress components due to the symmetry of the stress tensor following the balance of the rotational momentum.

%Stress tensor and stress vector are related by the following equation.
%\begin{equation}
% \boldsymbol{t} = \boldsymbol{n} \cdot \sigma
%\end{equation}

%At this point we have defined all unknown quantities that form the main problem of solid continuum mechanics. %\red{igor: Maybe some sentence here can be condensed.}
%Note that the symmetric property of stress and strain tensors hold only for small deformation theory. More complicated stress and strain measures are proposed in finite strain theory with derivation using displacement gradient and distinguishing between the Eulerian and Lagrangian coordinate systems \cite{Alt15}.

To find the unknown quantities, we consider the balance of momentum of a body $\mathcal{B}$ with boundary $\partial \mathcal{B}$ in its current configuration and
denote $\boldsymbol{g}$ as gravity acceleration and $\rho$ as mass density. The balance of the momentum postulates that the overall impulse by the deformation of a body is equal to the sum of all surface and volume forces acting on it, i.e.,

\begin{equation}
 \int_{\partial \mathcal{B}} \boldsymbol{n} \cdot \sigma \, \mathrm{d} \boldsymbol{x} + \int_{\mathcal{B}} \rho \boldsymbol{g} \, \mathrm{d} \boldsymbol{x} = \int_{\mathcal{B}} \rho \ddot{\boldsymbol{x}} \, \mathrm{d} \boldsymbol{x}.
 \label{eq:balance}
\end{equation}
Note that the stress vector is substituted by its relation to the stress tensor $\boldsymbol{t} = \boldsymbol{n} \cdot \sigma$. Applying the divergence theorem and using the fundamental assumption that the identity \eqref{eq:balance} must hold for each subpart of the body, we get the local impulse balance equation as follows:

\begin{equation}
 \nabla^\top \cdot \sigma + \rho \boldsymbol{g} = \rho \ddot{\boldsymbol{x}}.
 \label{eq:pdegov}
\end{equation}
%
%\subsection{Finite-Elements}
%
To solve \eqref{eq:pdegov}, appropriate numerical methods are needed. Arguably, the most popular approach for this is FEM, in which the main idea is to discretize the spatial domain into a finite number of simpler and smaller parts, namely finite elements, transforming the infinite-dimensional problem into a finite-dimensional one. The bridge to the finite elements is the weak formulation of~\eqref{eq:pdegov}. It requires the multiplication of the governing equation \eqref{eq:pdegov} with the virtual displacement $\delta \boldsymbol{x}$ and integration over the domain $\mathcal{B}$.
%\begin{equation}
% \int_{\mathcal{B}} \delta \boldsymbol{x} ^\top (\nabla^\top \cdot \sigma + \rho \boldsymbol{g} - \rho \ddot{\boldsymbol{x}}) \; \mathrm{d} \boldsymbol{x} = 0
%\end{equation}
Applying the divergence theorem once more, we get the weak form of the equation of motion in the current configuration:

\begin{equation}
 \label{eq:weak}
 \int_{\mathcal{B}} \left( \rho \, (\delta \boldsymbol{x}) ^\top \ddot{\boldsymbol{x}} + (\nabla \cdot \delta \boldsymbol{x} )^\top \sigma \right) \, \mathrm{d} \boldsymbol{x}
 = \int_{\mathcal{B}} \rho \, (\delta \boldsymbol{x}) ^\top \boldsymbol{g} \,  \mathrm{d} \boldsymbol{x}  +  \int_{\partial \mathcal{B}} \boldsymbol{x} ^\top \boldsymbol{t} \; \mathrm{d} \boldsymbol{x}.
\end{equation}
The domain $\mathcal{B}$ is discretized in space in $n^e$ elements

\begin{equation}
\label{eq:discr}
\mathcal{B} \longrightarrow \bigcup_{i=1}^{n^e} \mathcal{B}_i.
\end{equation}
Now, the continuous displacement field can be approximated element-wise as

\begin{equation}
 \label{eq:shapefun}
\boldsymbol{x} \approx \sum_{k = 1}^{n} \phi_k (\xi, \eta, \zeta) \mathbf{x}_k = H \mathbf{x}^e,
\end{equation}
where 
\begin{equation}
H = \begin{bmatrix}
\phi_1 & 0 & 0 & \cdots & \phi_{n} & 0 & 0\\
0 & \phi_1 & 0 & \cdots & 0 & \phi_{n} & 0 \\
0 & 0 & \phi_1 & \cdots & 0 & 0 & \phi_{n} \\
\end{bmatrix}.
\end{equation}
Here, $\phi_k $ are shape functions of an element with $n$ nodes, which depend on the so-called isoparameter local coordinates within an element $(\xi, \eta, \zeta)$. The element displacement vector is assembled from the displacement vectors at each node:
\begin{equation*}
 \mathbf{x}^e = \begin{pmatrix}
                 \mathbf{x}_1 \\
                 \mathbf{x}_2 \\
                 \vdots \\
                 \mathbf{x}_n
                \end{pmatrix}.
\end{equation*}
The \textit{global displacement vector} $\mathbf{x}$ is related to the element displacement vector by the location matrix $Z_e$, where the topology of the discretization is stored, i.e.,
\begin{equation}
 \label{eq:Ze}
 \mathbf{x}^e = Z^e \mathbf{x}.
\end{equation}
To replace  the action of the $\nabla \cdot$ operation, the additional auxiliary matrix $L$ of size ${6 \times 3}$ is defined
\begin{equation} \label{mat:L}
L^\top = \begin{bmatrix}
\displaystyle\frac{\partial}{\partial \mathrm{x}_1} & 0 & 0 & \displaystyle \frac{\partial}{\partial \mathrm{x}_2} & 0 & \displaystyle \frac{\partial}{\partial \mathrm{x}_3}\\
0 & \displaystyle \frac{\partial}{\partial \mathrm{x}_2} & 0 &\displaystyle \frac{\partial}{\partial \mathrm{x}_1}& \displaystyle \frac{\partial}{\partial \mathrm{x}_3}& 0 \\
0 & 0 & \displaystyle \frac{\partial}{\partial \mathrm{x}_3} & 0 & \displaystyle \frac{\partial}{\partial \mathrm{x}_2} & \displaystyle \frac{\partial}{\partial \mathrm{x}_1}
\end{bmatrix}.
\end{equation}
We denote the element domain $\mathcal{B}_e$ with the boundary $\partial \mathcal{B}_e$ for each of the $n^e$ elements. With \eqref{eq:shapefun} and \eqref{eq:Ze}, the weak form of the balance of momentum \eqref{eq:weak} can be reformulated as

\begin{equation}
\begin{aligned}
 \label{virtfem}
 &\sum_{e = 1}^{n^e} \int_{\mathcal{B}_e} \rho (H Z^e \delta \mathbf{x})^\top H Z^e \ddot{\mathbf{x}} \, \mathrm{d} \mathbf{x}_e
 + \sum_{e = 1}^{n^e} \int_{\mathcal{B}_e} (L H Z^e \delta \mathbf{x})^\top \sigma \, \mathrm{d} \mathbf{x}_e \\ 
 & \hspace{4cm} =
\sum_{e = 1}^{n^e} \int_{\mathcal{B}_e} \rho (H Z^e \delta \mathbf{x})^\top \boldsymbol{g} \, \mathrm{d} \mathbf{x}_e +
 \sum_{e = 1}^{n^e} \int_{\partial \mathcal{B}_e} (H Z^e \delta \mathbf{x})^\top \boldsymbol{t} \, \mathrm{d}  \mathbf{x}_e.
\end{aligned}
\end{equation}
The equation \eqref{virtfem} must hold for any virtual displacement, leading to the following ODE system:
\begin{equation}
\label{eq:ode1}
 M \ddot{\mathbf{x}} + \mathbf{f}_{\text{int}} = \mathbf{f}_{\text{ext}},
\end{equation}
where
\begin{align}
	\label{eq:mass}
 M = \sum_{e = 1}^{n^e} \left( Z^e \right) ^\top & \left( \int_{\mathcal{B}_e} \rho (H^\top  H \,\mathrm{d} \mathbf{x}_e \right) Z^e, \quad \quad  \mathbf{f}_{\text{int}} = \sum_{e = 1}^{n^e} \left( Z^e \right) ^\top \int_{\mathcal{B}_e} (L H)^\top \sigma \, \mathrm{d} \mathbf{x}_e , \quad  \text{and}
 \end{align}
\begin{equation*}
  \mathbf{f}_{\text{ext}} = \sum_{e = 1}^{n^e} \left( Z^e \right) ^\top  \int_{\mathcal{B}_e} \rho H^\top \boldsymbol{g} \, \mathrm{d} \mathbf{x}_e  +
 \sum_{e = 1}^{n^e} \left( Z^e \right) ^\top \int_{\partial \mathcal{B}_e} H^\top \boldsymbol{t} \, \mathrm{d} \mathbf{x}_e  \hspace{2.7cm}
 \label{massm}
\end{equation*}
%We denote matrix $Q = L H $. 
% Material
are the consistent mass matrix, the vector of internal forces, and the vector of external forces, respectively. Since the material model has not yet been defined, the equation \eqref{eq:pdegov} is valid for linear and nonlinear material behavior and arbitrarily large displacement gradients.
To complete the field equations, we add the so-called constitutive equations.
In many scenarios, the stress tensor can be written as a linear function of the displacements.
For example, in the case of the classical Hooke's law
\begin{equation}
 \sigma = D^{el} \varepsilon,
 \label{eq:elastic}
\end{equation}
where $D^{el}$ is a fourth-order elastic stiffness tensor. Using the Voigt notation,
\begin{equation}
D^{el} = \begin{bmatrix}
					\lambda + 2 \mu & \lambda	      & \lambda         & 0 & 0 & 0 \\
					\lambda			& \lambda + 2 \mu & \lambda         & 0 & 0 & 0 \\
					 \lambda	    & \lambda		  & \lambda + 2 \mu & 0   & 0 & 0 \\
					  0 		    &			  0   &	 0              & \mu & 0 & 0 \\
					   0 & 0 & 0 & 0 & \mu & 0 \\
					    0 & 0 & 0 & 0 & 0 & \mu
					  \end{bmatrix},
\end{equation}
where $\lambda$ and $\mu$ are the Lam\'e constants. In its turn, the deformation field is approximated using \eqref{mat:L}, \eqref{eq:shapefun}, and \eqref{eq:kinem} as follows:
\begin{equation}
\varepsilon \approx L H \mathbf{x} = Q \mathbf{x}.
\end{equation}
Thus, the internal force vector can be written as
\begin{equation}
  \mathbf{f}_{\text{int}} = \sum_{e = 1}^{n^e} \left( Z^e \right) ^\top \int_{\mathcal{B}_e} Q^\top D^{el} Q \mathbf{x} \mathrm{d} \mathcal{B}_e.
\end{equation}
The equation \eqref{eq:ode1} takes the form
\begin{equation}
\label{eq:ode}
M\ddot{\mathbf{x}}(t) +  K \mathbf{x}(t) = \mathbf{f}_{\mathrm{ext}}(t),
\end{equation}
where the stiffness matrix is
\begin{equation}
  K = \sum_{e = 1}^{n^e} \left( Z^e \right) ^\top \int_{\mathcal{B}_e} Q^\top D^{el} Q \mathrm{d} \mathcal{B}_e.
\end{equation}
The computation of the system matrices requires numerical integration over the element domain using an appropriate method (e.g., Gauss integration). Of course, the dissipation forces also play an important role and is hence important to be taken into account in the internal force vector. It is described with a damping matrix $E$ analogously to the elastic forces and stiffness matrix. Very common in engineering practice is the Rayleigh damping model, which allows representing the damping matrix as a linear combination of mass and stiffness matrix, where the factors $\alpha_R$ and $\beta_R$ damp the lower and higher frequencies, respectively:
\begin{equation}
\label{eq:ray}
 E = \alpha_R M + \beta_R K.
\end{equation}
However, there are other damping models that can be preferably for different cases; in this work we are not limited to any particular model. Thus, in a general case, we have the following system of ODEs:
\begin{equation}
\label{eq:original}
M\ddot{\mathbf{x}}(t) + E\dot{\mathbf{x}}(t) + K\mathbf{x}(t) = \mathbf{f}(t),
\end{equation}
where $\mathbf{x}(t)$ are the fundamental unknowns --- nodal displacements; $M,K,E \in \mathbb{R}^{n \times n}$ are the system mass, stiffness, and damping matrices,  respectively.  The external force vector $\mathbf{f}(t)$ can be formulated for some applications in terms of a certain control operator $B \in \mathbb{R}^{n \times m}$ and input vector $\mathbf{u}(t) \in \mathbb{R}^{m}$, consisting of $m$ input signals $u(t)$
\begin{equation}
\mathbf{f}(t) = B \mathbf{u}(t).
\end{equation}
 It is worth mentioning that $M$ and $K$ are typically symmetric positive definite, and $E$ is symmetric positive semidefinite. We will denote this conditions as $M \succ 0,  K \succ 0,  E \succeq 0$. Moreover, if those conditions hold, it is well known that the mechanical system is stable, see \cite{Mue72,Tho21,morSalEL06}.

Equations \eqref{eq:original} describe the dynamics of a system and are often inaccessible from the FEM software. The system dimension is usually very high, which is natural, considering the high number of elements and nodes needed to maintain structure geometry precisely. Each system matrix depends on the material parameters, element type, and other specific FEM settings. Our primary goal in this work is to identify smaller dimension surrogate models having the mechanical structure as in \eqref{eq:original} using simulated data information, which is described in the next section.

\begin{remark}
In the presence of geometric nonlinearities, the deformation gradient tensor, which describes the rotation and deformation of the body, is no longer equal to the identity tensor due to the loss of equivalence between the deformed and undeformed configuration. Therefore, other appropriate stress and strain measures should be used to describe the motion of the system. In particular, it is natural for solid mechanics to use the original reference configuration, namely the Green-Lagrange strain tensor and the Second Piola-Kirchhoff stress tensor. As a consequence, the governing equation becomes nonlinear. Another source of nonlinearity can be material behavior, introducing a nonlinear relationship between stress and strain tensor. For these cases, the solution of the governing system of equations has to be computed in an iterative manner. Given that, the reduction method has to be performed in a more involved way, which will be considered in our future work.
\end{remark}

\section{Data setup}
\label{sec:setup}

In this section, we explore the available data, including a time-integration solver description. We assume that the model of a mechanical system \eqref{eq:original} is given as a \emph{gray-box}, i.e., the underlying abstract model structure is known by utilizing the physical knowledge laid out in the previous section, but the system operators are unavailable. Instead, we have access to the simulation input and output data, which consist of the excitation signals $\mathbf{u}(t)$ and the nodal displacements in the state vector $\mathbf{x}(t)$.
The simulation is performed with the following time discretization $0 =  t_0 < t_1 < \dots < t_N = T$  of the time domain $[0, T ]$.
 Further, we assemble the snapshot matrix $X$ and the input signal matrix $U$ by collecting the inputs and the snapshots of the state at pre-defined time-steps:
\begin{equation}
\label{mat:snap}
U =
\begin{bmatrix}
| &  \dots & | \\
\mathbf{u}(t_1) & \dots & \mathbf{u}(t_N) \\
| & \dots & |
\end{bmatrix} \in \mathbb{R}^{m \times N}, \quad
 X = \begin{bmatrix}
| &  \dots & | \\
\mathbf{x}(t_1) & \dots & \mathbf{x}(t_N) \\
| & \dots & |
\end{bmatrix} \in \mathbb{R}^{n \times N}.
\end{equation}
The time-integration in FE-packages is usually performed by second-order integration methods, such as Newmark-$\beta$ \cite{New59}, Hilber-Hughes-Taylor (HHT) method \cite{HilHT77}, and Generalized-$\alpha$ method \cite{ChuH93}. The latter two methods are the generalizations of the Newmark method with controllable numerical damping, which is particularly important for the automatic time stepping scheme to reduce the effect of the high-frequency noise resulting from too large step size or a poor spatial discretization.
Using the HHT method, the equilibrium \eqref{eq:original} is replaced by the following discretized expression:

\begin{subequations}
	\label{eq:hht}
\begin{align}
	M &\ddot{\mathbf{x}}_{k+1} + E ( (1 + \alpha) \dot{\mathbf{x}}_{k+1} - \alpha \dot{\mathbf{x}}_k ) + K ( (1 + \alpha) \mathbf{x}_{k+1} - \alpha \mathbf{x}_k ) = \mathbf{f}_{k+1}, \\
	&\mathbf{x}_{k+1} = \mathbf{x}_k + \Delta t \dot{\mathbf{x}}_k + (\Delta t)^2 \left[ \left( \frac{1}{2} - \beta \right) \ddot{\mathbf{x}}_k + \beta \ddot{\mathbf{x}}_{k+1} \right] , \\
	&\dot{\mathbf{x}}_{k+1} = \dot{\mathbf{x}}_k + \Delta t \left[ (1 - \gamma) \ddot{\mathbf{x}}_k + \gamma \ddot{\mathbf{x}}_{k+1} \right].
\end{align}
\end{subequations}
Numerical damping is controlled by the parameter $\alpha \in \left[ - \frac{1}{3}, 0 \right]$ (negative $\alpha$-dissipation). The parameters $\gamma$ and $\beta$ govern the stability of the algorithm and are often chosen as $\gamma = \frac{1 - 2 \alpha}{2}, \; \beta = \frac{(1 - \alpha)^2}{4} $~\cite{GeR14}. Setting $\alpha = 0$ makes \eqref{eq:hht} equivalent to the Newmark-$\beta$ family of algorithms, which we use in our numerical simulations. Hence, the derivative data needed for the system identification can also be extracted from the integrator, which is assembled as follows:
\begin{equation}
\begin{aligned}
 \displaystyle \dot{X} &= \begin{bmatrix}
| &  \dots & | \\
\dot{\mathbf{x}}(t_1) & \dots & \dot{\mathbf{x}}(t_N) \\
| & \dots & |
\end{bmatrix} \in \mathbb{R}^{n \times N}, \quad
\displaystyle \ddot{X} &= \begin{bmatrix}
| &  \dots & | \\
\ddot{\mathbf{x}}(t_1) & \dots & \ddot{\mathbf{x}}(t_N) \\
| & \dots & |
\end{bmatrix} \in \mathbb{R}^{n \times N},
\end{aligned}
\end{equation}
where $\dot{X}$ and $\ddot{X}$ contain velocities and accelerations information. 
%The quality of the derivative data may play an important role for the optimization problem, since it can provide additional errors to the optimization. \todo[inline]{It is not clear which optimization problem are you talking about. Therefore, I would suggest you either remove or rephrase without mentioning optimization problem. Like `` the quality of these information plays an important role in learning dynamical systems, as we shall see shortly``} 
Since solvers can often provide the velocity and acceleration data, we will use these data in our work.  

With this, we aim to develop a data-driven framework to learn second-order dynamical systems to capture the dynamics present in the data. Particularly,  our focus lies in constructing low-dimensional dynamical models to achieve our goal.

%%%%%%%%%%%%%%%%%%%%%%%%%%%%%%%%%
%%%%%%%%%%%%%%%%%%%%%%%%%%%%%%%%%
\section{Intrusive POD reduction}
\label{sec:pod}

Before proceeding to the description of a non-intrusive operator inference approach, we briefly recapitulate the intrusive snapshot-based POD method that forms the basis for identifying low-dimensional subspaces for data or the compression step for the operator inference method.
The main feature of the POD method is to identify orthogonal modes that optimally capture the energy present in the snapshot matrix. These modes also capture most of the dynamics in the data. This can be achieved by employing the singular value decomposition (SVD) of the snapshot matrix $X$ \eqref{mat:snap}:
\begin{equation}
\label{eq:svd}
X = V \Sigma W^\top.
\end{equation}
 Recall that according to the Eckart-Schmidt-Young-Mirsky theorem, the truncated SVD provides the best rank-$r$ approximation of a given matrix in the Frobenius norm \cite{morVol12}. In order to get the low-dimensional representation of the system dynamics, we approximate \eqref{eq:svd} by truncating the small singular values. Hence, we construct the subspace basis $V_r$ by choosing the first $r$ dominant left singular vectors.
The system operators in \eqref{eq:original} can be projected onto the subspace $V_r$, yielding the following reduced POD system:
\begin{equation}
\label{eq:pod}
\widetilde{M} \ddot{\tilde{\mathbf{x}}} (t) + \widetilde{E} \dot{\tilde{\mathbf{x}}} (t) + \widetilde{K} \tilde{\mathbf{x}} (t) = \widetilde{B} \mathbf{u}(t),
\end{equation}
with the reduced system operators being defined as 
\begin{equation}\label{eq:PODopt}
\widetilde{M} = V_r^\top M V_r, \quad 
\widetilde{E} = V_r^\top E V_r, \quad
\widetilde{K} = V_r^\top K V_r, \quad
\widetilde{B} = V_r^\top B.
\end{equation}
%The columns of $\widetilde{X}$ contain the solution of simulation of the  reduced POD model \eqref{eq:pod}, which will be referred to later in the paper. 
Notice that if the original matrices $M, E$, and, $K$ are symmetric positive (semi)definite, then so are the reduced matrices $\widetilde{M}$, $\widetilde{E}$ and  $\widetilde{K}$. As a consequence, the intrusive POD model preserves the stability of the original one, as mentioned in \cite{morSalEL06} in the context of moment matching.  Except for the basis construction from snapshots, the reduction is performed intrusively, i.e., it requires the original matrices $M$, $E$, $K$, and $B$, which describe the dynamics of the original mechanical systems. 
%In the next section, we propose two different non-intrusive model reduction methods to obtain reduced mechanical systems.%

%%%%%%%%%%%%%%%%%%%%%%%%%%%%%%%%%%%%%%%%
%%%%%%%%%%%%%%%%%%%%%%%%%%%%%%%%%%%%
\section{Operator inference for mechanical systems}
\label{sec:opinf}

Instead of projecting the known system operators, our goal is to infer the reduced operators using the data available in \Cref{sec:setup}.
Towards learning low-dimensional systems from given high-dimensional data, we first need to prepare an appropriate low-dimensional data representation. To that end, we aim at finding a low-dimensional approximation of the snapshot matrix \eqref{mat:snap}, which is done as described in \Cref{sec:pod} by applying SVD and choosing the $r$ most dominant singular vectors as a projection basis. Using the obtained dominant subspace, we prepare the compressed low-dimensional data as follows:
\begin{equation}\label{eq:low_dimensionaldata}
\begin{aligned}
\widehat{X} = V_r^\top X, \qquad
\dot{\widehat{X}} = V_r^\top \dot{X}, \qquad
\ddot{\widehat{X}} = V_r^\top \ddot{X},
\end{aligned} 
\end{equation}
assuming we have access to the velocity and acceleration vectors as well.
Next,  we present an optimization-based formulation to infer reduced-order operators directly using the data~\eqref{eq:low_dimensionaldata}.

\subsection{Second-order formulation}
\label{sub:opinf_u}
First, we recall that the intrusive POD model \eqref{eq:pod} is represented by the matrices $\widetilde{M}$, $ \widetilde{E}$,   $\widetilde{K}$, and $\widetilde{B}$. These reduced matrices satisfy the following equation:
\begin{equation}
\label{eq:pod_matrix}
\widetilde{M} \ddot{\widetilde{X}} + \widetilde{E} \dot{\widetilde{X}} + \widetilde{K} \widetilde{X} = \widetilde{B} U,
\end{equation}
where $\widetilde{X}$, $\dot{\widetilde{X}}$ and  $\ddot{\widetilde{X}}$, respectively, are the snapshot matrix assembling $N$ snapshots of the reduced POD model \eqref{eq:pod}, its corresponding derivative, and second-order derivative matrices, i.e.,
\begin{equation} \label{mat:pod_snap}
\tilde X = \begin{bmatrix}
| &   & | \\
\tilde{\mathbf{x}}(t_1) & \dots & \tilde{\mathbf{x}}(t_N) \\
| &  & |
\end{bmatrix}, \; \dot{\tilde{X}}= \begin{bmatrix}
| &   & | \\
\dot{\tilde{\mathbf{x}}}(t_1) & \dots & \dot{\tilde{\mathbf{x}}}(t_N) \\
| &  & |
\end{bmatrix}, \; \ddot{\tilde{X}}= \begin{bmatrix}
| &   & | \\
\ddot{\tilde{\mathbf{x}}}(t_1) & \dots & \ddot{\tilde{\mathbf{x}}}(t_N) \\
| &  & |
\end{bmatrix}.
\end{equation}
Assuming the reduced mass matrix $\widetilde{M}$ is invertible, we multiply \eqref{eq:pod_matrix} by $\widetilde{M}^{-1}$ from the left, yielding the following differential system of equations
\begin{equation}
\label{eq:pod_inv}
\ddot{\widetilde{X}} = - \widetilde{M}^{-1} \widetilde{E} \dot{\widetilde{X}} - \widetilde{M}^{-1} \widetilde{K} \widetilde{X} + \widetilde{M}^{-1} \widetilde{B} U.
\end{equation}
Hence, the dynamics of the POD intrusive model is fully described by the matrices $\widetilde{M}^{-1} \widetilde{E}$,  $\widetilde{M}^{-1} \widetilde{K}$ and $\widetilde{M}^{-1} \widetilde{B}$. It is important to notice that $\widetilde{M}^{-1} \widetilde{E}$ and  $\widetilde{M}^{-1} \widetilde{K}$  may not be symmetric positive (semi)definite, even if $\widetilde{M}$,  $\widetilde{E}$, $\widetilde{K}$ are.
%The need for this transformation is argued at the end of the current section. 
The structure \eqref{eq:pod_inv} is used as a foundation to formulate a least-squares problem using the projected data \eqref{eq:low_dimensionaldata}.
Inspired by the structure \eqref{eq:pod_inv}, our next goal is to identify a second-order reduced model of the form as follows:
\begin{equation}\label{eq:inf_ROM}
\ddot{\hat{\mathbf{x}}} (t) + \widehat{E}_{\mathrm{M}} \dot{\hat{\mathbf{x}}} (t) +\widehat{K}_{\mathrm{M}} \hat{\mathbf{x}} (t) =  \widehat{B}_{\mathrm{M}} \mathbf{u}(t),
\end{equation}
using the projected data $\widehat{X}$, $\dot{\widehat{X}}$ and $\ddot{\widehat{X}}$ in \eqref{eq:low_dimensionaldata} and the input data $U$. In particular, we seek to determine the matrices or operators $\widehat{E}_{\mathrm{M}}$, $\widehat{K}_{\mathrm{M}}$, and $ \widehat{B}_{\mathrm{M}}$.
Hence,  we propose the following second-order inference problem:
% \begin{equation}\label{eq:opinf}
%   \underset{\substack{\widehat{E}_{\mathrm{M}}, \widehat{K}_{\mathrm{M}} \in \mathbb{R}^{r \times r}, \\ \phantom{\widehat{E}_{\mathrm{M},}} \widehat{B}_{\mathrm{M}} \in \mathbb{R}^{r \times m}}}{\text{minimize}} \left\| \ddot{\widehat{X}} +  \widehat{E}_{\mathrm{M}} \dot{\widehat{X}} +  \widehat{K}_{\mathrm{M}} \widehat{X} -  \widehat{B}_{\mathrm{M}} U \right\| ^2_F,
% \end{equation}
\begin{equation}\label{eq:opinf}
  \underset{\substack{\widehat{E}_{\mathrm{M}}, \widehat{K}_{\mathrm{M}}, \widehat{B}_{\mathrm{M}} }}{\text{minimize}} \left\| \ddot{\widehat{X}} +  \widehat{E}_{\mathrm{M}} \dot{\widehat{X}} +  \widehat{K}_{\mathrm{M}} \widehat{X} -  \widehat{B}_{\mathrm{M}} U \right\| ^2_F,
\end{equation}
where the matrices $\displaystyle \widehat{E}_{\mathrm{M}}, \, \widehat{K}_{\mathrm{M}} \in \mathbb{R}^{r \times r}$, and $\displaystyle \widehat{B}_{\mathrm{M}} \in \mathbb{R}^{r \times m}$ are the unknown operators. Since the intrusive matrices in \eqref{eq:pod_inv} $\widetilde{M}^{-1} \widetilde{E}$ and $\widetilde{M}^{-1} \widetilde{K}$ may not be symmetric positive definite, we expect the same for the inferred matrices $\widehat{E}_{\mathrm{M}}$ and $\widehat{K}_{\mathrm{M}}$.
In order to reformulate the optimization problem \eqref{eq:opinf} in a more compact way, we assemble the global data matrix:
\begin{equation} \label{mat:data_opinf}
\widehat{\cD} = \begin{bmatrix} \dot{\widehat{X}} {\,}^{\top},~ \widehat{X}^\top ,~ U^\top \end{bmatrix}^{\top}
\end{equation}
using the available project snapshot matrices, except for the second-order derivative matrix $\ddot{\skew4\widehat{X}}$, which plays the role of the right-hand side for the regression problem.
Finally, we state the optimization problem as follows:
	\begin{equation}
	\label{eq:opinf_compact}
	\underset{\widehat{P} \in \mathbb{R}^{r \times (2r+m)}}{\text{minimize}} \left\|  \widehat{P}\widehat{\cD}  - \ddot{\widehat{X}} \right\| ^2_F,
	\end{equation}
where the variable parameter matrix consists of all the unknown operators 
\begin{equation}\label{mat:inferred_operators}
\widehat{P} = \begin{bmatrix}
- \skew4\widehat{E}_M ,~   -\skew4\widehat{K}_M ,~ \skew4\widehat{B}_M
\end{bmatrix}.
\end{equation}
It is worth mentioning that the inferred model is obtained non-intuitively, i.e., the construction of the matrices $\widehat{E}_{\mathrm{M}}$, $\widehat{K}_{\mathrm{M}}$ and $ \widehat{B}_{\mathrm{M}}$ is based only on the provided data. Also, in this setup, the mass matrix of the inferred model \eqref{eq:inf_ROM} is assumed to be the identity by construction. 

\begin{remark}
One may argue that the mass matrix can also be identified using this approach. To this aim, one needs to include the mass matrix in the unknown operators
\[ \widehat{P}_{mod} = \begin{bmatrix}
- \skew4\widehat{M}	&- \skew4\widehat{E} &   -\skew4\widehat{K} & \skew4\widehat{B}
\end{bmatrix}, \]
and add the projected second derivative to the data matrix as follows
\begin{equation*}
	\widehat{\cD}_{mod}= \begin{bmatrix} \ddot{\widehat{X}}{\,}^\top& \dot{\widehat{X}} {\,}^{\top} & \widehat{X}^\top & U^\top \end{bmatrix}^{\top}.
\end{equation*}
Hence, to infer the reduced operators, one would have to solve the following least square problem 
	\begin{equation*}
	\underset{\widehat{P}_{mod} \in \mathbb{R}^{r \times (2r+m)}}{\text{minimize}} \left\| \widehat{P}_{mod}\widehat{\cD}_{mod}  \right\| ^2_F,
\end{equation*}
It consists of a least-squares problem without a right-hand side, for which zero is a trivial solution. Problems of this type are usually solved by constraining the size of the solution norm, which is not suitable for our case. In \Cref{subsec:ForceOpInf}, we will propose an approach enabling us to also infer the mass matrix, provided that some additional data is available.
\end{remark}

%Now, the necessity of apply the inverse of the mass matrix in \eqref{eq:pod_inv} can be justified by the fact that using the POD formulation \eqref{eq:pod} for developing the operator inference procedure, without making any assumptions on the mass matrix, leads to the least-squares problem without right-hand side. Problems of this type are usually solved by constraining the size of the solution norm, which is not suitable for our case. Hence, the inferred second-order operators \eqref{mat:inferred_operators} are the solution of the optimization problem \eqref{eq:opinf_compact}, formulated inspired by the original \eqref{eq:original} and intrusive POD system structure \eqref{eq:pod_inv}.

%Operator inference and POD method use the same basis for the low-dimensional subspace, used for the data compression in the non-intrusive case, and on the other hand for the projection of system matrices in the intrusive case. Studying closeness of these two models is a logical follow-up in developing the operator inference method. 

\subsubsection{Theoretical closeness of the intrusive and non-intrusive ROMs}
%%
% In this part, we establish the connection between the inference reduced model and POD reduced model.
Although the inferred and intrusive reduced operators are obtained with different procedures, we can show  an asymptotic closeness of these two models. The original paper on operator inference \cite{morPehW16} and some other articles, such as \cite{morBenGHetal22}, provide theoretical results for the first order reduced systems. They show, under  certain assumptions, that the inferred matrices are an approximation of the intrusive reduced matrices in the Frobenius norm. When such a result holds, the inferred system can inherit several useful properties of POD models, such as stability and error analysis.

Let the parametric matrix $\widehat{P}$ be the solution of the optimization problem \eqref{eq:opinf} with the corresponding matrix $\widehat{\cD}$, constructed from the available data. We denote $\mathbf{x}(t_i)$ as the continuous displacement at the time $t_i$, and $\mathbf{x}_i$ as the discretized displacement snapshot-vector. Further, we consider the following assumptions.

\begin{assumption}
Time-stepping scheme is convergent, i.e., $\displaystyle \lVert \mathbf{x}_i - \mathbf{x}(t_i) \rVert \rightarrow 0$ as $\Delta t \rightarrow 0$. %where $\mathbf{x}_i$ is the i-th column of the snapshot matrix $X$.	
\end{assumption}
\begin{assumption}
The discretized reduced derivative data converges to the continuous derivative data, i.e., $ \displaystyle \lVert \dot{\hat{\mathbf{x}}}_i - \frac{\mathrm{d}}{\mathrm{d}t}\hat{\mathbf{x}}(t_i) \rVert \rightarrow 0$ and $\displaystyle \lVert \ddot{\hat{\mathbf{x}}}_i - \frac{\mathrm{d}^2}{\mathrm{d}t^2}\hat{\mathbf{x}}(t_i) \rVert \rightarrow 0$ as $\Delta t \rightarrow 0$.%, where $\mathbf{x}_i$ is the i-th column of the snapshot matrix $\dot{X}$.
\end{assumption}
\begin{assumption}
The matrix $\widehat{\cD} \in \mathbb{R}^{N \times(2r+m)}$ has full rank, assuming that the dimension $r$ is much smaller than the number of time steps $N$.
\end{assumption}%
Using the above assumptions, we formulate the following theorem:
\begin{theorem}
	\label{theorem}
Let Assumptions 1,2,3 hold and $\widetilde{M},  \widetilde{E}, \widetilde{K},$  and $\widetilde{B}$ be the reduced-order operators obtained intrusively as in \eqref{eq:PODopt} using the POD basis $V_r$. Then, for every $ \varepsilon > 0 $ there exist a reduced order $r<n$ and a step size $ \Delta t > 0$ such that %the inferred operators are close to the operators, obtained intrusively by POD, i.e. 
\[ \lVert \widetilde{M}^{-1} \widetilde{E}  - \skew4\widehat{E}_M \rVert _F < \varepsilon, \quad  \lVert \widetilde{M}^{-1} \widetilde{K}  -  \skew4\widehat{K}_M \rVert _F < \varepsilon, \quad  \text{and} \quad  \lVert \widetilde{M}^{-1} \widetilde{B}  -  \skew4\widehat{B}_M \rVert _F < \varepsilon, \] 
where $\skew4\widehat{E}_M$, $\skew4\widehat{K}_M$ and $\skew4\widehat{B}_M$ are the inferred operators via the optimization problem \eqref{eq:opinf}.
\end{theorem}

\begin{proof}
%
%The aim is to show that the solution of the optimization problem \eqref{eq:opinf} tends to the POD-operators as the time-step tends to zero and the reduced order tends to the full order.
Recall that the intrusive POD reduced model has the form \eqref{eq:pod_inv}.
Let $\widetilde\cD = \begin{bmatrix}
	\dot{\widetilde{X}}{\,}^\top &
	\widetilde{X}^\top &
	U^\top
\end{bmatrix}^{\top}$ denote the corresponding data matrix for the system in \eqref{eq:pod_inv} with the POD snapshot matrices, defined in \eqref{mat:pod_snap}. Hence, the concatenated  intrusive reduced operators $\widetilde P = \begin{bmatrix}
-\widetilde{M}^{-1}\widetilde{E} & 	-\widetilde{M}^{-1}\widetilde{K} & \widetilde{M}^{-1}\widetilde{B}
\end{bmatrix}$ represent one solution of the least-squares problem
\begin{equation}\label{eq:intrusive_opinf}
\tilde P =	\text{arg}\vspace{-0.2cm}\min_P{ \left\|P \widetilde \cD - \ddot{\widetilde{X}} \right\| ^2_F}.
\end{equation} 
Moreover, it represents the unique solution if the matrix $\widetilde\cD$ has full rank. %\newline
Next, the projected matrix $\widehat\cD$ \eqref{mat:data_opinf} and the projected second order derivative  $\ddot{\widehat{X}}$ can be interpreted, respectively, as a disturbed POD data matrix $\widetilde{\cD}$ and disturbed second order POD derivative $\ddot{\widetilde{X}}$, i.e.,

\begin{equation}\label{eq:OpInf_pod_data}
\widehat{\cD} = \widetilde{\cD} + \delta \widetilde{\cD} \quad \text{and} \quad  \ddot{\widehat{X}} = \ddot{\widetilde{X}} +  \delta \ddot{\widetilde{X}}. 
\end{equation}
Indeed, the disturbing term $\delta \widetilde{\cD}$ comes from the time-sampling error of the solution data and from the approximation error considering $\displaystyle X \approx V_r \widetilde{X}$ and $\displaystyle X \approx V_r \widehat{X}$, which also holds for the first and second order derivative data.
Hence, $\delta \widetilde{\cD} \rightarrow 0$ and $\delta \ddot{\widetilde{X}} \rightarrow 0$ as $r \rightarrow n$ and $\Delta t \rightarrow 0$. Therefore,  this leads to the following asymptotic result for the least-squares problem
\begin{equation*}
\min_{\widehat{P}}  \left(  \lim_{\substack{\Delta t \to 0 \\ \phantom{\Delta}r \to n}} \left\| \, {\widehat{P}} \cdot \widehat{\cD} - \ddot{\widehat{X}} \, \right\| ^2_F \right)  =  \min_{\widetilde{P}} \left(  \lim_{\substack{\Delta t \to 0 \\ \phantom{\Delta}r \to n}} \left\| \, \widetilde{P} \cdot \left(\widetilde{\cD} +\delta \widetilde{\cD}\right)- \left(\ddot{\widetilde{X}} +  \delta \ddot{\widetilde{X}}\right) \, \right\| ^2_F \right) 
\end{equation*}
% =  \min_{\widetilde{P}} \left\| \, \widetilde{P} \cdot \widetilde{\cD}  - \ddot{\widetilde{X}} \, \right\| ^2_F   
In other words, the operator inference problem in \eqref{eq:opinf_compact} can be seen as a  perturbed version of the minimization problem in \eqref{eq:intrusive_opinf}. The pre-asymptotic case combined with the assumption that $\widehat{\cD}$ has full rank leads to the proof of the theorem.

\end{proof}
The above theorem states that if the least-squares problem is well-conditioned, then in the asymptotic case, when the time step converges to zero, and the reduced order converges to the full dimension, the operators obtained by POD are close to the inferred operators. This result is important because, for the broad class of mechanical systems, the POD method preserves stability by keeping symmetric positive definiteness of the system matrices due to one-sided projection. Therefore, the inferred model will also be stable in case it is close enough to the POD one.  However, the relevant properties can be inherited only for the asymptotic case. Moreover, in \cite{morPeh20},  for discrete-time linear first-order systems, it has been shown that it is possible to exactly recover the intrusive operators for any order using a re-projection scheme. %Such an analysis would certainly be interesting to extend to second-order mechanical systems, which will be studied in our future work.

%
%It is important to note that the accuracy of the operator inference model can be still even better than the accuracy of POD model as it is shown in some papers for particular cases \cite{morBenGKetal20}. 
%

%One of the assumptions requires the data matrix to have full rank. It means, the results hold only if the optimization problem is well-posed. 
In many applications, the data matrix is numerically rank deficient and the corresponding least-squares problem becomes ill-conditioned. Therefore, it is necessary to use  appropriate regularization techniques. 
Among different methods (such as truncated SVD or truncated QR) \cite{morYilGBetal20},  the Tikhonov regularization \cite{TikGetal77} is one of the widely used techniques. The optimization problem \eqref{eq:opinf} is replaced by the following regularized problem:
\begin{equation}\label{eq:reg}
	\widehat{P} = \text{arg}\vspace{-0.2cm}\min_P \left( \rVert P \cdot \widehat{\cD} - \ddot{\widehat{X}} \lVert ^2_F + \lambda \lVert P \rVert ^2_F \right),
\end{equation}
where $\lambda$ is a penalty parameter. The choice of $\lambda$ plays an important role in obtaining a good solution. One of the criteria is to ensure a minimal residual $\displaystyle \rVert \widehat{P} \cdot \widehat{\cD} - \ddot{\widehat{X}} \lVert_F$ for the smallest operator norm $\displaystyle \rVert \widehat{P} \lVert _F $.
In this work, we use the Tikhonov regularization by penalizing all the operators with the same regularization parameter. 
% \cite{morSawKP21,morMcqHW21}

\subsubsection{Separating the operators}
From the above theorem, we conclude that the inferred operators are close to the POD matrices in \eqref{eq:pod_inv}, assuming that the matrix $\widetilde{M}$ is "absorbed" in other operators. Further, we assume that the inferred operators can be decomposed as follows:
	\begin{equation}\label{eq:coupled}
	\widehat{E}_M = \widehat{M}^{-1} \widehat{E}, \quad
	\widehat{K}_M = \widehat{M}^{-1} \widehat{K}, \quad
	\widehat{B}_M = \hspace{1pt} \widehat{M}^{-1} \widehat{B}.
	\end{equation}
In order to obtain a ROM with second-order structure as in \eqref{eq:original}, we may think about some post-processing method to separate the inferred operators. The suggested procedure below uses the transformation of the operator inference model to the modal coordinates. The generalized eigenvalue problem can be written as

\begin{equation}
\widehat{K}_M \Phi = \Phi \Omega^2,
\end{equation}
where $\Omega$ is a diagonal matrix with the natural eigenfrequencies on the diagonal and $\Phi$ are the eigenmodes of the operator inference system.
The reduced stiffness matrix in modal realization is equal to $\Omega^2$, while the reduced modal mass matrix is identity. Using the fact that the modal stiffness is defined as
\begin{equation}
\Omega^2 = \Phi^\top \widehat{K} \Phi,
 \label{eq:kmod}
\end{equation}
 we can extract the reduced stiffness matrix from \eqref{eq:kmod}, using the eigenfrequencies and eigenmodes
\begin{equation}
 \widehat{K} =\Phi^{-\top} \Omega^2 \Phi^{-1}.
\end{equation}
Then, we can separate the reduced mass matrix, and the damping matrix from \eqref{eq:coupled} as
\begin{equation}
 \widehat{M} = \widehat{K} \widehat{K}_M ^{-1}, \quad  \widehat{E} = \widehat{M} \widehat{E}_M.
\end{equation}
We would like to stress that there are no guarantees that the separated operators will satisfy stability properties. A possible remedy to ensure the stability of the learned model is by performing post-processing by finding the nearest symmetric positive definite matrix as in \cite{Hig88}.
This can be done for mass, stiffness, and damping matrices if needed, but it would be at the expense of losing the accuracy of the learned models.

\subsection{Force-informed operator inference}\label{subsec:ForceOpInf}	
%\todo[inline]{In this section we are presenting another formulation of operator inference. We have to make clear what is the data available and also and an story line. 
%\newline 
%
%Briefly summarize  what was presented before. Say that we cannot guarantee stability and improve symmetry conditions on the operators, because of remark . Recall that spd conditions typically happen for mechanical system is this is enough for stability. That is why in this section we will propose another method . 
%\newline
%
%In contrast to what was presented before, here we will assume not only the data (recall) but also the full force data.
%\newline
%
%Emphasis the optimization problem we have in this method. Say the differences regarding the previous one.
%
% }
In the previous subsections, we have defined the second-order operator inference method for learning the reduced mechanical models of structure \eqref{eq:pod_inv}, using the state and derivative data, as explained in \Cref{sec:setup}. As discussed previously, we were not able to impose the symmetric positive definiteness of the inferred operators in this formulation, even if the intrusive reduced model possesses this structure. 

In this section, we present a alternative operator inference methodology, enabling us to enforce the system's matrices to be symmetric positive definite. To this aim, we will use additional information from the full-order model. Hence, in this section, we will assume we have access to all the external forces and their positions, meaning the vector $\mathbf{f}(t)$ in \eqref{eq:original} is given.  
%Then, the given external force-vector $\mathbf{f}(t)$ can be used in the optimization problem as a right-hand side. 

In many engineering applications, an analysis of a system response under a certain load is required. In these scenarios, the forces acting on the system are known and can be extracted from the simulation software (for example, using the input-file defining the simulation setup). Moreover, for some simulations, the load data may come from experimental measurements of real working conditions, given as force values at certain time-space points. The force matrix can be constructed from the force snapshots at the pre-defined time steps:
\begin{equation} \label{force_data}
F = \begin{bmatrix}
| &   & | \\
\mathbf{f}(t_1) & \dots & \mathbf{f}(t_N) \\
| &  & |
\end{bmatrix}. 
\end{equation}
%Similar to the projected trajectory data in \eqref{eq:low_dimensionaldata}, the force data can be therefore projected onto dominant POD subspace $\skew4\widehat{F} = V^\top F$. 
As a consequence, the POD reduced model satisfies the following equation:
\begin{equation}
	\label{eq:pod_matrix_withforce}
	\widetilde{M} \ddot{\widetilde{X}} + \widetilde{E} \dot{\widetilde{X}} + \widetilde{K} \widetilde{X} = V^\top F,
\end{equation}
where, once again,  $\widetilde{X}$, $\dot{\widetilde{X}}$ and  $\ddot{\widetilde{X}}$  are, respectively, the snapshot matrix assembling $N$ snapshots of the projected reduced POD model in \eqref{eq:pod}, its corresponding derivative and second-order derivative matrices as in \eqref{mat:pod_snap}. We also recall that the intrusive reduced operators $\widetilde{M}$,  $\widetilde{E}$ and  $\widetilde{K}$ are typically symmetric positive (semi)definite, implying that the intrusive model is stable. Hence, our goal in this section is to infer the second-order operators $\widehat{M}$, $\widehat{E}$ and $\widehat{K}$ of a reduced model of the form
\begin{equation}\label{eq:inf_ROM_const}
	 \widehat{M} \ddot{\hat{\mathbf{x}}} (t) + \widehat{E} \dot{\hat{\mathbf{x}}} (t) +\widehat{K}\hat{\mathbf{x}} (t) =  V^{\top}\mathbf{f}(t),
\end{equation}
such that
\begin{equation}\label{constraints}
	\widehat{M} \succ 0,\quad  \widehat{K} \succ 0, \quad \widehat{E} \succeq 0.
\end{equation}
For this, similar to the projected trajectory data in \eqref{eq:low_dimensionaldata}, the force data can be projected onto the dominant POD subspace $\skew4\widehat{F} = V^\top F$.  Moreover, let the new data matrix include the state and derivative data as follows:
\begin{equation}
\widehat{\cD} = \begin{bmatrix} \ddot{\widehat{X}} {\,}^{\top} &\dot{\widehat{X}} {\,}^{\top} & \widehat{X}^\top \end{bmatrix} ^\top.
\end{equation}
%
%Our goal is now to postulate the optimization problem imposing constraints related to the original system structure. Recall that very often, mechanical system operators have the following properties:
% 
%\begin{align}
%M \succ 0,\quad  K \succ 0, \quad E \succeq 0 .
%\label{constraints}
%\end{align}
%%This is true for structural dynamical systems with linear material law and modal damping.
%Given the s.p.d. mass matrix, which follows from \eqref{eq:mass}, the system \eqref{eq:original} with the properties \eqref{constraints} is guaranteed to be stable \cite{Mue72,Tho21}. For this reason we are aiming to construct reduced operators, which preserve the properties \eqref{constraints} by adding the constraints to the optimization problem. Note that in the general case the second-order formulation given in \Cref{sub:opinf_u} does not allow to add similar constraints, because the provided (coupled) operators are not symmetric.
%
The operator inference optimization problem with constraints \eqref{constraints} using the external force data \eqref{force_data} is formulated as follows:

\begin{equation}
	\label{eq:optim_c2}
	\underset{\substack{\skew4\widehat{M} \succ 0 ~\widehat{E} \succeq 0 ~ \widehat{K} \succ 0}}{\text{minimize}} \left\| ~ \left[\widehat{M} ~~ \widehat{E} ~~ \widehat{K}\right]  \widehat{\cD} - \widehat{F} \, \right\|^2_F.
\end{equation}
In practice, it is not possible to add rigid constraints to the optimization problem, therefore the reformulation  $ \widehat{M} - \omega I \succeq 0 $ and $ \widehat{K} - \omega I \succeq 0 $ with a small positive threshold $\omega > 0$ can be used to ensure the strict positive definiteness.
The operator inference formulation \eqref{eq:optim_c2} is a convex optimization problem, which can be solved using semidefinite programming algorithms, e.g., \cite{Boy94}. 
In contrast to the optimization problem \eqref{eq:opinf_compact}, which has an analytical solution via Moore-Penrose inverse, the problem \eqref{eq:optim_c2} requires linear matrix inequality solvers, which are computationally more expensive. However, since the computations are done in the POD-reduced dimension, they can still be performed in moderate time. Moreover, this methodology has the advantage of preserving the symmetric positive definite structure of the inferred system's operators, which implies that the inferred model is also stable.
%However, for the benchmarks presented in this work, it was not necessary, since even the positive semidefinite constraints have provided s.p.d. matrices.
%\begin{remark}
%In the following, we provide the results only for the displacement field. Since the relationships between displacement, deformation, and stress remain unknown from the proprietary code, the results can be extended only after these dependencies are approximated or identified.
%\end{remark}

\section{Numerical results}
\label{sec:num}
In this section, we study the performance of the proposed operator inference methodologies for mechanical systems to learn reduced-order models directly from data and present a comparison with the intrusive POD approach.
For this purpose, numerical experiments are performed, namely for international space station \cite{slicot_iss}, butterfly gyroscope \cite{morwiki_gyro} benchmarks, and vibrating plate model \cite{supAumW22}.
%\cite{morKorR05}.
 The first model is used for the analysis of vibrations caused by the docking of an incoming spaceship. The model is given in a first-order state-space realization, which originates from the second-order form, and can thus be transformed back to the second-order form. The second benchmark is a finite element structural model of a vibrating micro-mechanical gyroscope for inertial navigation applications. For a more detailed description of the model, we refer to \cite{morKorR05}. The latter example is a finite-element model for analysis of a vibration response of the aluminium plate exited by a point load.

The time integrator for the simulations of the full-order model and reduced-order models is described in \Cref{sec:setup}. The Newmark parameters in \eqref{eq:hht} are chosen as $\gamma = \frac{1}{2}$ and $\beta = \frac{1}{4}$, which are based on the average constant acceleration assumption ensuring the unconditional stability of the method.
For the implementation of the optimization with linear matrix inequality constraints, the YALMIP Toolbox \cite{Lof04} is used together with the SeDuMi solver\footnote{\url{https://sedumi.ie.lehigh.edu/}}.

The quality analysis of the ROMs is done by comparing the state trajectories and inspecting the relative state error, which is given by the relation to the maximum norm of the state vector, $\max \lVert \mathbf{x}(t) \rVert _2$. This is,
\begin{equation}
\epsilon_{\mathrm{err}} (t_i) = \frac{\left\| \mathbf{x}(t_i) - \hat{\mathbf{x}}(t_i)\right\|_2}{\max_{t\in [t_1, \; t_N]} \left\| \mathbf{x}(t) \right\|_2}.
\end{equation}
A comparison is performed for the original full-order model (\texttt{FOM}), the POD-reduced model (\texttt{POD}), the operator inference model in the second-order formulation (\texttt{OpInf}), and the force-informed operator inference model with constraints (\texttt{cOpInf}). All experiments were performed using \matlab~(2021a) running on an HP Probook 430 G3, 2.30 GHz \intel~\coreithree-6100U CPU, 8GB of RAM. \newline \newline
\textbf{Code availability} \newline
The source code of the implementations and the raw data are available at \newline
\url{https://gitlab.mpi-magdeburg.mpg.de/filanova/mechopinf}.

\subsection{International space station}
\begin{figure}[h!]
	\begin{center}
		\input{fig/tikz/iss_svd2.tex}
		\caption{Singular value decay for ISS model with the excitation signal $u(t) = \sin (t)$. The black dot denotes the singular value for the truncation order $r  = 4$. }
		\label{fig:iss_svd}
	\end{center}
\end{figure}
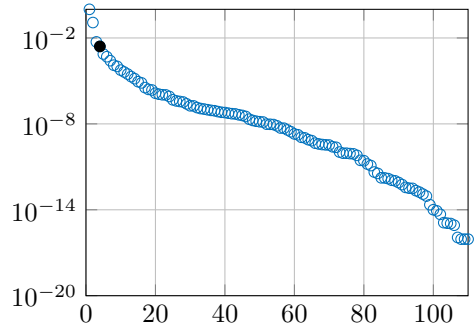

The structural model of the international space station \cite{morAntSG01} is a second-order system used for vibration analysis with the state dimension $n = 135$. The benchmark data are available in \cite{slicot_iss}. As a first step towards learning  intrusive POD and operator inference reduced models, we collect the training data in the time-interval $[ 0, 7s ]$  with the time step $\Delta t = 0.01s$ and input signal $u(t) = \sin (t)$. In \Cref{fig:iss_svd}, the normalized singular value decay is depicted for the collected snapshot matrix $X$, defined in \eqref{mat:snap}. The black dot denotes the singular value, corresponding to the order $r = 4$,  at which the truncation is done. The reduced order is selected so that the approximation error is at least below the threshold $10^{-2}$. 

\begin{figure}[h!]
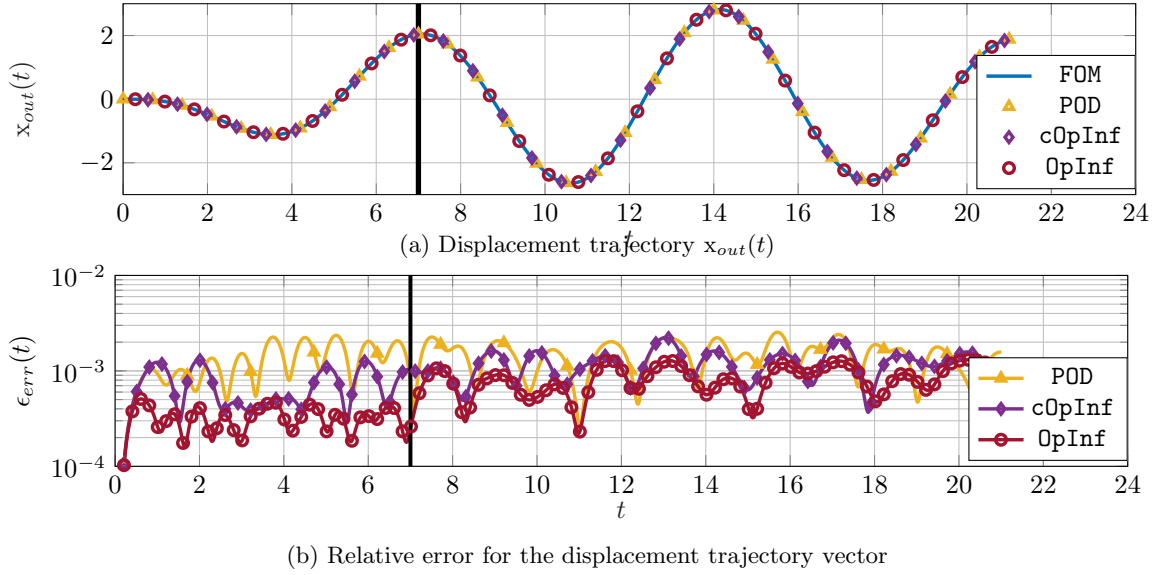

	\centering
\begin{subfigure}{\textwidth}
	\input{fig/tikz/iss_d_c.tex}
	\vspace{-15pt}
	\caption{Displacement trajectory $\mathrm{x}_{out} (t)$}
	\label{fig:iss_x_c}
	\vspace{-0.1cm}
\end{subfigure}
\vfill
\begin{subfigure}{\textwidth}
	\input{fig/tikz/iss_er_c.tex}
	\caption{Relative error for the displacement trajectory vector}
	\label{fig:iss_xe_c}
	%\vspace{-0.6cm}
\end{subfigure}
\caption{ISS benchmark: a comparison of FOM and ROMs of order $r = 4$. Black vertical line denotes the training period, used for constructing \texttt{POD, cOpInf}, and \texttt{OpInf} models.}
\end{figure}

The testing is performed on a time-interval $[ 0, 21s ]$ with the same time step and input as for the training phase.
In \Cref{fig:iss_x_c}, the trajectory for the second component of the displacement vector $\mathbf{x} (t)$ is shown. The curves show a good capture of the dynamics in the training phase and in the testing phase. \Cref{fig:iss_xe_c} shows it more clearly in the comparison of the relative error of the state trajectories. The operator inference reduction without constraints shows slightly better accuracy in the training phase, but in the testing phase, all these methods yield similar errors.

In general, both formulations of the operator inference methodology show good results for this example.

\subsection{Butterfly Gyroscope}
Our next example is the butterfly gyroscope \cite{morBil05} which is a linear second-order model with the state dimension $n = 17 132$.  The benchmark data are available in \cite{morwiki_gyro}. The model contains \texttt{s.p.d.} mass and stiffness matrices, the damping is modeled using the Rayleigh assumption -- a model with pure stiffness damping, where the coefficients are $\alpha_R = 0$ and $\beta_R = 10^{-6}$, see the equation \eqref{eq:ray}.  The training data is obtained by the simulation of the system on $t = [ 0  , \, 10^{-3} ]s$ with the time step $\Delta t = 10^{-6}s$ and input signal $u(t) = \sin (2 \pi f t)$ with $f =1 \; kHz$. In \Cref{fig:g_svd}, we depict the singular value decay of the collected snapshot matrix $X$, defined in \eqref{mat:snap}. The reduced order is selected as $r = 6$.

\begin{figure}[h!]
	\begin{center}
	\input{fig/tikz/g_svd.tex}
	\caption{Singular value decay for butterfly gyroscope model with the excitation signal $u(t) = \sin (2 \pi f t)$. The black dot denotes the singular value for the truncation order $r  = 6$.}
	\label{fig:g_svd}
	\end{center}
\end{figure}

\begin{figure}[h!]
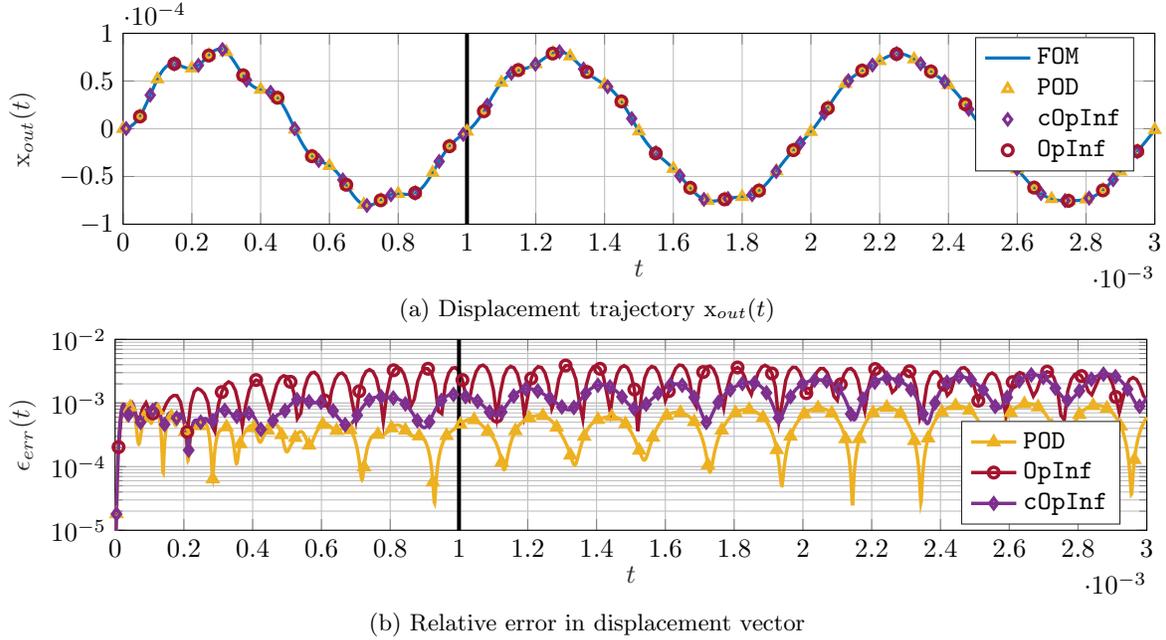

	\centering
	\begin{subfigure}{\textwidth}
	\input{fig/tikz/g_d_c.tex}
	\vspace{-15pt}
 	\caption{Displacement trajectory $\mathrm{x}_{out} (t)$}
	\label{fig:g_top}
	\vspace{-0.1cm}
	\end{subfigure}
	\vfill
	\begin{subfigure}{\textwidth}
	\input{fig/tikz/g_er_c.tex}
 	\caption{Relative error in displacement vector}
 	\label{fig:g_bot}
	\end{subfigure}
	\caption{Butterfly gyroscope benchmark: a comparison of FOM and ROMs of order $r = 6$. Black vertical line denotes the training period, used for constructing \texttt{POD, cOpInf}, and \texttt{OpInf} models.}
\end{figure}

The testing is performed for the same time-step and input load over a longer time interval $t = [ 0, \; 3 \cdot 10^{-3} ]s$.
The qualitative comparison of the trajectories for ROMs of reduced order $r = 6$ is presented in \Cref{fig:g_top} for the displacement component $\mathrm{x}_{out} = \mathrm{x}_{3143}$, which corresponds to one of the degrees of freedom, where the external force is applied. Over the whole simulation time, the ROMs are able to capture oscillations of the original system. To analyze the performance of the ROMs for the displacement field, we demonstrate the relative error for the state trajectories in \Cref{fig:g_bot}. As for the previous benchmark, the error does not exceed 1\%; therefore, we can ensure a good match of the state trajectories. For the whole simulation time, the force-informed formulation has slightly better accuracy than the operator inference formulation without constraints. However, the POD model shows a better performance than all non-intrusive ROMs. 

The deterioration in the accuracy of the operator inference model compared to the previous benchmark may be explained by a more ill-conditioned least-squares problem resulting from high-frequency loading and higher state dimension. %Although intrusive MOR provides better results for this example, it requires access to the full-order model and training data. In contrast, the non-intrusive operator inference formulation can predict the behavior with acceptable accuracy by constructing reduced-order models using only data.

\subsection{Vibrating plate}
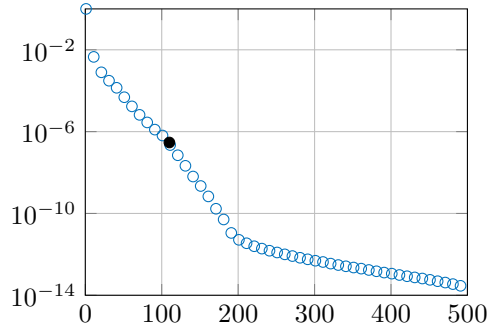
\begin{figure}[h!]
	\begin{center}
		\input{fig/tikz/p_svd.tex}
		\caption{Singular value decay for the vibrating plate model with the excitation signal $u(t) = \sin (2 \pi f t)$. The black dot denotes the singular value for the truncation order $r  = 110$.}
		\label{fig:p_svd}
	\end{center}
\end{figure}

Finally, we present the results for a model of a simply supported strutted plate excited by a point load \cite{morAumW22}. The model data is available from \cite{supAumW22}. This is a linear second-order model of state dimension $n = 201900$.  The damping is modeled using the Rayleigh assumption, where the coefficients are $\alpha_R = 0.01$ and $\beta_R = 10^{-4}$, see the equation \eqref{eq:ray}.  The training data is obtained by the simulation of the system on $t = [ 0  , \, 0.5 ]s$ with the time step $\Delta t = 10^{-3}s$ and input signal $u(t) = \sin (2 \pi f t)$ with $f =10 \; Hz$. In \Cref{fig:p_svd}, we depict the singular value decay of the collected snapshot matrix $X$, defined in \eqref{mat:snap}. To ensure the desired accuracy the reduced order is selected as $r = 110$.

\begin{figure}[h!]
	\centering
	\begin{subfigure}{\textwidth}
		\input{fig/tikz/p_x.tex}
		\vspace{-15pt}
		\caption{Displacement trajectory $\mathrm{x}_{out} (t)$}
		\label{fig:p_top}
	\end{subfigure}
	\vfill
	\begin{subfigure}{\textwidth}
		\input{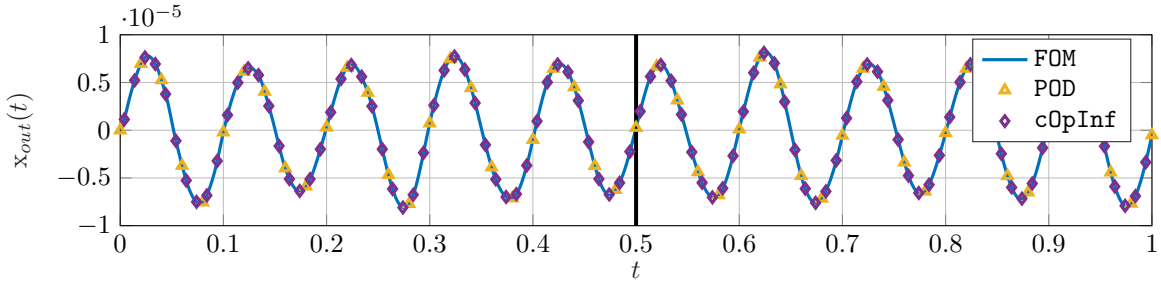}
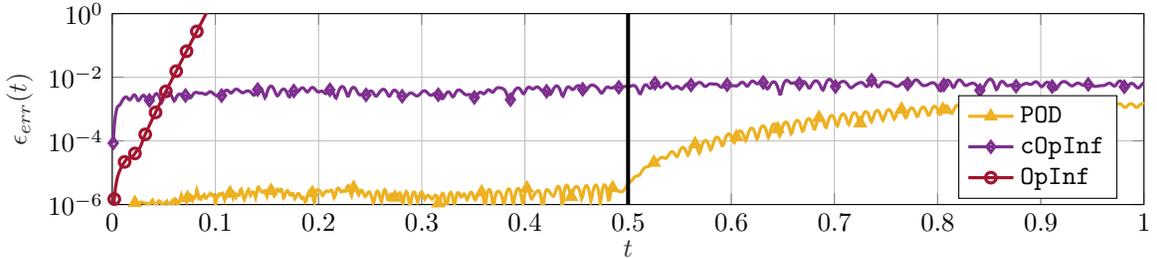
		\caption{Relative error in displacement vector}
		\label{fig:p_bot}
	\end{subfigure}
	\caption{Vibrating plate model: a comparison of FOM and ROMs of order $r = 110$. Black vertical line denotes the training period, used for constructing \texttt{POD, cOpInf}, and \texttt{OpInf} models.}
\end{figure}

The testing is performed for the same time-step and input load over a longer time interval $t = [ 0, \; 1 ]s$.
The qualitative comparison of the trajectories for ROMs of reduced order $r = 110$ is presented in \Cref{fig:p_top} for the displacement component $\mathrm{x}_{out} = \mathrm{x}_{176544}$; the relative error for the state trajectories is demonstrated in \Cref{fig:p_bot}. We can observe that the second-order operator inference formulation without constraints does not provide meaningful results for this example: the relative error blows up already during the training phase. In contrast, the force-informed operator inference model leads to a stable model with relative error below $1 \%$.
Although the POD model performs an order of magnitude better than the operator inference model, the accuracy of the POD model changes intermittently in the testing phase and reaches the force-informed operator inference level. This confirms the need to preserve the specific mathematical properties of the original system operators. Moreover, in the force-informed formulation, the stability of the model is guaranteed by the imposed constraints, which is not the case for the unconstrained version. 

\section{Conclusions}
\label{sec:conc}

In this paper, we have discussed extensions of the operator inference method incorporating the mechanical system structure of the governing equations. We presented a second-order formulation of operator inference, where the unknown operators can be identified using data. The asymptotic closeness of the inferred model to the corresponding intrusive POD model is also shown. An alternative formulation, as an optimization problem with positive semidefinite constraints for system operators, is proposed for the special case when the external force-data is available. The latter formulation allows ensuring stability of the inferred model. Both versions of operator inference provide reduced-order models that capture system dynamics very well. 

In this work, we provide the results only for the displacement field. However, the identification of stress-strain state might also be of interest. For this task, the access to the corresponding deformation data is necessary. Using the empirical knowledge about the strain-displacement relationship, it can be learned from the given deformation snapshots. Moreover, so far we assumed to have derivative data (e.g., velocity and acceleration) which may not be accessible. Therefore, in our future work, we explore approaches to use numerical approximation tools to approximate these quantities and analyze the effect of these on learning the operators.

\addcontentsline{toc}{section}{References}
\bibliographystyle{plainurl}
\bibliography{mor,mor_new,other}
  
\end{document}

%% file: fig/tikz/iss_svd2.tex
\definecolor{mycolor1}{rgb}{0.00000,0.44700,0.74100}%
\begin{tikzpicture}

\begin{axis}[%
width=2in,
height=1.5in,
at={(0.758in,0.481in)},
scale only axis,
xmin=0,
xmax=110,
ymode=log,
ymin=1e-20,
ymax=1,
yminorticks=true,
axis background/.style={fill=white},
xmajorgrids,
ymajorgrids,
yminorgrids
]
\addplot [color=mycolor1, only marks, mark=o, mark options={solid, mycolor1}, forget plot]
  table[row sep=crcr]{%
1	1\\
2	0.119535457052299\\
3	0.00530991533747063\\
4	0.0025679144904397\\
5	0.000759736409275429\\
6	0.000489661040213013\\
7	0.000259425484704114\\
8	0.000125527370027385\\
9	0.00010287134181834\\
10	5.69307950998576e-05\\
11	4.20222495521957e-05\\
12	3.00557374630174e-05\\
13	2.0044532439337e-05\\
14	1.43620680506833e-05\\
15	8.84845153691228e-06\\
16	7.29689945982626e-06\\
17	3.41887523104476e-06\\
18	2.56111984022317e-06\\
19	2.32121502039089e-06\\
20	1.38601514902331e-06\\
21	1.18294903420064e-06\\
22	1.08027457757727e-06\\
23	1.04959896045532e-06\\
24	8.18578921438449e-07\\
25	4.78980974210229e-07\\
26	4.02682228669751e-07\\
27	3.6712763457939e-07\\
28	3.35270524666118e-07\\
29	2.48929053413828e-07\\
30	1.83313195794847e-07\\
31	1.79463124581428e-07\\
32	1.37243655845408e-07\\
33	1.17478666803263e-07\\
34	1.0817617324015e-07\\
35	9.74202736034599e-08\\
36	8.78010947996829e-08\\
37	8.2785357086446e-08\\
38	7.19457716125854e-08\\
39	6.47027491052945e-08\\
40	6.17989134063995e-08\\
41	5.58092624831037e-08\\
42	5.14065672251498e-08\\
43	4.94933950137031e-08\\
44	4.30200059232518e-08\\
45	3.88584344303281e-08\\
46	3.28915453171006e-08\\
47	2.12891737436494e-08\\
48	1.73664833085644e-08\\
49	1.53201914335902e-08\\
50	1.41311056011432e-08\\
51	1.34483635464382e-08\\
52	9.55996119062379e-09\\
53	9.37265662914416e-09\\
54	8.96320852169209e-09\\
55	7.27547026329657e-09\\
56	5.15928248540938e-09\\
57	4.97866605029544e-09\\
58	3.753759861677e-09\\
59	2.82142063548867e-09\\
60	1.98777648092311e-09\\
61	1.73744718618774e-09\\
62	1.16014467200577e-09\\
63	1.05087964110994e-09\\
64	7.86597632397744e-10\\
65	6.72781861141762e-10\\
66	4.48174033548561e-10\\
67	3.96390437544521e-10\\
68	3.60917404533187e-10\\
69	3.42857652813416e-10\\
70	3.33256098213571e-10\\
71	2.51554228385249e-10\\
72	2.31431395951837e-10\\
73	1.07724134133728e-10\\
74	8.97313746487037e-11\\
75	8.92182610130979e-11\\
76	8.4106034502711e-11\\
77	8.0514828467606e-11\\
78	6.32646750747104e-11\\
79	2.84394903800374e-11\\
80	2.68025336372522e-11\\
81	1.56705551483321e-11\\
82	1.22258684616844e-11\\
83	4.31337961080628e-12\\
84	3.46693419782574e-12\\
85	1.70851816794524e-12\\
86	1.67607894975591e-12\\
87	1.49827200216022e-12\\
88	1.16184473739298e-12\\
89	1.04744905662499e-12\\
90	8.02445677154499e-13\\
91	5.75075021559618e-13\\
92	3.66826485925679e-13\\
93	3.24098815521555e-13\\
94	3.11540431973594e-13\\
95	2.16415462038961e-13\\
96	1.71181139086114e-13\\
97	1.24510677907462e-13\\
98	9.03381328105177e-14\\
99	2.25441072604364e-14\\
100	1.00911575090852e-14\\
101	8.3867148389163e-15\\
102	4.81207838033318e-15\\
103	1.26347076819285e-15\\
104	1.15966883603692e-15\\
105	1.0962149440695e-15\\
106	8.3734833436303e-16\\
107	1.16975356313621e-16\\
108	8.8974794921277e-17\\
109	8.8974794921277e-17\\
110	8.8974794921277e-17\\
111	8.8974794921277e-17\\
112	8.8974794921277e-17\\
113	8.8974794921277e-17\\
114	8.8974794921277e-17\\
115	8.8974794921277e-17\\
116	8.8974794921277e-17\\
117	8.8974794921277e-17\\
118	8.8974794921277e-17\\
119	8.8974794921277e-17\\
120	8.8974794921277e-17\\
121	8.8974794921277e-17\\
122	8.8974794921277e-17\\
123	8.8974794921277e-17\\
124	8.8974794921277e-17\\
125	8.8974794921277e-17\\
126	8.8974794921277e-17\\
127	8.8974794921277e-17\\
128	8.8974794921277e-17\\
129	8.8974794921277e-17\\
130	8.8974794921277e-17\\
131	8.8974794921277e-17\\
132	8.8974794921277e-17\\
133	8.8974794921277e-17\\
134	8.8974794921277e-17\\
135	4.67600276145841e-17\\
};
\addplot [black, only marks, mark=*, mark options={fill=black}, forget plot]
table[row sep=crcr]{%
4	0.0025679144904397\\
};
\end{axis}
\end{tikzpicture}%

%% file: fig/tikz/g_svd.tex
\definecolor{mycolor1}{rgb}{0.00000,0.44700,0.74100}%
\begin{tikzpicture}

\begin{axis}[%
width=2in,
height=1.5in,
at={(0.758in,0.481in)},
scale only axis,
xmin=0,
xmax=100,
ymode=log,
ymin=1e-16,
ymax=1,
%yminorticks=true,
axis background/.style={fill=white},
xmajorgrids,
ymajorgrids,
yminorgrids,
ytick = {1e-2,1e-6,1e-10,1e-14}
]
\addplot [color=mycolor1, only marks, mark=o, mark options={solid, mycolor1}, forget plot]
  table[row sep=crcr]{%
1	1\\
2	0.0438863163484598\\
3	0.014372229802962\\
4	0.0046089556091031\\
5	0.0031292149420254\\
6	0.00157416246222934\\
7	0.000395417586371964\\
8	0.00033168516375086\\
9	0.000314635232467864\\
10	3.76425183682347e-05\\
11	2.3301676851985e-05\\
12	9.69022639279845e-06\\
13	7.72672221940767e-06\\
14	5.10308225370075e-06\\
15	5.43650251453015e-07\\
16	2.77288740449891e-07\\
17	2.331710682864e-07\\
18	5.37810019526641e-08\\
19	1.50905551249969e-08\\
20	1.22372608923131e-08\\
21	2.0516792250315e-09\\
22	6.1415901357167e-10\\
23	2.78398219044244e-10\\
24	1.2380775770205e-10\\
25	8.10961368054017e-11\\
26	5.14316947094726e-11\\
27	2.0077051176641e-11\\
28	1.56908209181018e-11\\
29	5.02132696297727e-12\\
30	3.64501360758046e-12\\
31	1.44809453053852e-12\\
32	7.68233370112931e-13\\
33	4.79397239553136e-13\\
34	2.5909346780076e-13\\
35	1.08788443272731e-13\\
36	7.09056454999974e-14\\
37	4.40685363619248e-14\\
38	3.40083086375308e-14\\
39	2.50202729517311e-14\\
40	2.08502165166148e-14\\
41	1.94886883770732e-14\\
42	1.6686748103203e-14\\
43	1.1449655027433e-14\\
44	1.05497838926877e-14\\
45	9.64116024144564e-15\\
46	5.29693021645572e-15\\
47	4.635067605322e-15\\
48	4.53952654129456e-15\\
49	4.29967661707043e-15\\
50	4.12591781105002e-15\\
51	4.07073797758644e-15\\
52	3.59329351418621e-15\\
53	3.46211148318938e-15\\
54	3.07088365185288e-15\\
55	3.05611731828729e-15\\
56	2.86059362717659e-15\\
57	2.77241913565208e-15\\
58	2.54893133839472e-15\\
59	2.35440396446382e-15\\
60	2.17011456969961e-15\\
61	2.061029923772e-15\\
62	1.98736829479524e-15\\
63	1.61724582828537e-15\\
64	1.58781000024208e-15\\
65	1.35148385941082e-15\\
66	1.33716893682318e-15\\
67	1.27574611567117e-15\\
68	1.10371677068452e-15\\
69	1.06733470647702e-15\\
70	1.01754995086264e-15\\
71	9.61534591650569e-16\\
72	9.36159286063495e-16\\
73	9.21876372168555e-16\\
74	9.09841109684694e-16\\
75	8.88178419700125e-16\\
76	8.79512599550788e-16\\
77	8.52678861301162e-16\\
78	7.62813165128125e-16\\
79	7.58747802955007e-16\\
80	7.18705355866905e-16\\
81	7.030205389505e-16\\
82	6.78313382403548e-16\\
83	6.34838426754455e-16\\
84	5.61667790475436e-16\\
85	5.55098177716646e-16\\
86	5.38718561872341e-16\\
87	5.1089699649128e-16\\
88	5.02854181437433e-16\\
89	4.64592559205835e-16\\
90	4.58974085612485e-16\\
91	4.21103451915699e-16\\
92	4.07520789369127e-16\\
93	3.86257186719564e-16\\
94	3.78896219956811e-16\\
95	3.78778017021649e-16\\
96	3.70513525371957e-16\\
97	3.65462166156272e-16\\
98	3.58875258098307e-16\\
99	3.47509241310364e-16\\
100	3.39381222174533e-16\\
101	3.31600576767758e-16\\
102	3.24824233753664e-16\\
103	3.23877327357232e-16\\
104	3.18381393838586e-16\\
105	3.09403599923217e-16\\
106	2.99918165165058e-16\\
107	2.87785919990403e-16\\
108	2.79651121947173e-16\\
109	2.7740229046906e-16\\
110	2.75331719897197e-16\\
111	2.68048082210884e-16\\
112	2.56883333783535e-16\\
113	2.44538547700357e-16\\
114	2.39664287768889e-16\\
115	2.31750913029852e-16\\
116	2.25536907816903e-16\\
117	2.21631935041303e-16\\
118	2.14266902300709e-16\\
119	2.071938641013e-16\\
120	2.04850924967841e-16\\
121	1.98944913290206e-16\\
122	1.94434526353307e-16\\
123	1.82554915833586e-16\\
124	1.75432123083937e-16\\
125	1.66017216834219e-16\\
126	1.63656338065025e-16\\
127	1.61900927888506e-16\\
128	1.61048349828062e-16\\
129	1.60834645671826e-16\\
130	1.5305761320244e-16\\
131	1.42057604069267e-16\\
132	1.3456778807496e-16\\
133	1.28477555787648e-16\\
134	1.20687509327799e-16\\
135	1.11778672511181e-16\\
136	1.01393689543438e-16\\
137	9.99143591400114e-17\\
138	9.99143154578862e-17\\
139	9.99143154578862e-17\\
140	9.99143154578862e-17\\
141	9.99143154578862e-17\\
142	9.99143154578862e-17\\
143	9.99143154578862e-17\\
144	9.99143154578862e-17\\
145	9.99143154578862e-17\\
146	9.99143154578862e-17\\
147	9.99143154578862e-17\\
148	9.99143154578862e-17\\
149	9.99143154578862e-17\\
150	9.99143154578862e-17\\
151	9.99143154578862e-17\\
152	9.99143154578862e-17\\
153	9.99143154578862e-17\\
154	9.99143154578862e-17\\
155	9.99143154578862e-17\\
156	9.99143154578862e-17\\
157	9.99143154578862e-17\\
158	9.99143154578862e-17\\
159	9.99143154578862e-17\\
160	9.99143154578862e-17\\
161	9.99143154578862e-17\\
162	9.99143154578862e-17\\
163	9.99143154578862e-17\\
164	9.99143154578862e-17\\
165	9.99143154578862e-17\\
166	9.99143154578862e-17\\
167	9.99143154578862e-17\\
168	9.99143154578862e-17\\
169	9.99143154578862e-17\\
170	9.99143154578862e-17\\
171	9.99143154578862e-17\\
172	9.99143154578862e-17\\
173	9.99143154578862e-17\\
174	9.99143154578862e-17\\
175	9.99143154578862e-17\\
176	9.99143154578862e-17\\
177	9.99143154578862e-17\\
178	9.99143154578862e-17\\
179	9.99143154578862e-17\\
180	9.99143154578862e-17\\
181	9.99143154578862e-17\\
182	9.99143154578862e-17\\
183	9.99143154578862e-17\\
184	9.99143154578862e-17\\
185	9.99143154578862e-17\\
186	9.99143154578862e-17\\
187	9.99143154578862e-17\\
188	9.99143154578862e-17\\
189	9.99143154578862e-17\\
190	9.99143154578862e-17\\
191	9.99143154578862e-17\\
192	9.99143154578862e-17\\
193	9.99143154578862e-17\\
194	9.99143154578862e-17\\
195	9.99143154578862e-17\\
196	9.99143154578862e-17\\
197	9.99143154578862e-17\\
198	9.99143154578862e-17\\
199	9.99143154578862e-17\\
200	9.99143154578862e-17\\
201	9.99143154578862e-17\\
202	9.99143154578862e-17\\
203	9.99143154578862e-17\\
204	9.99143154578862e-17\\
205	9.99143154578862e-17\\
206	9.99143154578862e-17\\
207	9.99143154578862e-17\\
208	9.99143154578862e-17\\
209	9.99143154578862e-17\\
210	9.99143154578862e-17\\
211	9.99143154578862e-17\\
212	9.99143154578862e-17\\
213	9.99143154578862e-17\\
214	9.99143154578862e-17\\
215	9.99143154578862e-17\\
216	9.99143154578862e-17\\
217	9.99143154578862e-17\\
218	9.99143154578862e-17\\
219	9.99143154578862e-17\\
220	9.99143154578862e-17\\
221	9.99143154578862e-17\\
222	9.99143154578862e-17\\
223	9.99143154578862e-17\\
224	9.99143154578862e-17\\
225	9.99143154578862e-17\\
226	9.99143154578862e-17\\
227	9.99143154578862e-17\\
228	9.99143154578862e-17\\
229	9.99143154578862e-17\\
230	9.99143154578862e-17\\
231	9.99143154578862e-17\\
232	9.99143154578862e-17\\
233	9.99143154578862e-17\\
234	9.99143154578862e-17\\
235	9.99143154578862e-17\\
236	9.99143154578862e-17\\
237	9.99143154578862e-17\\
238	9.99143154578862e-17\\
239	9.99143154578862e-17\\
240	9.99143154578862e-17\\
241	9.99143154578862e-17\\
242	9.99143154578862e-17\\
243	9.99143154578862e-17\\
244	9.99143154578862e-17\\
245	9.99143154578862e-17\\
246	9.99143154578862e-17\\
247	9.99143154578862e-17\\
248	9.99143154578862e-17\\
249	9.99143154578862e-17\\
250	9.99143154578862e-17\\
251	9.99143154578862e-17\\
252	9.99143154578862e-17\\
253	9.99143154578862e-17\\
254	9.99143154578862e-17\\
255	9.99143154578862e-17\\
256	9.99143154578862e-17\\
257	9.99143154578862e-17\\
258	9.99143154578862e-17\\
259	9.99143154578862e-17\\
260	9.99143154578862e-17\\
261	9.99143154578862e-17\\
262	9.99143154578862e-17\\
263	9.99143154578862e-17\\
264	9.99143154578862e-17\\
265	9.99143154578862e-17\\
266	9.99143154578862e-17\\
267	9.99143154578862e-17\\
268	9.99143154578862e-17\\
269	9.99143154578862e-17\\
270	9.99143154578862e-17\\
271	9.99143154578862e-17\\
272	9.99143154578862e-17\\
273	9.99143154578862e-17\\
274	9.99143154578862e-17\\
275	9.99143154578862e-17\\
276	9.99143154578862e-17\\
277	9.99143154578862e-17\\
278	9.99143154578862e-17\\
279	9.99143154578862e-17\\
280	9.99143154578862e-17\\
281	9.99143154578862e-17\\
282	9.99143154578862e-17\\
283	9.99143154578862e-17\\
284	9.99143154578862e-17\\
285	9.99143154578862e-17\\
286	9.99143154578862e-17\\
287	9.99143154578862e-17\\
288	9.99143154578862e-17\\
289	9.99143154578862e-17\\
290	9.99143154578862e-17\\
291	9.99143154578862e-17\\
292	9.99143154578862e-17\\
293	9.99143154578862e-17\\
294	9.99143154578862e-17\\
295	9.99143154578862e-17\\
296	9.99143154578862e-17\\
297	9.99143154578862e-17\\
298	9.99143154578862e-17\\
299	9.99143154578862e-17\\
300	9.99143154578862e-17\\
301	9.99143154578862e-17\\
302	9.99143154578862e-17\\
303	9.99143154578862e-17\\
304	9.99143154578862e-17\\
305	9.99143154578862e-17\\
306	9.99143154578862e-17\\
307	9.99143154578862e-17\\
308	9.99143154578862e-17\\
309	9.99143154578862e-17\\
310	9.99143154578862e-17\\
311	9.99143154578862e-17\\
312	9.99143154578862e-17\\
313	9.99143154578862e-17\\
314	9.99143154578862e-17\\
315	9.99143154578862e-17\\
316	9.99143154578862e-17\\
317	9.99143154578862e-17\\
318	9.99143154578862e-17\\
319	9.99143154578862e-17\\
320	9.99143154578862e-17\\
321	9.99143154578862e-17\\
322	9.99143154578862e-17\\
323	9.99143154578862e-17\\
324	9.99143154578862e-17\\
325	9.99143154578862e-17\\
326	9.99143154578862e-17\\
327	9.99143154578862e-17\\
328	9.99143154578862e-17\\
329	9.99143154578862e-17\\
330	9.99143154578862e-17\\
331	9.99143154578862e-17\\
332	9.99143154578862e-17\\
333	9.99143154578862e-17\\
334	9.99143154578862e-17\\
335	9.99143154578862e-17\\
336	9.99143154578862e-17\\
337	9.99143154578862e-17\\
338	9.99143154578862e-17\\
339	9.99143154578862e-17\\
340	9.99143154578862e-17\\
341	9.99143154578862e-17\\
342	9.99143154578862e-17\\
343	9.99143154578862e-17\\
344	9.99143154578862e-17\\
345	9.99143154578862e-17\\
346	9.99143154578862e-17\\
347	9.99143154578862e-17\\
348	9.99143154578862e-17\\
349	9.99143154578862e-17\\
350	9.99143154578862e-17\\
351	9.99143154578862e-17\\
352	9.99143154578862e-17\\
353	9.99143154578862e-17\\
354	9.99143154578862e-17\\
355	9.99143154578862e-17\\
356	9.99143154578862e-17\\
357	9.99143154578862e-17\\
358	9.99143154578862e-17\\
359	9.99143154578862e-17\\
360	9.99143154578862e-17\\
361	9.99143154578862e-17\\
362	9.99143154578862e-17\\
363	9.99143154578862e-17\\
364	9.99143154578862e-17\\
365	9.99143154578862e-17\\
366	9.99143154578862e-17\\
367	9.99143154578862e-17\\
368	9.99143154578862e-17\\
369	9.99143154578862e-17\\
370	9.99143154578862e-17\\
371	9.99143154578862e-17\\
372	9.99143154578862e-17\\
373	9.99143154578862e-17\\
374	9.99143154578862e-17\\
375	9.99143154578862e-17\\
376	9.99143154578862e-17\\
377	9.99143154578862e-17\\
378	9.99143154578862e-17\\
379	9.99143154578862e-17\\
380	9.99143154578862e-17\\
381	9.99143154578862e-17\\
382	9.99143154578862e-17\\
383	9.99143154578862e-17\\
384	9.99143154578862e-17\\
385	9.99143154578862e-17\\
386	9.99143154578862e-17\\
387	9.99143154578862e-17\\
388	9.99143154578862e-17\\
389	9.99143154578862e-17\\
390	9.99143154578862e-17\\
391	9.99143154578862e-17\\
392	9.99143154578862e-17\\
393	9.99143154578862e-17\\
394	9.99143154578862e-17\\
395	9.99143154578862e-17\\
396	9.99143154578862e-17\\
397	9.99143154578862e-17\\
398	9.99143154578862e-17\\
399	9.99143154578862e-17\\
400	9.99143154578862e-17\\
401	9.99143154578862e-17\\
402	9.99143154578862e-17\\
403	9.99143154578862e-17\\
404	9.99143154578862e-17\\
405	9.99143154578862e-17\\
406	9.99143154578862e-17\\
407	9.99143154578862e-17\\
408	9.99143154578862e-17\\
409	9.99143154578862e-17\\
410	9.99143154578862e-17\\
411	9.99143154578862e-17\\
412	9.99143154578862e-17\\
413	9.99143154578862e-17\\
414	9.99143154578862e-17\\
415	9.99143154578862e-17\\
416	9.99143154578862e-17\\
417	9.99143154578862e-17\\
418	9.99143154578862e-17\\
419	9.99143154578862e-17\\
420	9.99143154578862e-17\\
421	9.99143154578862e-17\\
422	9.99143154578862e-17\\
423	9.99143154578862e-17\\
424	9.99143154578862e-17\\
425	9.99143154578862e-17\\
426	9.99143154578862e-17\\
427	9.99143154578862e-17\\
428	9.99143154578862e-17\\
429	9.99143154578862e-17\\
430	9.99143154578862e-17\\
431	9.99143154578862e-17\\
432	9.99143154578862e-17\\
433	9.99143154578862e-17\\
434	9.99143154578862e-17\\
435	9.99143154578862e-17\\
436	9.99143154578862e-17\\
437	9.99143154578862e-17\\
438	9.99143154578862e-17\\
439	9.99143154578862e-17\\
440	9.99143154578862e-17\\
441	9.99143154578862e-17\\
442	9.99143154578862e-17\\
443	9.99143154578862e-17\\
444	9.99143154578862e-17\\
445	9.99143154578862e-17\\
446	9.99143154578862e-17\\
447	9.99143154578862e-17\\
448	9.99143154578862e-17\\
449	9.99143154578862e-17\\
450	9.99143154578862e-17\\
451	9.99143154578862e-17\\
452	9.99143154578862e-17\\
453	9.99143154578862e-17\\
454	9.99143154578862e-17\\
455	9.99143154578862e-17\\
456	9.99143154578862e-17\\
457	9.99143154578862e-17\\
458	9.99143154578862e-17\\
459	9.99143154578862e-17\\
460	9.99143154578862e-17\\
461	9.99143154578862e-17\\
462	9.99143154578862e-17\\
463	9.99143154578862e-17\\
464	9.99143154578862e-17\\
465	9.99143154578862e-17\\
466	9.99143154578862e-17\\
467	9.99143154578862e-17\\
468	9.99143154578862e-17\\
469	9.99143154578862e-17\\
470	9.99143154578862e-17\\
471	9.99143154578862e-17\\
472	9.99143154578862e-17\\
473	9.99143154578862e-17\\
474	9.99143154578862e-17\\
475	9.99143154578862e-17\\
476	9.99143154578862e-17\\
477	9.99143154578862e-17\\
478	9.99143154578862e-17\\
479	9.99143154578862e-17\\
480	9.99143154578862e-17\\
481	9.99143154578862e-17\\
482	9.99143154578862e-17\\
483	9.99143154578862e-17\\
484	9.99143154578862e-17\\
485	9.99143154578862e-17\\
486	9.99143154578862e-17\\
487	9.99143154578862e-17\\
488	9.99143154578862e-17\\
489	9.99143154578862e-17\\
490	9.99143154578862e-17\\
491	9.99143154578862e-17\\
492	9.99143154578862e-17\\
493	9.99143154578862e-17\\
494	9.99143154578862e-17\\
495	9.99143154578862e-17\\
496	9.99143154578862e-17\\
497	9.99143154578862e-17\\
498	9.99143154578862e-17\\
499	9.99143154578862e-17\\
500	9.99143154578862e-17\\
501	9.99143154578862e-17\\
502	9.99143154578862e-17\\
503	9.99143154578862e-17\\
504	9.99143154578862e-17\\
505	9.99143154578862e-17\\
506	9.99143154578862e-17\\
507	9.99143154578862e-17\\
508	9.99143154578862e-17\\
509	9.99143154578862e-17\\
510	9.99143154578862e-17\\
511	9.99143154578862e-17\\
512	9.99143154578862e-17\\
513	9.99143154578862e-17\\
514	9.99143154578862e-17\\
515	9.99143154578862e-17\\
516	9.99143154578862e-17\\
517	9.99143154578862e-17\\
518	9.99143154578862e-17\\
519	9.99143154578862e-17\\
520	9.99143154578862e-17\\
521	9.99143154578862e-17\\
522	9.99143154578862e-17\\
523	9.99143154578862e-17\\
524	9.99143154578862e-17\\
525	9.99143154578862e-17\\
526	9.99143154578862e-17\\
527	9.99143154578862e-17\\
528	9.99143154578862e-17\\
529	9.99143154578862e-17\\
530	9.99143154578862e-17\\
531	9.99143154578862e-17\\
532	9.99143154578862e-17\\
533	9.99143154578862e-17\\
534	9.99143154578862e-17\\
535	9.99143154578862e-17\\
536	9.99143154578862e-17\\
537	9.99143154578862e-17\\
538	9.99143154578862e-17\\
539	9.99143154578862e-17\\
540	9.99143154578862e-17\\
541	9.99143154578862e-17\\
542	9.99143154578862e-17\\
543	9.99143154578862e-17\\
544	9.99143154578862e-17\\
545	9.99143154578862e-17\\
546	9.99143154578862e-17\\
547	9.99143154578862e-17\\
548	9.99143154578862e-17\\
549	9.99143154578862e-17\\
550	9.99143154578862e-17\\
551	9.99143154578862e-17\\
552	9.99143154578862e-17\\
553	9.99143154578862e-17\\
554	9.99143154578862e-17\\
555	9.99143154578862e-17\\
556	9.99143154578862e-17\\
557	9.99143154578862e-17\\
558	9.99143154578862e-17\\
559	9.99143154578862e-17\\
560	9.99143154578862e-17\\
561	9.99143154578862e-17\\
562	9.99143154578862e-17\\
563	9.99143154578862e-17\\
564	9.99143154578862e-17\\
565	9.99143154578862e-17\\
566	9.99143154578862e-17\\
567	9.99143154578862e-17\\
568	9.99143154578862e-17\\
569	9.99143154578862e-17\\
570	9.99143154578862e-17\\
571	9.99143154578862e-17\\
572	9.99143154578862e-17\\
573	9.99143154578862e-17\\
574	9.99143154578862e-17\\
575	9.99143154578862e-17\\
576	9.99143154578862e-17\\
577	9.99143154578862e-17\\
578	9.99143154578862e-17\\
579	9.99143154578862e-17\\
580	9.99143154578862e-17\\
581	9.99143154578862e-17\\
582	9.99143154578862e-17\\
583	9.99143154578862e-17\\
584	9.99143154578862e-17\\
585	9.99143154578862e-17\\
586	9.99143154578862e-17\\
587	9.99143154578862e-17\\
588	9.99143154578862e-17\\
589	9.99143154578862e-17\\
590	9.99143154578862e-17\\
591	9.99143154578862e-17\\
592	9.99143154578862e-17\\
593	9.99143154578862e-17\\
594	9.99143154578862e-17\\
595	9.99143154578862e-17\\
596	9.99143154578862e-17\\
597	9.99143154578862e-17\\
598	9.99143154578862e-17\\
599	9.99143154578862e-17\\
600	9.99143154578862e-17\\
601	9.99143154578862e-17\\
602	9.99143154578862e-17\\
603	9.99143154578862e-17\\
604	9.99143154578862e-17\\
605	9.99143154578862e-17\\
606	9.99143154578862e-17\\
607	9.99143154578862e-17\\
608	9.99143154578862e-17\\
609	9.99143154578862e-17\\
610	9.99143154578862e-17\\
611	9.99143154578862e-17\\
612	9.99143154578862e-17\\
613	9.99143154578862e-17\\
614	9.99143154578862e-17\\
615	9.99143154578862e-17\\
616	9.99143154578862e-17\\
617	9.99143154578862e-17\\
618	9.99143154578862e-17\\
619	9.99143154578862e-17\\
620	9.99143154578862e-17\\
621	9.99143154578862e-17\\
622	9.99143154578862e-17\\
623	9.99143154578862e-17\\
624	9.99143154578862e-17\\
625	9.99143154578862e-17\\
626	9.99143154578862e-17\\
627	9.99143154578862e-17\\
628	9.99143154578862e-17\\
629	9.99143154578862e-17\\
630	9.99143154578862e-17\\
631	9.99143154578862e-17\\
632	9.99143154578862e-17\\
633	9.99143154578862e-17\\
634	9.99143154578862e-17\\
635	9.99143154578862e-17\\
636	9.99143154578862e-17\\
637	9.99143154578862e-17\\
638	9.99143154578862e-17\\
639	9.99143154578862e-17\\
640	9.99143154578862e-17\\
641	9.99143154578862e-17\\
642	9.99143154578862e-17\\
643	9.99143154578862e-17\\
644	9.99143154578862e-17\\
645	9.99143154578862e-17\\
646	9.99143154578862e-17\\
647	9.99143154578862e-17\\
648	9.99143154578862e-17\\
649	9.99143154578862e-17\\
650	9.99143154578862e-17\\
651	9.99143154578862e-17\\
652	9.99143154578862e-17\\
653	9.99143154578862e-17\\
654	9.99143154578862e-17\\
655	9.99143154578862e-17\\
656	9.99143154578862e-17\\
657	9.99143154578862e-17\\
658	9.99143154578862e-17\\
659	9.99143154578862e-17\\
660	9.99143154578862e-17\\
661	9.99143154578862e-17\\
662	9.99143154578862e-17\\
663	9.99143154578862e-17\\
664	9.99143154578862e-17\\
665	9.99143154578862e-17\\
666	9.99143154578862e-17\\
667	9.99143154578862e-17\\
668	9.99143154578862e-17\\
669	9.99143154578862e-17\\
670	9.99143154578862e-17\\
671	9.99143154578862e-17\\
672	9.99143154578862e-17\\
673	9.99143154578862e-17\\
674	9.99143154578862e-17\\
675	9.99143154578862e-17\\
676	9.99143154578862e-17\\
677	9.99143154578862e-17\\
678	9.99143154578862e-17\\
679	9.99143154578862e-17\\
680	9.99143154578862e-17\\
681	9.99143154578862e-17\\
682	9.99143154578862e-17\\
683	9.99143154578862e-17\\
684	9.99143154578862e-17\\
685	9.99143154578862e-17\\
686	9.99143154578862e-17\\
687	9.99143154578862e-17\\
688	9.99143154578862e-17\\
689	9.99143154578862e-17\\
690	9.99143154578862e-17\\
691	9.99143154578862e-17\\
692	9.99143154578862e-17\\
693	9.99143154578862e-17\\
694	9.99143154578862e-17\\
695	9.99143154578862e-17\\
696	9.99143154578862e-17\\
697	9.99143154578862e-17\\
698	9.99143154578862e-17\\
699	9.99143154578862e-17\\
700	9.99143154578862e-17\\
701	9.99143154578862e-17\\
702	9.99143154578862e-17\\
703	9.99143154578862e-17\\
704	9.99143154578862e-17\\
705	9.99143154578862e-17\\
706	9.99143154578862e-17\\
707	9.99143154578862e-17\\
708	9.99143154578862e-17\\
709	9.99143154578862e-17\\
710	9.99143154578862e-17\\
711	9.99143154578862e-17\\
712	9.99143154578862e-17\\
713	9.99143154578862e-17\\
714	9.99143154578862e-17\\
715	9.99143154578862e-17\\
716	9.99143154578862e-17\\
717	9.99143154578862e-17\\
718	9.99143154578862e-17\\
719	9.99143154578862e-17\\
720	9.99143154578862e-17\\
721	9.99143154578862e-17\\
722	9.99143154578862e-17\\
723	9.99143154578862e-17\\
724	9.99143154578862e-17\\
725	9.99143154578862e-17\\
726	9.99143154578862e-17\\
727	9.99143154578862e-17\\
728	9.99143154578862e-17\\
729	9.99143154578862e-17\\
730	9.99143154578862e-17\\
731	9.99143154578862e-17\\
732	9.99143154578862e-17\\
733	9.99143154578862e-17\\
734	9.99143154578862e-17\\
735	9.99143154578862e-17\\
736	9.99143154578862e-17\\
737	9.99143154578862e-17\\
738	9.99143154578862e-17\\
739	9.99143154578862e-17\\
740	9.99143154578862e-17\\
741	9.99143154578862e-17\\
742	9.99143154578862e-17\\
743	9.99143154578862e-17\\
744	9.99143154578862e-17\\
745	9.99143154578862e-17\\
746	9.99143154578862e-17\\
747	9.99143154578862e-17\\
748	9.99143154578862e-17\\
749	9.99143154578862e-17\\
750	9.99143154578862e-17\\
751	9.99143154578862e-17\\
752	9.99143154578862e-17\\
753	9.99143154578862e-17\\
754	9.99143154578862e-17\\
755	9.99143154578862e-17\\
756	9.99143154578862e-17\\
757	9.99143154578862e-17\\
758	9.99143154578862e-17\\
759	9.99143154578862e-17\\
760	9.99143154578862e-17\\
761	9.99143154578862e-17\\
762	9.99143154578862e-17\\
763	9.99143154578862e-17\\
764	9.99143154578862e-17\\
765	9.99143154578862e-17\\
766	9.99143154578862e-17\\
767	9.99143154578862e-17\\
768	9.99143154578862e-17\\
769	9.99143154578862e-17\\
770	9.99143154578862e-17\\
771	9.99143154578862e-17\\
772	9.99143154578862e-17\\
773	9.99143154578862e-17\\
774	9.99143154578862e-17\\
775	9.99143154578862e-17\\
776	9.99143154578862e-17\\
777	9.99143154578862e-17\\
778	9.99143154578862e-17\\
779	9.99143154578862e-17\\
780	9.99143154578862e-17\\
781	9.99143154578862e-17\\
782	9.99143154578862e-17\\
783	9.99143154578862e-17\\
784	9.99143154578862e-17\\
785	9.99143154578862e-17\\
786	9.99143154578862e-17\\
787	9.99143154578862e-17\\
788	9.99143154578862e-17\\
789	9.99143154578862e-17\\
790	9.99143154578862e-17\\
791	9.99143154578862e-17\\
792	9.99143154578862e-17\\
793	9.99143154578862e-17\\
794	9.99143154578862e-17\\
795	9.99143154578862e-17\\
796	9.99143154578862e-17\\
797	9.99143154578862e-17\\
798	9.99143154578862e-17\\
799	9.99143154578862e-17\\
800	9.99143154578862e-17\\
801	9.99143154578862e-17\\
802	9.99143154578862e-17\\
803	9.99143154578862e-17\\
804	9.99143154578862e-17\\
805	9.99143154578862e-17\\
806	9.99143154578862e-17\\
807	9.99143154578862e-17\\
808	9.99143154578862e-17\\
809	9.99143154578862e-17\\
810	9.99143154578862e-17\\
811	9.99143154578862e-17\\
812	9.99143154578862e-17\\
813	9.99143154578862e-17\\
814	9.99143154578862e-17\\
815	9.99143154578862e-17\\
816	9.99143154578862e-17\\
817	9.99143154578862e-17\\
818	9.99143154578862e-17\\
819	9.99143154578862e-17\\
820	9.99143154578862e-17\\
821	9.99143154578862e-17\\
822	9.99143154578862e-17\\
823	9.99143154578862e-17\\
824	9.99143154578862e-17\\
825	9.99143154578862e-17\\
826	9.99143154578862e-17\\
827	9.99143154578862e-17\\
828	9.99143154578862e-17\\
829	9.99143154578862e-17\\
830	9.99143154578862e-17\\
831	9.99143154578862e-17\\
832	9.99143154578862e-17\\
833	9.99143154578862e-17\\
834	9.99143154578862e-17\\
835	9.99143154578862e-17\\
836	9.99143154578862e-17\\
837	9.99143154578862e-17\\
838	9.99143154578862e-17\\
839	9.99143154578862e-17\\
840	9.99143154578862e-17\\
841	9.99143154578862e-17\\
842	9.99143154578862e-17\\
843	9.99143154578862e-17\\
844	9.99143154578862e-17\\
845	9.99143154578862e-17\\
846	9.99143154578862e-17\\
847	9.99143154578862e-17\\
848	9.99143154578862e-17\\
849	9.99143154578862e-17\\
850	9.99143154578862e-17\\
851	9.99143154578862e-17\\
852	9.99143154578862e-17\\
853	9.99143154578862e-17\\
854	9.99143154578862e-17\\
855	9.99143154578862e-17\\
856	9.99143154578862e-17\\
857	9.99143154578862e-17\\
858	9.99143154578862e-17\\
859	9.99143154578862e-17\\
860	9.99143154578862e-17\\
861	9.99143154578862e-17\\
862	9.99143154578862e-17\\
863	9.99143154578862e-17\\
864	9.99143154578862e-17\\
865	9.99143154578862e-17\\
866	9.99143154578862e-17\\
867	9.99143154578862e-17\\
868	9.99143154578862e-17\\
869	9.99143154578862e-17\\
870	9.99143154578862e-17\\
871	9.99143154578862e-17\\
872	9.99143154578862e-17\\
873	9.99143154578862e-17\\
874	9.99143154578862e-17\\
875	9.99143154578862e-17\\
876	9.99143154578862e-17\\
877	9.99143154578862e-17\\
878	9.99143154578862e-17\\
879	9.99143154578862e-17\\
880	9.99143154578862e-17\\
881	9.99143154578862e-17\\
882	9.99143154578862e-17\\
883	9.99143154578862e-17\\
884	9.99143154578862e-17\\
885	9.99143154578862e-17\\
886	9.99143154578862e-17\\
887	9.99143154578862e-17\\
888	9.99143154578862e-17\\
889	9.99143154578862e-17\\
890	9.99143154578862e-17\\
891	9.99143154578862e-17\\
892	9.99143154578862e-17\\
893	9.99143154578862e-17\\
894	9.99143154578862e-17\\
895	9.99143154578862e-17\\
896	9.99143154578862e-17\\
897	9.99143154578862e-17\\
898	9.99143154578862e-17\\
899	9.99143154578862e-17\\
900	9.99143154578862e-17\\
901	9.99143154578862e-17\\
902	9.99143154578862e-17\\
903	9.99143154578862e-17\\
904	9.99143154578862e-17\\
905	9.99143154578862e-17\\
906	9.99143154578862e-17\\
907	9.99143154578862e-17\\
908	9.99143154578862e-17\\
909	9.99143154578862e-17\\
910	9.99143154578862e-17\\
911	9.99143154578862e-17\\
912	9.99143154578862e-17\\
913	9.99143154578862e-17\\
914	9.99143154578862e-17\\
915	9.99143154578862e-17\\
916	9.99143154578862e-17\\
917	9.99143154578862e-17\\
918	9.99143154578862e-17\\
919	9.99143154578862e-17\\
920	9.99143154578862e-17\\
921	9.99143154578862e-17\\
922	9.99143154578862e-17\\
923	9.99143154578862e-17\\
924	9.99143154578862e-17\\
925	9.99143154578862e-17\\
926	9.99143154578862e-17\\
927	9.99143154578862e-17\\
928	9.99143154578862e-17\\
929	9.99143154578862e-17\\
930	9.99143154578862e-17\\
931	9.99143154578862e-17\\
932	9.99143154578862e-17\\
933	9.99143154578862e-17\\
934	9.99143154578862e-17\\
935	9.99143154578862e-17\\
936	9.99143154578862e-17\\
937	9.99143154578862e-17\\
938	9.99143154578862e-17\\
939	9.99143154578862e-17\\
940	9.99143154578862e-17\\
941	9.99143154578862e-17\\
942	9.99143154578862e-17\\
943	9.99143154578862e-17\\
944	9.99143154578862e-17\\
945	9.99143154578862e-17\\
946	9.99143154578862e-17\\
947	9.99143154578862e-17\\
948	9.99143154578862e-17\\
949	9.99143154578862e-17\\
950	9.99143154578862e-17\\
951	9.99143154578862e-17\\
952	9.99143154578862e-17\\
953	9.99143154578862e-17\\
954	9.99143154578862e-17\\
955	9.99143154578862e-17\\
956	9.99143154578862e-17\\
957	9.99143154578862e-17\\
958	9.99143154578862e-17\\
959	9.99143154578862e-17\\
960	9.99143154578862e-17\\
961	9.99143154578862e-17\\
962	9.99143154578862e-17\\
963	9.99143154578862e-17\\
964	9.99143154578862e-17\\
965	9.99143154578862e-17\\
966	9.99143154578862e-17\\
967	9.99143154578862e-17\\
968	9.99143154578862e-17\\
969	9.99143154578862e-17\\
970	9.99143154578862e-17\\
971	9.99143154578862e-17\\
972	9.99143154578862e-17\\
973	9.99143154578862e-17\\
974	9.99143154578862e-17\\
975	9.99143154578862e-17\\
976	9.99143154578862e-17\\
977	9.99143154578862e-17\\
978	9.99143154578862e-17\\
979	9.99143154578862e-17\\
980	9.99143154578862e-17\\
981	9.99143154578862e-17\\
982	9.99143154578862e-17\\
983	9.99143154578862e-17\\
984	9.99143154578862e-17\\
985	9.99143154578862e-17\\
986	9.99143154578862e-17\\
987	9.99143154578862e-17\\
988	9.99143154578862e-17\\
989	9.44849419724004e-17\\
990	8.64111667610022e-17\\
991	7.81400466128003e-17\\
992	6.91288013046745e-17\\
993	6.17903137508243e-17\\
994	6.16603657780933e-17\\
995	6.06550500105508e-17\\
996	5.07101598539977e-17\\
997	4.2272633726845e-17\\
998	3.22046746130135e-17\\
999	2.34908888560978e-17\\
1000	1.51687123006812e-17\\
1001	7.5124988487464e-18\\
};
\addplot [only marks, mark=*, mark options={solid, black}, forget plot]
table[row sep=crcr]{%
6	0.00157416246222934\\
};
\end{axis}
\end{tikzpicture}%

%% file: fig/tikz/p_svd.tex
\definecolor{mycolor1}{rgb}{0.00000,0.44700,0.74100}%
\begin{tikzpicture}

\begin{axis}[%
width=2in,
height=1.5in,
at={(0.772in,0.473in)},
scale only axis,
xmin=0,
xmax=500,
ymode=log,
ymin=1e-14,
ymax=1,
yminorticks=true,
axis background/.style={fill=white},
xmajorgrids,
ymajorgrids,
yminorgrids
]
\addplot [color=mycolor1, only marks, mark=o, mark options={solid, mycolor1},mark repeat = 10, forget plot]
  table[row sep=crcr]{%
1	1\\
2	0.507069157544121\\
3	0.214516072318151\\
4	0.112736356971856\\
5	0.0598023588922573\\
6	0.0315714464471682\\
7	0.0196657964795942\\
8	0.0175032136678707\\
9	0.0147394431734811\\
10	0.00764068949094215\\
11	0.00450895058147984\\
12	0.00418781125722247\\
13	0.00385432033965219\\
14	0.00301278884931264\\
15	0.00214526494139428\\
16	0.00179800720289494\\
17	0.00155992583841016\\
18	0.00104655804243617\\
19	0.00100210574246961\\
20	0.000902819853390924\\
21	0.000794008975203279\\
22	0.000759415385418377\\
23	0.000651751915867258\\
24	0.000572628597949653\\
25	0.000520658634368858\\
26	0.000450912202300554\\
27	0.000416371242242538\\
28	0.000401782442337834\\
29	0.000382120092681456\\
30	0.000321572797339776\\
31	0.000308684136743461\\
32	0.000306141056686156\\
33	0.000290173132126495\\
34	0.000255882320412523\\
35	0.000237422930597657\\
36	0.000230893010079759\\
37	0.000206500532846231\\
38	0.000187222571236034\\
39	0.000161855460378124\\
40	0.00015457606682326\\
41	0.000139159390110192\\
42	0.000124814744857087\\
43	0.000118469391779567\\
44	0.00010175831412503\\
45	8.79296998874787e-05\\
46	8.21591744711146e-05\\
47	6.9917311457414e-05\\
48	5.82174000376791e-05\\
49	5.7279591833801e-05\\
50	5.08130826774939e-05\\
51	4.80618563755919e-05\\
52	4.40374041871566e-05\\
53	4.08530373655117e-05\\
54	3.73741109071205e-05\\
55	3.16114795149033e-05\\
56	3.04504775176699e-05\\
57	2.6187165385634e-05\\
58	2.54375077904963e-05\\
59	2.15988168814602e-05\\
60	2.07424671339307e-05\\
61	1.71714327367829e-05\\
62	1.67244146011892e-05\\
63	1.38727151913369e-05\\
64	1.34087238421955e-05\\
65	1.23474197328915e-05\\
66	1.09207647126036e-05\\
67	9.73111838821385e-06\\
68	8.74000648000066e-06\\
69	7.98080501643417e-06\\
70	7.46807604856168e-06\\
71	6.60890846491079e-06\\
72	6.05935448450636e-06\\
73	5.74654981522535e-06\\
74	5.46518633216974e-06\\
75	4.95058510487894e-06\\
76	4.33511338104887e-06\\
77	3.80079533483429e-06\\
78	3.31074844352472e-06\\
79	3.18449530295393e-06\\
80	3.08045328029361e-06\\
81	2.81861493061016e-06\\
82	2.61181038409031e-06\\
83	2.34462667221395e-06\\
84	2.18700706985714e-06\\
85	1.91167105154531e-06\\
86	1.75503880055508e-06\\
87	1.64539020110128e-06\\
88	1.63445016830087e-06\\
89	1.49591712543722e-06\\
90	1.398024407424e-06\\
91	1.25508792714958e-06\\
92	1.18388574572144e-06\\
93	1.14199071472899e-06\\
94	1.07438393037058e-06\\
95	9.98983641254546e-07\\
96	9.08286839973214e-07\\
97	8.67114975332366e-07\\
98	7.65754665369489e-07\\
99	6.86472074100301e-07\\
100	6.6957433272037e-07\\
101	6.45597915136069e-07\\
102	5.90860358653478e-07\\
103	5.20562691152307e-07\\
104	4.9050557494193e-07\\
105	4.52771649947493e-07\\
106	4.21501592091129e-07\\
107	3.98473247553185e-07\\
108	3.30946831581498e-07\\
109	3.14250686450415e-07\\
110	2.94451454443168e-07\\
111	2.20649258256232e-07\\
112	2.05742251935463e-07\\
113	1.72064912273951e-07\\
114	1.53361213667657e-07\\
115	1.30241744625303e-07\\
116	1.20323202523983e-07\\
117	1.14653721883292e-07\\
118	1.05125927957809e-07\\
119	8.93467828772563e-08\\
120	7.87793880727286e-08\\
121	6.99737882131014e-08\\
122	6.31251661875894e-08\\
123	6.16833867458991e-08\\
124	5.10799289977206e-08\\
125	4.8331317616829e-08\\
126	4.61732140823476e-08\\
127	4.19463565527835e-08\\
128	3.3469319871294e-08\\
129	3.02160683859947e-08\\
130	2.88045081200309e-08\\
131	2.13582777745297e-08\\
132	2.06541064722318e-08\\
133	1.60826052956386e-08\\
134	1.49512527010626e-08\\
135	1.25241896045188e-08\\
136	1.11031358468476e-08\\
137	9.91636648341412e-09\\
138	8.96155399136269e-09\\
139	8.33729537362478e-09\\
140	7.59849841820696e-09\\
141	6.38897483613771e-09\\
142	5.58749420058907e-09\\
143	5.42864255546898e-09\\
144	4.75289761990671e-09\\
145	4.72140541575106e-09\\
146	4.33879642204696e-09\\
147	4.18214284964412e-09\\
148	3.49905622366273e-09\\
149	2.78489689699663e-09\\
150	2.63876022284108e-09\\
151	2.19697941714575e-09\\
152	1.78899303252251e-09\\
153	1.53320114449002e-09\\
154	1.46660540205249e-09\\
155	1.30456299136624e-09\\
156	1.21847519581125e-09\\
157	1.03132118898206e-09\\
158	9.78417282610371e-10\\
159	8.51768226941747e-10\\
160	7.06926231279197e-10\\
161	6.80755175325687e-10\\
162	5.89335859056942e-10\\
163	5.26151346177802e-10\\
164	4.49682196965331e-10\\
165	3.77915948393852e-10\\
166	3.27382286961665e-10\\
167	2.74335574140244e-10\\
168	2.57120148110347e-10\\
169	2.36894010411019e-10\\
170	1.98580340488852e-10\\
171	1.69621846683797e-10\\
172	1.61860277926366e-10\\
173	1.37214126760352e-10\\
174	1.25665272065343e-10\\
175	1.05727410243858e-10\\
176	1.01583627777776e-10\\
177	7.33709200436226e-11\\
178	6.74083099215106e-11\\
179	5.88161745526511e-11\\
180	5.47421081157629e-11\\
181	5.07347992564445e-11\\
182	4.33040016125631e-11\\
183	3.83707626113574e-11\\
184	3.12348425317245e-11\\
185	2.86264623670156e-11\\
186	2.54680990474844e-11\\
187	2.27782153981778e-11\\
188	1.90484152656464e-11\\
189	1.55765977855488e-11\\
190	1.17497662672927e-11\\
191	1.11507861396613e-11\\
192	9.23651016746346e-12\\
193	8.65347535063068e-12\\
194	7.6298130735066e-12\\
195	7.04575538728807e-12\\
196	6.60840399139732e-12\\
197	6.20977511671388e-12\\
198	5.75238783499149e-12\\
199	5.60374364132724e-12\\
200	5.35797051018272e-12\\
201	5.13976080320241e-12\\
202	4.98918727554102e-12\\
203	4.70973186096215e-12\\
204	4.39358867901492e-12\\
205	4.34115471099475e-12\\
206	4.16771979458869e-12\\
207	3.92454879980369e-12\\
208	3.86202738404203e-12\\
209	3.67454861231486e-12\\
210	3.54294807072096e-12\\
211	3.4907825969965e-12\\
212	3.42534269077898e-12\\
213	3.18468467785964e-12\\
214	3.15949499764621e-12\\
215	3.03991189554762e-12\\
216	2.98629309144585e-12\\
217	2.74901656161523e-12\\
218	2.70141475782004e-12\\
219	2.62986991202277e-12\\
220	2.510536220512e-12\\
221	2.49360265358814e-12\\
222	2.43250593440931e-12\\
223	2.35078180238262e-12\\
224	2.29616796430195e-12\\
225	2.25347518831143e-12\\
226	2.23587750238927e-12\\
227	2.1775251748572e-12\\
228	2.06557830877017e-12\\
229	2.01559444719602e-12\\
230	1.96563359505928e-12\\
231	1.9097393098164e-12\\
232	1.89592805118648e-12\\
233	1.8198680563348e-12\\
234	1.79276579076155e-12\\
235	1.75975321733814e-12\\
236	1.69670410651583e-12\\
237	1.66986235507783e-12\\
238	1.6275207266866e-12\\
239	1.57104425981969e-12\\
240	1.54923356460098e-12\\
241	1.50735233996179e-12\\
242	1.49602496801235e-12\\
243	1.43258377750171e-12\\
244	1.40225468112816e-12\\
245	1.39393088903812e-12\\
246	1.38819697333085e-12\\
247	1.34894382680372e-12\\
248	1.30136130211e-12\\
249	1.29005287237621e-12\\
250	1.25694348496414e-12\\
251	1.23577914821538e-12\\
252	1.19770214601689e-12\\
253	1.16986829010124e-12\\
254	1.13610467994187e-12\\
255	1.11713911833995e-12\\
256	1.10164180628115e-12\\
257	1.07650964101442e-12\\
258	1.05812117792395e-12\\
259	1.03605752850912e-12\\
260	1.00845386605268e-12\\
261	1.00097498928801e-12\\
262	9.73325318041723e-13\\
263	9.54603934442593e-13\\
264	9.43618929567182e-13\\
265	9.26338521911166e-13\\
266	9.08427852378769e-13\\
267	8.86939392518466e-13\\
268	8.74472635656689e-13\\
269	8.45033930964385e-13\\
270	8.37655718596848e-13\\
271	8.24767693294043e-13\\
272	8.15519630165e-13\\
273	7.69899627279325e-13\\
274	7.65292742048292e-13\\
275	7.36007221361868e-13\\
276	7.30010185052914e-13\\
277	7.28450582376535e-13\\
278	7.22385690101485e-13\\
279	7.08587266409296e-13\\
280	6.97648332453383e-13\\
281	6.80840231482421e-13\\
282	6.73166357308373e-13\\
283	6.57776507391203e-13\\
284	6.46556705630074e-13\\
285	6.40607498572094e-13\\
286	6.24651235227115e-13\\
287	6.18305373574906e-13\\
288	6.11060025047812e-13\\
289	5.92198370542494e-13\\
290	5.88930284124814e-13\\
291	5.72990490517443e-13\\
292	5.62365791792987e-13\\
293	5.57234352243362e-13\\
294	5.42291559671223e-13\\
295	5.34571080278209e-13\\
296	5.23025372590778e-13\\
297	5.14274971404183e-13\\
298	5.07565308791795e-13\\
299	4.99461216469546e-13\\
300	4.98782858203235e-13\\
301	4.88430523815452e-13\\
302	4.86284031494737e-13\\
303	4.76123480183794e-13\\
304	4.70482710462595e-13\\
305	4.63534403094682e-13\\
306	4.59670366454174e-13\\
307	4.48447207406145e-13\\
308	4.42947850842402e-13\\
309	4.29634518128385e-13\\
310	4.2462974674483e-13\\
311	4.23174395365503e-13\\
312	4.14265499427121e-13\\
313	4.09931027191373e-13\\
314	4.02785394194091e-13\\
315	3.95450499153582e-13\\
316	3.88448820893365e-13\\
317	3.76670457903908e-13\\
318	3.72140202793712e-13\\
319	3.64111423595699e-13\\
320	3.58962041022022e-13\\
321	3.52997661215218e-13\\
322	3.5006663544649e-13\\
323	3.38609529382984e-13\\
324	3.31893523236229e-13\\
325	3.30776580102096e-13\\
326	3.28892750416247e-13\\
327	3.21731466312713e-13\\
328	3.13109985411965e-13\\
329	3.06452309486171e-13\\
330	3.0447588898952e-13\\
331	3.00799443072747e-13\\
332	2.95386790269155e-13\\
333	2.93631365188645e-13\\
334	2.85563210470126e-13\\
335	2.85083775430325e-13\\
336	2.81825984353878e-13\\
337	2.76783441473453e-13\\
338	2.72426355869106e-13\\
339	2.69213718514653e-13\\
340	2.68204648534015e-13\\
341	2.61160784118822e-13\\
342	2.56021930104878e-13\\
343	2.55519034060834e-13\\
344	2.47343485213727e-13\\
345	2.45917846282904e-13\\
346	2.3788429018084e-13\\
347	2.3638488145801e-13\\
348	2.34299507142354e-13\\
349	2.31788472369547e-13\\
350	2.28094640716624e-13\\
351	2.27213997928152e-13\\
352	2.23073765126185e-13\\
353	2.20873419986822e-13\\
354	2.18908094651199e-13\\
355	2.16284344055929e-13\\
356	2.15075275920519e-13\\
357	2.10562513456209e-13\\
358	2.07344754836968e-13\\
359	2.04755979687809e-13\\
360	2.02889202691854e-13\\
361	2.00345869805028e-13\\
362	1.96143945287641e-13\\
363	1.93069073530383e-13\\
364	1.89775326617694e-13\\
365	1.87599821069014e-13\\
366	1.86667220737986e-13\\
367	1.83096396032322e-13\\
368	1.77376710958245e-13\\
369	1.7611817495932e-13\\
370	1.75210617510892e-13\\
371	1.7229787814164e-13\\
372	1.70547820372701e-13\\
373	1.65614870110834e-13\\
374	1.64955122470787e-13\\
375	1.60168001746652e-13\\
376	1.58316160526206e-13\\
377	1.58121675367988e-13\\
378	1.55707452782389e-13\\
379	1.53374275033215e-13\\
380	1.5069771702902e-13\\
381	1.48867684323951e-13\\
382	1.46561906989923e-13\\
383	1.46084673335226e-13\\
384	1.43827449883708e-13\\
385	1.42367050798513e-13\\
386	1.41214581328477e-13\\
387	1.38090083294681e-13\\
388	1.368027925719e-13\\
389	1.34792947844996e-13\\
390	1.29838478004679e-13\\
391	1.28755119444101e-13\\
392	1.27408703111879e-13\\
393	1.26317129731859e-13\\
394	1.24001292337963e-13\\
395	1.22186221510934e-13\\
396	1.21756750941872e-13\\
397	1.19200830853702e-13\\
398	1.1819327072675e-13\\
399	1.16798552472448e-13\\
400	1.14654095777006e-13\\
401	1.14174818398932e-13\\
402	1.12784703308515e-13\\
403	1.10653589576885e-13\\
404	1.08738350584181e-13\\
405	1.07299435128123e-13\\
406	1.05539256498076e-13\\
407	1.04612401911334e-13\\
408	1.03047147765578e-13\\
409	1.00103943917065e-13\\
410	9.94653219387121e-14\\
411	9.79029689119392e-14\\
412	9.72810098690872e-14\\
413	9.58386218262595e-14\\
414	9.46595404209341e-14\\
415	9.39743163349447e-14\\
416	9.18259294877534e-14\\
417	9.01616251429815e-14\\
418	8.96139931622013e-14\\
419	8.70959812228858e-14\\
420	8.63548448459545e-14\\
421	8.51426414695881e-14\\
422	8.50044969983282e-14\\
423	8.44581178712887e-14\\
424	8.21193639375538e-14\\
425	8.15551667922655e-14\\
426	8.14187922017663e-14\\
427	7.92408356797145e-14\\
428	7.79921675910448e-14\\
429	7.6926352220327e-14\\
430	7.54238380747295e-14\\
431	7.46086462391255e-14\\
432	7.39077158856274e-14\\
433	7.35539241758579e-14\\
434	7.17326179944084e-14\\
435	7.0736796351822e-14\\
436	6.99548876600073e-14\\
437	6.79082411363376e-14\\
438	6.77425850438754e-14\\
439	6.73375421829927e-14\\
440	6.61256884401758e-14\\
441	6.59490538648995e-14\\
442	6.41714266222718e-14\\
443	6.35001273579667e-14\\
444	6.31836307086431e-14\\
445	6.1481308910654e-14\\
446	6.06685461561511e-14\\
447	6.0494974727812e-14\\
448	5.94134867599083e-14\\
449	5.91224479301898e-14\\
450	5.72949663155623e-14\\
451	5.68108557053636e-14\\
452	5.55117475955296e-14\\
453	5.4576615771245e-14\\
454	5.45601322309148e-14\\
455	5.2767860665399e-14\\
456	5.17834186143138e-14\\
457	5.1383908214602e-14\\
458	5.08760014305158e-14\\
459	5.00509671533791e-14\\
460	4.93983978019267e-14\\
461	4.871559788185e-14\\
462	4.8661868805137e-14\\
463	4.82642834681934e-14\\
464	4.75647505256348e-14\\
465	4.55042831345592e-14\\
466	4.52689081133511e-14\\
467	4.4502162895065e-14\\
468	4.41407479383546e-14\\
469	4.27090764891092e-14\\
470	4.22366663045862e-14\\
471	4.18568034831152e-14\\
472	4.02462122294185e-14\\
473	4.00362247772816e-14\\
474	3.84936136518624e-14\\
475	3.78804155002259e-14\\
476	3.78647716803049e-14\\
477	3.75607376224552e-14\\
478	3.6912679622307e-14\\
479	3.5611614417727e-14\\
480	3.54879257901668e-14\\
481	3.51436949444966e-14\\
482	3.43790861666598e-14\\
483	3.3774369365008e-14\\
484	3.31394650705703e-14\\
485	3.18114590248882e-14\\
486	3.12511980316994e-14\\
487	3.11314643435325e-14\\
488	3.05351616879133e-14\\
489	2.97537936699694e-14\\
490	2.93824830020461e-14\\
491	2.9001156526301e-14\\
492	2.74670871105068e-14\\
493	2.72319914039578e-14\\
494	2.57121209886142e-14\\
495	2.56742959331076e-14\\
496	2.55342115907784e-14\\
497	2.41744369804471e-14\\
498	2.33201853361566e-14\\
499	2.25549445315224e-14\\
500	2.09708371841434e-14\\
501	9.43517007158469e-17\\
};

\addplot [black, only marks, mark=*, mark options={fill=black}, forget plot]
table[row sep=crcr]{%
110	2.94451454443168e-07\\
};
\end{axis}
\end{tikzpicture}%

%% file: fig/tikz/p_x.tex
\definecolor{mycolor1}{rgb}{0.00000,0.44700,0.74100}%
\definecolor{mycolor2}{rgb}{0.92900,0.69400,0.12500}%
\definecolor{mycolor3}{rgb}{0.4940 0.1840 0.5560}%
\begin{tikzpicture}

\begin{axis}[%
width=5.4in,
height=1in,
at={(0.772in,0.473in)},
scale only axis,
xmin=0,
xmax=1,
xlabel style={font=\color{white!15!black}, yshift = 5pt},
xlabel={$t$},
ymin=-1e-05,
ymax=1e-05,
ylabel style={font=\color{white!15!black},yshift=3pt},
ylabel={$\mathrm{x}_{out} (t)$},
%axis background/.style={fill=white},
xmajorgrids,
ymajorgrids,
legend style={legend cell align=left, align=left, draw=white!15!black}
]
\addplot [color=mycolor1, line width=1.2pt]
  table[row sep=crcr]{%
0	0\\
0.001	8.73631141906023e-08\\
0.002	3.42755726924542e-07\\
0.003	7.11446099230902e-07\\
0.004	1.13347375055768e-06\\
0.005	1.58034486037069e-06\\
0.006	2.03157321486548e-06\\
0.007	2.47905563796971e-06\\
0.008	2.9182471932098e-06\\
0.009	3.34479914268689e-06\\
0.01	3.75547409082163e-06\\
0.011	4.14726485058565e-06\\
0.012	4.52050809166011e-06\\
0.013	4.8779162117043e-06\\
0.014	5.2201846088241e-06\\
0.015	5.54900539829328e-06\\
0.016	5.86486558611591e-06\\
0.017	6.16717785703916e-06\\
0.018	6.45362662296823e-06\\
0.019	6.7210552169162e-06\\
0.02	6.96633470252419e-06\\
0.021	7.18586477011675e-06\\
0.022	7.3740399425379e-06\\
0.023	7.52570591052723e-06\\
0.024	7.63958447403166e-06\\
0.025	7.71787566873387e-06\\
0.026	7.76297905587643e-06\\
0.027	7.77522201703884e-06\\
0.028	7.75277099430104e-06\\
0.029	7.69386408051573e-06\\
0.03	7.59867294746579e-06\\
0.031	7.46951947788475e-06\\
0.032	7.31071971122607e-06\\
0.033	7.12874423100196e-06\\
0.034	6.92975716267138e-06\\
0.035	6.71518598008844e-06\\
0.036	6.48022786619259e-06\\
0.037	6.21784871988824e-06\\
0.038	5.92442860317214e-06\\
0.039	5.6025867017129e-06\\
0.04	5.25979431355554e-06\\
0.041	4.90438555043755e-06\\
0.042	4.54125705867541e-06\\
0.043	4.16964878187255e-06\\
0.044	3.78371523226118e-06\\
0.045	3.37513075577719e-06\\
0.046	2.93638195280231e-06\\
0.047	2.46376295662244e-06\\
0.048	1.95924874435102e-06\\
0.049	1.4311116987317e-06\\
0.05	8.92952897122485e-07\\
0.051	3.60689356296373e-07\\
0.052	-1.52294130172779e-07\\
0.053	-6.3983374434473e-07\\
0.054	-1.10462830900824e-06\\
0.055	-1.55502073202551e-06\\
0.056	-1.99879593675932e-06\\
0.057	-2.4386639333413e-06\\
0.058	-2.87265615174145e-06\\
0.059	-3.29777585414265e-06\\
0.06	-3.7127676256596e-06\\
0.061	-4.11788205252944e-06\\
0.062	-4.51295647507458e-06\\
0.063	-4.89592255279782e-06\\
0.064	-5.26268680210359e-06\\
0.065	-5.60810823020196e-06\\
0.066	-5.92752324631245e-06\\
0.067	-6.21811808388991e-06\\
0.068	-6.47953197657667e-06\\
0.069	-6.71337772536108e-06\\
0.07	-6.92199821030349e-06\\
0.071	-7.10704503131857e-06\\
0.072	-7.26852238338108e-06\\
0.073	-7.40473123300529e-06\\
0.074	-7.51326280767995e-06\\
0.075	-7.59249296015144e-06\\
0.076	-7.64247151297061e-06\\
0.077	-7.6642343132762e-06\\
0.078	-7.65782312189246e-06\\
0.079	-7.62055334576063e-06\\
0.08	-7.54715882811332e-06\\
0.081	-7.4319876214263e-06\\
0.082	-7.27191049509854e-06\\
0.083	-7.06821260433892e-06\\
0.084	-6.82658995074646e-06\\
0.085	-6.55535903191434e-06\\
0.086	-6.26269985503242e-06\\
0.087	-5.95406043043822e-06\\
0.088	-5.63084722772793e-06\\
0.089	-5.29088522492315e-06\\
0.09	-4.93020878529138e-06\\
0.091	-4.5450788666514e-06\\
0.092	-4.13331101771518e-06\\
0.093	-3.69469831334903e-06\\
0.094	-3.23095379388445e-06\\
0.095	-2.74558661680445e-06\\
0.096	-2.24377534451928e-06\\
0.097	-1.73204404316062e-06\\
0.098	-1.21768619410191e-06\\
0.099	-7.08025240306481e-07\\
0.1	-2.09619403986046e-07\\
0.101	2.72541146696844e-07\\
0.102	7.35649286517166e-07\\
0.103	1.17920340641535e-06\\
0.104	1.60445134737455e-06\\
0.105	2.01337140741326e-06\\
0.106	2.40764707742467e-06\\
0.107	2.78802620260328e-06\\
0.108	3.15414463879667e-06\\
0.109	3.50471644927555e-06\\
0.11	3.83796260474987e-06\\
0.111	4.15220688103081e-06\\
0.112	4.44648982753439e-06\\
0.113	4.72095087585181e-06\\
0.114	4.9767105283067e-06\\
0.115	5.21523001835408e-06\\
0.116	5.43744222252473e-06\\
0.117	5.64313520089185e-06\\
0.118	5.83090813697774e-06\\
0.119	5.99868019642326e-06\\
0.12	6.14440844181928e-06\\
0.121	6.26661866658633e-06\\
0.122	6.36452668004777e-06\\
0.123	6.43784132446649e-06\\
0.124	6.48652775632788e-06\\
0.125	6.51078989241244e-06\\
0.126	6.51128150774054e-06\\
0.127	6.48931878449149e-06\\
0.128	6.44678321101788e-06\\
0.129	6.38557572504027e-06\\
0.13	6.30674823118483e-06\\
0.131	6.20968089193723e-06\\
0.132	6.09172235211688e-06\\
0.133	5.94857546159052e-06\\
0.134	5.7753943813056e-06\\
0.135	5.56824265139449e-06\\
0.136	5.3253718973228e-06\\
0.137	5.04786555434189e-06\\
0.138	4.73946281640653e-06\\
0.139	4.40571709956377e-06\\
0.14	4.05284292460921e-06\\
0.141	3.68663593323306e-06\\
0.142	3.31171100146999e-06\\
0.143	2.93113310938599e-06\\
0.144	2.54637001019322e-06\\
0.145	2.15746057289763e-06\\
0.146	1.76331807243915e-06\\
0.147	1.36215819667116e-06\\
0.148	9.52060401621975e-07\\
0.149	5.31643347615537e-07\\
0.15	1.00741635262014e-07\\
0.151	-3.39096530114819e-07\\
0.152	-7.84483828644769e-07\\
0.153	-1.23052900145709e-06\\
0.154	-1.67156369693204e-06\\
0.155	-2.10203008042822e-06\\
0.156	-2.5172500291696e-06\\
0.157	-2.9138640293977e-06\\
0.158	-3.28989514828428e-06\\
0.159	-3.64453318209548e-06\\
0.16	-3.97779807658801e-06\\
0.161	-4.29019457767658e-06\\
0.162	-4.58240083049412e-06\\
0.163	-4.85498290767348e-06\\
0.164	-5.10814611151061e-06\\
0.165	-5.34156285614403e-06\\
0.166	-5.55433423111726e-06\\
0.167	-5.74510108240943e-06\\
0.168	-5.91227014771364e-06\\
0.169	-6.05427689532696e-06\\
0.17	-6.1698166475079e-06\\
0.171	-6.25800141054113e-06\\
0.172	-6.31843725537718e-06\\
0.173	-6.35122558481614e-06\\
0.174	-6.35690182977585e-06\\
0.175	-6.33632662915597e-06\\
0.176	-6.29056081996998e-06\\
0.177	-6.22075066820445e-06\\
0.178	-6.12803838323097e-06\\
0.179	-6.01348343235793e-06\\
0.18	-5.87797616340454e-06\\
0.181	-5.7221367909726e-06\\
0.182	-5.54623054803884e-06\\
0.183	-5.35014794200007e-06\\
0.184	-5.13349915513395e-06\\
0.185	-4.89583230059066e-06\\
0.186	-4.63694139620258e-06\\
0.187	-4.35718170189168e-06\\
0.188	-4.05769445334288e-06\\
0.189	-3.74045190408773e-06\\
0.19	-3.40808742460342e-06\\
0.191	-3.0635424101071e-06\\
0.192	-2.70963356536848e-06\\
0.193	-2.34867221699656e-06\\
0.194	-1.98225243833287e-06\\
0.195	-1.61125595930698e-06\\
0.196	-1.23604516121739e-06\\
0.197	-8.56749347697663e-07\\
0.198	-4.73536055524429e-07\\
0.199	-8.67840827919664e-08\\
0.2	3.02864801574603e-07\\
0.201	6.94549205227965e-07\\
0.202	1.08726910444396e-06\\
0.203	1.47995522903621e-06\\
0.204	1.8714966927221e-06\\
0.205	2.260723851816e-06\\
0.206	2.64635225974613e-06\\
0.207	3.0269017758436e-06\\
0.208	3.40060683182427e-06\\
0.209	3.7653493983899e-06\\
0.21	4.11865403567556e-06\\
0.211	4.45778465119347e-06\\
0.212	4.77995109833813e-06\\
0.213	5.08259054712135e-06\\
0.214	5.36364169175188e-06\\
0.215	5.62171907012679e-06\\
0.216	5.85611931338696e-06\\
0.217	6.06665453151824e-06\\
0.218	6.25336889786528e-06\\
0.219	6.4162359831194e-06\\
0.22	6.55492842778282e-06\\
0.221	6.66871754587239e-06\\
0.222	6.75650758878945e-06\\
0.223	6.81697353701785e-06\\
0.224	6.84875285201965e-06\\
0.225	6.85065129773374e-06\\
0.226	6.82183694607469e-06\\
0.227	6.76201317727044e-06\\
0.228	6.67156061681744e-06\\
0.229	6.55163146219124e-06\\
0.23	6.40416413266908e-06\\
0.231	6.23178497621776e-06\\
0.232	6.03757500589813e-06\\
0.233	5.82471558434253e-06\\
0.234	5.59606675230043e-06\\
0.235	5.3537689723085e-06\\
0.236	5.09896559044664e-06\\
0.237	4.83172238796729e-06\\
0.238	4.55116779937788e-06\\
0.239	4.25582052896668e-06\\
0.24	3.94402043502661e-06\\
0.241	3.61435883940868e-06\\
0.242	3.26601037368551e-06\\
0.243	2.89890478657754e-06\\
0.244	2.51372204515635e-06\\
0.245	2.11174339142875e-06\\
0.246	1.6946224510361e-06\\
0.247	1.26415529462341e-06\\
0.248	8.22114923788856e-07\\
0.249	3.70188401455907e-07\\
0.25	-8.99849694640802e-08\\
0.251	-5.56706591060949e-07\\
0.252	-1.02809531659459e-06\\
0.253	-1.50199805933938e-06\\
0.254	-1.97597019953521e-06\\
0.255	-2.4473344644559e-06\\
0.256	-2.91330385835475e-06\\
0.257	-3.37113357441573e-06\\
0.258	-3.81826136810953e-06\\
0.259	-4.25239826787187e-06\\
0.26	-4.67154873188804e-06\\
0.261	-5.07395913300868e-06\\
0.262	-5.45801786519921e-06\\
0.263	-5.82214395741431e-06\\
0.264	-6.16470600611847e-06\\
0.265	-6.48399958689181e-06\\
0.266	-6.77829090732367e-06\\
0.267	-7.04590811343819e-06\\
0.268	-7.28534705740111e-06\\
0.269	-7.49535422605324e-06\\
0.27	-7.67496352437147e-06\\
0.271	-7.82348162213078e-06\\
0.272	-7.9404365777291e-06\\
0.273	-8.02551123915184e-06\\
0.274	-8.07848159046445e-06\\
0.275	-8.09916733467389e-06\\
0.276	-8.08739232678374e-06\\
0.277	-8.04294571077128e-06\\
0.278	-7.96554058558972e-06\\
0.279	-7.85477601501765e-06\\
0.28	-7.7101208077081e-06\\
0.281	-7.53094016448542e-06\\
0.282	-7.31658229228718e-06\\
0.283	-7.06652567566892e-06\\
0.284	-6.78057065116143e-06\\
0.285	-6.45904124890758e-06\\
0.286	-6.1029571110374e-06\\
0.287	-5.71413661689948e-06\\
0.288	-5.29520667340038e-06\\
0.289	-4.84951190925754e-06\\
0.29	-4.38093778648767e-06\\
0.291	-3.89367737357098e-06\\
0.292	-3.39198288012108e-06\\
0.293	-2.87994230752167e-06\\
0.294	-2.36131482182039e-06\\
0.295	-1.83944164286983e-06\\
0.296	-1.31723240288433e-06\\
0.297	-7.97209023435324e-07\\
0.298	-2.81580618483056e-07\\
0.299	2.27679797734444e-07\\
0.3	7.28777744112773e-07\\
0.301	1.22006965822945e-06\\
0.302	1.7000627688983e-06\\
0.303	2.1674330269655e-06\\
0.304	2.62104699788674e-06\\
0.305	3.05997624872379e-06\\
0.306	3.48349512241945e-06\\
0.307	3.89105836897908e-06\\
0.308	4.28225731163761e-06\\
0.309	4.65675753495056e-06\\
0.31	5.01422270612811e-06\\
0.311	5.35423380048819e-06\\
0.312	5.67621494461716e-06\\
0.313	5.97938031887912e-06\\
0.314	6.26271435115154e-06\\
0.315	6.52499403626743e-06\\
0.316	6.7648527865416e-06\\
0.317	6.98087657374463e-06\\
0.318	7.17171296722428e-06\\
0.319	7.33617000461886e-06\\
0.32	7.47328183516346e-06\\
0.321	7.58232718673952e-06\\
0.322	7.66279846696831e-06\\
0.323	7.71433516780674e-06\\
0.324	7.73664599318084e-06\\
0.325	7.72945051961541e-06\\
0.326	7.69246616145287e-06\\
0.327	7.62545480769251e-06\\
0.328	7.52832417716543e-06\\
0.329	7.40126146656471e-06\\
0.33	7.24486224029768e-06\\
0.331	7.06021466860496e-06\\
0.332	6.84890542117829e-06\\
0.333	6.61293194715819e-06\\
0.334	6.35452791073668e-06\\
0.335	6.07593278051047e-06\\
0.336	5.77915249854633e-06\\
0.337	5.46576534542715e-06\\
0.338	5.13681900783749e-06\\
0.339	4.79284683872841e-06\\
0.34	4.43400340235038e-06\\
0.341	4.06029243174296e-06\\
0.342	3.67183725867402e-06\\
0.343	3.26913456982864e-06\\
0.344	2.85323553543889e-06\\
0.345	2.42581683977492e-06\\
0.346	1.98912905093496e-06\\
0.347	1.54583830046719e-06\\
0.348	1.09879849044788e-06\\
0.349	6.50804266217701e-07\\
0.35	2.04373988239427e-07\\
0.351	-2.38397743987538e-07\\
0.352	-6.75896987409243e-07\\
0.353	-1.10695454242469e-06\\
0.354	-1.53073370202755e-06\\
0.355	-1.94656994276114e-06\\
0.356	-2.35379811844518e-06\\
0.357	-2.75160150755935e-06\\
0.358	-3.13891181770068e-06\\
0.359	-3.51437754221352e-06\\
0.36	-3.87640543360054e-06\\
0.361	-4.22326493811796e-06\\
0.362	-4.55323445357469e-06\\
0.363	-4.86475996329427e-06\\
0.364	-5.15659609211354e-06\\
0.365	-5.42790352648505e-06\\
0.366	-5.67828728073078e-06\\
0.367	-5.90777148286235e-06\\
0.368	-6.1167188207848e-06\\
0.369	-6.30571057280061e-06\\
0.37	-6.47540825731949e-06\\
0.371	-6.62641677974918e-06\\
0.372	-6.7591659092568e-06\\
0.373	-6.87382024134131e-06\\
0.374	-6.97022265328314e-06\\
0.375	-7.04787073239361e-06\\
0.376	-7.10592382317491e-06\\
0.377	-7.1432368579833e-06\\
0.378	-7.15841855665292e-06\\
0.379	-7.14991152712981e-06\\
0.38	-7.11609242903094e-06\\
0.381	-7.05538783667561e-06\\
0.382	-6.96639904095576e-06\\
0.383	-6.8480244146125e-06\\
0.384	-6.699565633889e-06\\
0.385	-6.52080229670651e-06\\
0.386	-6.3120224250822e-06\\
0.387	-6.07400112536926e-06\\
0.388	-5.8079286554796e-06\\
0.389	-5.51529768264447e-06\\
0.39	-5.19776857110191e-06\\
0.391	-4.85703604587529e-06\\
0.392	-4.49472204201304e-06\\
0.393	-4.11231450450101e-06\\
0.394	-3.71116397838302e-06\\
0.395	-3.29253779599985e-06\\
0.396	-2.85772081294913e-06\\
0.397	-2.40814149351231e-06\\
0.398	-1.94549752570874e-06\\
0.399	-1.47185421391044e-06\\
0.4	-9.89694421431753e-07\\
0.401	-5.01906917433805e-07\\
0.402	-1.17116595231858e-08\\
0.403	4.77469054158315e-07\\
0.404	9.6217468007331e-07\\
0.405	1.43908167212917e-06\\
0.406	1.90516343646389e-06\\
0.407	2.35781201586135e-06\\
0.408	2.79490451583377e-06\\
0.409	3.21480734467123e-06\\
0.41	3.61632117174424e-06\\
0.411	3.99857993901138e-06\\
0.412	4.36092388565412e-06\\
0.413	4.70277065893261e-06\\
0.414	5.02350694861403e-06\\
0.415	5.32241853786854e-06\\
0.416	5.59866769905593e-06\\
0.417	5.85131774533742e-06\\
0.418	6.0793947054279e-06\\
0.419	6.28196989460463e-06\\
0.42	6.45824340886542e-06\\
0.421	6.60761031950792e-06\\
0.422	6.72969574853233e-06\\
0.423	6.8243533613559e-06\\
0.424	6.89162990060694e-06\\
0.425	6.93170669491654e-06\\
0.426	6.9448334426108e-06\\
0.427	6.93127131584588e-06\\
0.428	6.89125896447709e-06\\
0.429	6.82500947046568e-06\\
0.43	6.7327378436613e-06\\
0.431	6.61471151427754e-06\\
0.432	6.47130992716537e-06\\
0.433	6.30307743399001e-06\\
0.434	6.11075442108569e-06\\
0.435	5.89527698629458e-06\\
0.436	5.65774210911237e-06\\
0.437	5.39934384151905e-06\\
0.438	5.12129254170458e-06\\
0.439	4.82473426344157e-06\\
0.44	4.51068778762155e-06\\
0.441	4.18001434966034e-06\\
0.442	3.83342819064666e-06\\
0.443	3.47154817360113e-06\\
0.444	3.0949815030444e-06\\
0.445	2.70442419399435e-06\\
0.446	2.30075877913194e-06\\
0.447	1.88513081620261e-06\\
0.448	1.458989741082e-06\\
0.449	1.02408752893696e-06\\
0.45	5.82436597199841e-07\\
0.451	1.36236576158056e-07\\
0.452	-3.12215740923644e-07\\
0.453	-7.60616598533663e-07\\
0.454	-1.20673105711752e-06\\
0.455	-1.64844350589747e-06\\
0.456	-2.0837782388534e-06\\
0.457	-2.51089063493892e-06\\
0.458	-2.9280346422141e-06\\
0.459	-3.3335146954979e-06\\
0.46	-3.7256318049567e-06\\
0.461	-4.10263266348176e-06\\
0.462	-4.46266979103437e-06\\
0.463	-4.80377856327865e-06\\
0.464	-5.1238758236375e-06\\
0.465	-5.42078280993549e-06\\
0.466	-5.69227414317122e-06\\
0.467	-5.93615239389174e-06\\
0.468	-6.15034589631143e-06\\
0.469	-6.33302398130791e-06\\
0.47	-6.48272074251862e-06\\
0.471	-6.59845446357296e-06\\
0.472	-6.67982759183362e-06\\
0.473	-6.72709067278367e-06\\
0.474	-6.74115562643819e-06\\
0.475	-6.72354771683035e-06\\
0.476	-6.6762932347065e-06\\
0.477	-6.60174864298392e-06\\
0.478	-6.50238716269201e-06\\
0.479	-6.38056691325689e-06\\
0.48	-6.23831077836531e-06\\
0.481	-6.07712891970035e-06\\
0.482	-5.8979113760844e-06\\
0.483	-5.70090889087383e-06\\
0.484	-5.48580795590869e-06\\
0.485	-5.25189124710909e-06\\
0.486	-4.99826133593893e-06\\
0.487	-4.72409449846943e-06\\
0.488	-4.42888635163862e-06\\
0.489	-4.11265135957772e-06\\
0.49	-3.77604577669069e-06\\
0.491	-3.42039555680806e-06\\
0.492	-3.04762679754157e-06\\
0.493	-2.6601120069584e-06\\
0.494	-2.26045971675004e-06\\
0.495	-1.85128358776171e-06\\
0.496	-1.43499023795383e-06\\
0.497	-1.01362050171561e-06\\
0.498	-5.88769206337997e-07\\
0.499	-1.61594101160366e-07\\
0.5	2.67090488722824e-07\\
0.501	6.96654279305939e-07\\
0.502	1.12645498509419e-06\\
0.503	1.55564029578825e-06\\
0.504	1.98298150071945e-06\\
0.505	2.40676670402546e-06\\
0.506	2.82477095090016e-06\\
0.507	3.23430844255356e-06\\
0.508	3.63235883152708e-06\\
0.509	4.01574886278732e-06\\
0.51	4.3813623454808e-06\\
0.511	4.72634814329004e-06\\
0.512	5.04829637013895e-06\\
0.513	5.34535851455203e-06\\
0.514	5.61629541415988e-06\\
0.515	5.86044810794457e-06\\
0.516	6.07763754925889e-06\\
0.517	6.26800986581674e-06\\
0.518	6.43185134311619e-06\\
0.519	6.5694019458817e-06\\
0.52	6.68069587420774e-06\\
0.521	6.76545384730793e-06\\
0.522	6.8230436588984e-06\\
0.523	6.85251567939897e-06\\
0.524	6.8527084809989e-06\\
0.525	6.82240975489831e-06\\
0.526	6.76054926817684e-06\\
0.527	6.66639636882621e-06\\
0.528	6.53973367341352e-06\\
0.529	6.38098245312345e-06\\
0.53	6.19126182667484e-06\\
0.531	5.97237348115452e-06\\
0.532	5.72671347769626e-06\\
0.533	5.45712270606235e-06\\
0.534	5.16669503834498e-06\\
0.535	4.85856758119225e-06\\
0.536	4.53571869097044e-06\\
0.537	4.20079787168027e-06\\
0.538	3.85600632860703e-06\\
0.539	3.50303992712454e-06\\
0.54	3.1430973138311e-06\\
0.541	2.77694756050144e-06\\
0.542	2.40504376414064e-06\\
0.543	2.02766394763529e-06\\
0.544	1.64505774069979e-06\\
0.545	1.25757829116307e-06\\
0.546	8.65782161009456e-07\\
0.547	4.70486258996638e-07\\
0.548	7.27778470994631e-08\\
0.549	-3.26018563377054e-07\\
0.55	-7.24407525114467e-07\\
0.551	-1.12080759586116e-06\\
0.552	-1.51363895070028e-06\\
0.553	-1.90139851797541e-06\\
0.554	-2.28271141584712e-06\\
0.555	-2.65635248356473e-06\\
0.556	-3.02123802759291e-06\\
0.557	-3.37639316860212e-06\\
0.558	-3.72090463446255e-06\\
0.559	-4.05387056024947e-06\\
0.56	-4.37435879467181e-06\\
0.561	-4.68138234408629e-06\\
0.562	-4.97389672720038e-06\\
0.563	-5.25081882332142e-06\\
0.564	-5.51106248048169e-06\\
0.565	-5.75358239791293e-06\\
0.566	-5.97741645658438e-06\\
0.567	-6.1817168343714e-06\\
0.568	-6.36576297512811e-06\\
0.569	-6.52895308241056e-06\\
0.57	-6.67077567140477e-06\\
0.571	-6.79076673848435e-06\\
0.572	-6.88846158851228e-06\\
0.573	-6.96335147364648e-06\\
0.574	-7.01485479508723e-06\\
0.575	-7.04230966119397e-06\\
0.576	-7.04499076873944e-06\\
0.577	-7.02214857198619e-06\\
0.578	-6.97306464110048e-06\\
0.579	-6.89711376248215e-06\\
0.58	-6.79382247317224e-06\\
0.581	-6.66291445216859e-06\\
0.582	-6.50433639920346e-06\\
0.583	-6.31826195551462e-06\\
0.584	-6.10507606908631e-06\\
0.585	-5.86534589461341e-06\\
0.586	-5.59978711573284e-06\\
0.587	-5.30923488472261e-06\\
0.588	-4.99462739148063e-06\\
0.589	-4.65700662874093e-06\\
0.59	-4.29753697027818e-06\\
0.591	-3.91753750297785e-06\\
0.592	-3.51852059282512e-06\\
0.593	-3.10222657646087e-06\\
0.594	-2.67064427990092e-06\\
0.595	-2.22600831615797e-06\\
0.596	-1.7707676226958e-06\\
0.597	-1.30752390034537e-06\\
0.598	-8.38943830831726e-07\\
0.599	-3.67653342370592e-07\\
0.6	1.03873929690359e-07\\
0.601	5.73420101564443e-07\\
0.602	1.03910425322584e-06\\
0.603	1.49942949952427e-06\\
0.604	1.95328761244502e-06\\
0.605	2.39991833613942e-06\\
0.606	2.83882616344849e-06\\
0.607	3.2696633710107e-06\\
0.608	3.69209273153006e-06\\
0.609	4.10564682697993e-06\\
0.61	4.50960175425661e-06\\
0.611	4.90288198767551e-06\\
0.612	5.28400932389821e-06\\
0.613	5.65110357706978e-06\\
0.614	6.00193573571178e-06\\
0.615	6.33402752187444e-06\\
0.616	6.64478472690685e-06\\
0.617	6.93164724390262e-06\\
0.618	7.19223613326996e-06\\
0.619	7.42447870827896e-06\\
0.62	7.62669556788648e-06\\
0.621	7.79763928771271e-06\\
0.622	7.93648138287157e-06\\
0.623	8.04275209604032e-06\\
0.624	8.11624445418428e-06\\
0.625	8.15689972278455e-06\\
0.626	8.16469407302429e-06\\
0.627	8.13954639269024e-06\\
0.628	8.08126372039005e-06\\
0.629	7.98953512371771e-06\\
0.63	7.86397691261042e-06\\
0.631	7.70422393464017e-06\\
0.632	7.51005372403534e-06\\
0.633	7.28152462352791e-06\\
0.634	7.01910572023063e-06\\
0.635	6.72377698595921e-06\\
0.636	6.39708152664097e-06\\
0.637	6.04111871075826e-06\\
0.638	5.65847530097743e-06\\
0.639	5.25210116556597e-06\\
0.64	4.82514439707208e-06\\
0.641	4.38076727340447e-06\\
0.642	3.92196758353173e-06\\
0.643	3.45142973245225e-06\\
0.644	2.97142586311702e-06\\
0.645	2.48378042129444e-06\\
0.646	1.98990230080542e-06\\
0.647	1.49087914612319e-06\\
0.648	9.87619170805023e-07\\
0.649	4.81019173168747e-07\\
0.65	-2.78665262263854e-08\\
0.651	-5.37681316046262e-07\\
0.652	-1.04666730616794e-06\\
0.653	-1.55262668587881e-06\\
0.654	-2.05294925757632e-06\\
0.655	-2.54470207374388e-06\\
0.656	-3.02476834407734e-06\\
0.657	-3.49001567636915e-06\\
0.658	-3.93747010993448e-06\\
0.659	-4.36447199849638e-06\\
0.66	-4.76879337561283e-06\\
0.661	-5.14870257399628e-06\\
0.662	-5.50297033659045e-06\\
0.663	-5.8308203155388e-06\\
0.664	-6.13183513495471e-06\\
0.665	-6.40583507911124e-06\\
0.666	-6.65274989282838e-06\\
0.667	-6.87250396465123e-06\\
0.668	-7.06493223858265e-06\\
0.669	-7.22973837874239e-06\\
0.67	-7.36649974985088e-06\\
0.671	-7.4747160419188e-06\\
0.672	-7.55389180305053e-06\\
0.673	-7.6036379484276e-06\\
0.674	-7.62377492181132e-06\\
0.675	-7.61442009931423e-06\\
0.676	-7.57604487974534e-06\\
0.677	-7.50949142308527e-06\\
0.678	-7.41594518301858e-06\\
0.679	-7.29686557070077e-06\\
0.68	-7.15388313858107e-06\\
0.681	-6.98867605009534e-06\\
0.682	-6.8028415091223e-06\\
0.683	-6.59777822646214e-06\\
0.684	-6.37459465552452e-06\\
0.685	-6.13405417766742e-06\\
0.686	-5.87656391105034e-06\\
0.687	-5.60220824754712e-06\\
0.688	-5.31082302171553e-06\\
0.689	-5.00210132511539e-06\\
0.69	-4.6757186858206e-06\\
0.691	-4.33146327685803e-06\\
0.692	-3.96935690226152e-06\\
0.693	-3.58975395457895e-06\\
0.694	-3.19340878718954e-06\\
0.695	-2.78150582979073e-06\\
0.696	-2.35565146139203e-06\\
0.697	-1.91783083754934e-06\\
0.698	-1.47033671378225e-06\\
0.699	-1.01567967709032e-06\\
0.7	-5.56490550530547e-07\\
0.701	-9.54253231297026e-08\\
0.702	3.64918376129354e-07\\
0.703	8.22066882156964e-07\\
0.704	1.27371461188119e-06\\
0.705	1.71774635844648e-06\\
0.706	2.15223994344832e-06\\
0.707	2.57545358145505e-06\\
0.708	2.98580326046399e-06\\
0.709	3.38183582777341e-06\\
0.71	3.76220272686346e-06\\
0.711	4.12563835765372e-06\\
0.712	4.47094542297512e-06\\
0.713	4.796988323591e-06\\
0.714	5.10269421310213e-06\\
0.715	5.38706052691793e-06\\
0.716	5.64916703016387e-06\\
0.717	5.88819031885596e-06\\
0.718	6.10341850869508e-06\\
0.719	6.29426412710084e-06\\
0.72	6.46027325891291e-06\\
0.721	6.60112940916698e-06\\
0.722	6.71665070072186e-06\\
0.723	6.80677962631396e-06\\
0.724	6.87156502709007e-06\\
0.725	6.91113695298646e-06\\
0.726	6.92567590183386e-06\\
0.727	6.91537918521912e-06\\
0.728	6.8804280280184e-06\\
0.729	6.82095990631024e-06\\
0.73	6.73705071833364e-06\\
0.731	6.62871112599676e-06\\
0.732	6.49590008373908e-06\\
0.733	6.33855685885787e-06\\
0.734	6.15665028442073e-06\\
0.735	5.95024150937923e-06\\
0.736	5.71955388718664e-06\\
0.737	5.46504190467751e-06\\
0.738	5.18744995743553e-06\\
0.739	4.88785220721001e-06\\
0.74	4.56766618219145e-06\\
0.741	4.22863570800685e-06\\
0.742	3.87278226473766e-06\\
0.743	3.50232809377697e-06\\
0.744	3.11959815379349e-06\\
0.745	2.72691137182853e-06\\
0.746	2.32647344819908e-06\\
0.747	1.920284010987e-06\\
0.748	1.51006944882173e-06\\
0.749	1.09724999674551e-06\\
0.75	6.82945363924255e-07\\
0.751	2.68018503392109e-07\\
0.752	-1.46847908367706e-07\\
0.753	-5.61051831821469e-07\\
0.754	-9.73961814677207e-07\\
0.755	-1.38481054944657e-06\\
0.756	-1.79260576830611e-06\\
0.757	-2.19606983197942e-06\\
0.758	-2.59361579020873e-06\\
0.759	-2.98336282664435e-06\\
0.76	-3.36318906220304e-06\\
0.761	-3.73081479242416e-06\\
0.762	-4.08390558242906e-06\\
0.763	-4.42018215963981e-06\\
0.764	-4.73752359914584e-06\\
0.765	-5.03405138557337e-06\\
0.766	-5.30818493639998e-06\\
0.767	-5.55866311425784e-06\\
0.768	-5.78453113941435e-06\\
0.769	-5.98509693704022e-06\\
0.77	-6.1598652316367e-06\\
0.771	-6.30846052072622e-06\\
0.772	-6.43055157969667e-06\\
0.773	-6.52578963475686e-06\\
0.774	-6.59377039508917e-06\\
0.775	-6.63402648463742e-06\\
0.776	-6.64605247690253e-06\\
0.777	-6.62935977429349e-06\\
0.778	-6.58355418610659e-06\\
0.779	-6.50842532777429e-06\\
0.78	-6.40403494553627e-06\\
0.781	-6.27079075009242e-06\\
0.782	-6.10949394497609e-06\\
0.783	-5.92135166959511e-06\\
0.784	-5.70795018078924e-06\\
0.785	-5.47118964447845e-06\\
0.786	-5.21318672729331e-06\\
0.787	-4.93615552857629e-06\\
0.788	-4.64228066464242e-06\\
0.789	-4.33359751804335e-06\\
0.79	-4.01189410700627e-06\\
0.791	-3.67864626318757e-06\\
0.792	-3.33499369602563e-06\\
0.793	-2.98175912521066e-06\\
0.794	-2.61950722482187e-06\\
0.795	-2.24863490634381e-06\\
0.796	-1.86948068670356e-06\\
0.797	-1.48243855829445e-06\\
0.798	-1.08806169568591e-06\\
0.799	-6.87142982323139e-07\\
0.8	-2.80762939705454e-07\\
0.801	1.29699837408225e-07\\
0.802	5.42595367270967e-07\\
0.803	9.5606092464199e-07\\
0.804	1.36810233093804e-06\\
0.805	1.77668609720567e-06\\
0.806	2.17982960180471e-06\\
0.807	2.57567746567778e-06\\
0.808	2.96255458421537e-06\\
0.809	3.3389899756113e-06\\
0.81	3.70370968628724e-06\\
0.811	4.05560135662759e-06\\
0.812	4.39365664509901e-06\\
0.813	4.71690056938549e-06\\
0.814	5.02431816722087e-06\\
0.815	5.31478907403626e-06\\
0.816	5.58703924133618e-06\\
0.817	5.83961678849138e-06\\
0.818	6.0708957540225e-06\\
0.819	6.27910817756027e-06\\
0.82	6.46240145144871e-06\\
0.821	6.61891509350266e-06\\
0.822	6.7468688452193e-06\\
0.823	6.84465297684416e-06\\
0.824	6.91091150167448e-06\\
0.825	6.94461010265195e-06\\
0.826	6.94508234153906e-06\\
0.827	6.91205038547584e-06\\
0.828	6.84561927992367e-06\\
0.829	6.74624684810655e-06\\
0.83	6.61469380673764e-06\\
0.831	6.45196082408219e-06\\
0.832	6.25922032139687e-06\\
0.833	6.03775116043148e-06\\
0.834	5.78888350636267e-06\\
0.835	5.51395965744715e-06\\
0.836	5.21431423409426e-06\\
0.837	4.89127457686496e-06\\
0.838	4.54617939274572e-06\\
0.839	4.1804113929358e-06\\
0.84	3.79543774174308e-06\\
0.841	3.39285125186623e-06\\
0.842	2.97440508318853e-06\\
0.843	2.54203467567611e-06\\
0.844	2.09786221031826e-06\\
0.845	1.64418127003066e-06\\
0.846	1.18342183132106e-06\\
0.847	7.18098373505689e-07\\
0.848	2.50745997206819e-07\\
0.849	-2.16148830687125e-07\\
0.85	-6.80215571555425e-07\\
0.851	-1.13925757157049e-06\\
0.852	-1.59129435451552e-06\\
0.853	-2.03458352328686e-06\\
0.854	-2.46761992488409e-06\\
0.855	-2.88911221402961e-06\\
0.856	-3.29793965987484e-06\\
0.857	-3.69309425384377e-06\\
0.858	-4.07361505812745e-06\\
0.859	-4.43852264991082e-06\\
0.86	-4.78676175964137e-06\\
0.861	-5.11715932250413e-06\\
0.862	-5.42840364751665e-06\\
0.863	-5.71904800030515e-06\\
0.864	-5.98753926161223e-06\\
0.865	-6.23226934249258e-06\\
0.866	-6.4516444502873e-06\\
0.867	-6.64416501067219e-06\\
0.868	-6.80850773554378e-06\\
0.869	-6.94360075582089e-06\\
0.87	-7.04868340621636e-06\\
0.871	-7.12334370633307e-06\\
0.872	-7.16752908126114e-06\\
0.873	-7.18152873138444e-06\\
0.874	-7.16592936109815e-06\\
0.875	-7.12154893335438e-06\\
0.876	-7.04935576884852e-06\\
0.877	-6.95038193340885e-06\\
0.878	-6.82564064342577e-06\\
0.879	-6.67605687612194e-06\\
0.88	-6.5024189390664e-06\\
0.881	-6.30535621405221e-06\\
0.882	-6.08534533537149e-06\\
0.883	-5.8427436589645e-06\\
0.884	-5.57784582265063e-06\\
0.885	-5.29095647522665e-06\\
0.886	-4.9824705594394e-06\\
0.887	-4.652951679953e-06\\
0.888	-4.30319954286694e-06\\
0.889	-3.93429876084729e-06\\
0.89	-3.54764367670652e-06\\
0.891	-3.14493659458576e-06\\
0.892	-2.72815995093618e-06\\
0.893	-2.29952575710913e-06\\
0.894	-1.8614081567258e-06\\
0.895	-1.41626650097209e-06\\
0.896	-9.66567158735357e-07\\
0.897	-5.14711920053322e-07\\
0.898	-6.29797613392247e-08\\
0.899	3.86513299360614e-07\\
0.9	8.31833882006133e-07\\
0.901	1.27123081579392e-06\\
0.902	1.70311737316588e-06\\
0.903	2.12604110726685e-06\\
0.904	2.53864753195313e-06\\
0.905	2.93964451980985e-06\\
0.906	3.32777383079839e-06\\
0.907	3.70179503433896e-06\\
0.908	4.06048503777863e-06\\
0.909	4.40265409484994e-06\\
0.91	4.7271764680186e-06\\
0.911	5.03303156513518e-06\\
0.912	5.31934933447185e-06\\
0.913	5.58545262271546e-06\\
0.914	5.83088890069072e-06\\
0.915	6.05544469102887e-06\\
0.916	6.25913777829443e-06\\
0.917	6.44218501432772e-06\\
0.918	6.60494658791352e-06\\
0.919	6.74785098123233e-06\\
0.92	6.87130769659764e-06\\
0.921	6.9756171491366e-06\\
0.922	7.06088819233813e-06\\
0.923	7.12697370627966e-06\\
0.924	7.17343311254102e-06\\
0.925	7.19952806593211e-06\\
0.926	7.20425384938667e-06\\
0.927	7.18640493206075e-06\\
0.928	7.14466891271199e-06\\
0.929	7.07773955701085e-06\\
0.93	6.98443697860402e-06\\
0.931	6.86382192838128e-06\\
0.932	6.71529150023488e-06\\
0.933	6.53864566350576e-06\\
0.934	6.33411732833271e-06\\
0.935	6.10236308039387e-06\\
0.936	5.84441639947135e-06\\
0.937	5.56160983605986e-06\\
0.938	5.25547632450133e-06\\
0.939	4.9276424030384e-06\\
0.94	4.57972691529426e-06\\
0.941	4.2132579913045e-06\\
0.942	3.82961851143125e-06\\
0.943	3.43002650004648e-06\\
0.944	3.01555213687292e-06\\
0.945	2.58716825317869e-06\\
0.946	2.14582660730614e-06\\
0.947	1.69254883332807e-06\\
0.948	1.22851885423178e-06\\
0.949	7.55163353877031e-07\\
0.95	2.74208321887267e-07\\
0.951	-2.12297145046215e-07\\
0.952	-7.01994409047349e-07\\
0.953	-1.19226586355185e-06\\
0.954	-1.68031735498851e-06\\
0.955	-2.16328300518324e-06\\
0.956	-2.63834008575062e-06\\
0.957	-3.10282016166611e-06\\
0.958	-3.55430326212871e-06\\
0.959	-3.99068397140673e-06\\
0.96	-4.41020209130693e-06\\
0.961	-4.81143510060446e-06\\
0.962	-5.19325472145959e-06\\
0.963	-5.55475451090204e-06\\
0.964	-5.89515922162769e-06\\
0.965	-6.21372889896188e-06\\
0.966	-6.50967131078573e-06\\
0.967	-6.78207496483811e-06\\
0.968	-7.02987214657011e-06\\
0.969	-7.25183719470842e-06\\
0.97	-7.44662050279589e-06\\
0.971	-7.612813831338e-06\\
0.972	-7.74903840562595e-06\\
0.973	-7.85404420634812e-06\\
0.974	-7.92680752812757e-06\\
0.975	-7.96661415088354e-06\\
0.976	-7.97311757859697e-06\\
0.977	-7.94636513880595e-06\\
0.978	-7.88678912531078e-06\\
0.979	-7.7951646971079e-06\\
0.98	-7.67254059951828e-06\\
0.981	-7.52015207005557e-06\\
0.982	-7.3393273982597e-06\\
0.983	-7.1313999811657e-06\\
0.984	-6.89763663349727e-06\\
0.985	-6.63919027931111e-06\\
0.986	-6.357081662047e-06\\
0.987	-6.05221056496746e-06\\
0.988	-5.72539311617545e-06\\
0.989	-5.37741831155773e-06\\
0.99	-5.00911465138311e-06\\
0.991	-4.62141674514058e-06\\
0.992	-4.21542228040949e-06\\
0.993	-3.79243145830781e-06\\
0.994	-3.35396387398562e-06\\
0.995	-2.90175114283356e-06\\
0.996	-2.43770714482924e-06\\
0.997	-1.96388079448004e-06\\
0.998	-1.48239860374905e-06\\
0.999	-9.95405407730043e-07\\
1	-5.0501163416396e-07\\
};
\addlegendentry{\texttt{FOM}}

\addplot [color=mycolor2, line width=1.2pt, only marks, mark=triangle, mark options={solid, mycolor2}]
  table[row sep=crcr]{%
0	0\\
0.02	6.96633649111858e-06\\
0.04	5.25979314808078e-06\\
0.06	-3.7127681416782e-06\\
0.08	-7.54715906615857e-06\\
0.1	-2.09619512100536e-07\\
0.12	6.14440839676789e-06\\
0.14	4.05284295320072e-06\\
0.16	-3.97779809115877e-06\\
0.18	-5.87797617398232e-06\\
0.2	3.02864857451924e-07\\
0.22	6.55492842775137e-06\\
0.24	3.94402043389634e-06\\
0.26	-4.67154868186574e-06\\
0.28	-7.71012084139998e-06\\
0.3	7.28777753832037e-07\\
0.32	7.47328184033809e-06\\
0.34	4.43400336298117e-06\\
0.36	-3.87640546627763e-06\\
0.38	-7.11609246555727e-06\\
0.4	-9.89694399783804e-07\\
0.42	6.45824342060687e-06\\
0.44	4.51068782462745e-06\\
0.46	-3.72563176032402e-06\\
0.48	-6.23831078010235e-06\\
0.5	2.67090496053796e-07\\
0.52	6.68069587849558e-06\\
0.54	3.14309738096966e-06\\
0.56	-4.37435869955038e-06\\
0.58	-6.79382235270881e-06\\
0.6	1.03874038304545e-07\\
0.62	7.62669559316282e-06\\
0.64	4.82514436072429e-06\\
0.66	-4.7687934503507e-06\\
0.68	-7.15388324358737e-06\\
0.7	-5.564905889487e-07\\
0.72	6.46027322518951e-06\\
0.74	4.56766616484616e-06\\
0.76	-3.36318906547015e-06\\
0.78	-6.40403495543159e-06\\
0.8	-2.8076291644173e-07\\
0.82	6.46240149559919e-06\\
0.84	3.79543783240324e-06\\
0.86	-4.78676165227613e-06\\
0.88	-6.50241886671919e-06\\
0.9	8.3183389695137e-07\\
0.92	6.87130762647662e-06\\
0.94	4.57972678912174e-06\\
0.96	-4.41020215622803e-06\\
0.98	-7.67254058814925e-06\\
1	-5.05011577981135e-07\\
};
\addlegendentry{\texttt{POD}}

\addplot [color=mycolor3, line width=1.2pt, only marks, mark=diamond, mark options={solid, mycolor3}]
  table[row sep=crcr]{%
0.004	1.13302093765536e-06\\
0.014	5.22288546407969e-06\\
0.024	7.63374009457883e-06\\
0.034	6.93791573709061e-06\\
0.044	3.782786860726e-06\\
0.054	-1.11087826493024e-06\\
0.064	-5.26258788325982e-06\\
0.074	-7.50719300047004e-06\\
0.084	-6.82568917467021e-06\\
0.094	-3.22088604987535e-06\\
0.104	1.58958844461398e-06\\
0.114	4.98327403322285e-06\\
0.124	6.49312813388649e-06\\
0.134	5.77235369407325e-06\\
0.144	2.54035725895789e-06\\
0.154	-1.65424305006849e-06\\
0.164	-5.11240074412119e-06\\
0.174	-6.36075647379069e-06\\
0.184	-5.12374322608832e-06\\
0.194	-1.9911208039434e-06\\
0.204	1.87853507598538e-06\\
0.214	5.37226916125401e-06\\
0.224	6.83107431569173e-06\\
0.234	5.60961157533489e-06\\
0.244	2.51452525298569e-06\\
0.254	-1.98409649642212e-06\\
0.264	-6.14996757810794e-06\\
0.274	-8.09739349622943e-06\\
0.284	-6.76460831384415e-06\\
0.294	-2.36090395957609e-06\\
0.304	2.60594243513878e-06\\
0.314	6.27160340642839e-06\\
0.324	7.73013018889953e-06\\
0.334	6.35753857021553e-06\\
0.344	2.85325780985171e-06\\
0.354	-1.53321176901449e-06\\
0.364	-5.14297691498292e-06\\
0.374	-6.99417305772554e-06\\
0.384	-6.69197116977943e-06\\
0.394	-3.69464096518136e-06\\
0.404	9.48955976020423e-07\\
0.414	5.0271878172748e-06\\
0.424	6.8917756037183e-06\\
0.434	6.1121830623378e-06\\
0.444	3.09971869776549e-06\\
0.454	-1.22861821300409e-06\\
0.464	-5.1076424470468e-06\\
0.474	-6.73179790536367e-06\\
0.484	-5.51130050397904e-06\\
0.494	-2.24670279675386e-06\\
0.504	1.97680861965575e-06\\
0.514	5.62616639746541e-06\\
0.524	6.83632628373523e-06\\
0.534	5.1814909308471e-06\\
0.544	1.65038080668591e-06\\
0.554	-2.2906059960225e-06\\
0.564	-5.501245712787e-06\\
0.574	-7.03415903475314e-06\\
0.584	-6.08075421413504e-06\\
0.594	-2.68763817644987e-06\\
0.604	1.94966069552931e-06\\
0.614	6.01812178890934e-06\\
0.624	8.1025381050442e-06\\
0.634	7.01313242359836e-06\\
0.644	2.97939887909952e-06\\
0.654	-2.06078065986398e-06\\
0.664	-6.12303494150892e-06\\
0.674	-7.60838416960795e-06\\
0.684	-6.40682755796444e-06\\
0.694	-3.16411542863946e-06\\
0.704	1.26703815627796e-06\\
0.714	5.08778893381341e-06\\
0.724	6.89868628870492e-06\\
0.734	6.13329704980114e-06\\
0.744	3.1368189478307e-06\\
0.754	-9.82114962214452e-07\\
0.764	-4.7570733248148e-06\\
0.774	-6.57560337946082e-06\\
0.784	-5.70674944983353e-06\\
0.794	-2.63629379729978e-06\\
0.804	1.3774065086764e-06\\
0.814	5.03749167451363e-06\\
0.824	6.89729226144395e-06\\
0.834	5.78973315457094e-06\\
0.844	2.10496772450408e-06\\
0.854	-2.46205933336002e-06\\
0.864	-5.98706849754269e-06\\
0.874	-7.17214944081975e-06\\
0.884	-5.58216262737545e-06\\
0.894	-1.84404119857908e-06\\
0.904	2.50924998229518e-06\\
0.914	5.82860543015909e-06\\
0.924	7.20326439769442e-06\\
0.934	6.31154934743075e-06\\
0.944	3.01726401870149e-06\\
0.954	-1.67608961608216e-06\\
0.964	-5.8983339463752e-06\\
0.974	-7.91738874194517e-06\\
0.984	-6.90008426374302e-06\\
0.994	-3.35444209076913e-06\\
};
\addlegendentry{\texttt{cOpInf}}

\addplot [color=black,line width=1.5pt, forget plot]
table[row sep=crcr]{%
	0.5	-1e-5\\
	0.5	1e-5\\
};

\end{axis}
\end{tikzpicture}%

%% file: MechOpInf.bbl
\begin{thebibliography}{10}

\bibitem{Alt15}
H.~Altenbach.
\newblock {\em Kontinuumsmechanik}.
\newblock Springer Vieweg Berlin, Heidelberg, 2015.
\newblock \href {https://doi.org/https://doi.org/10.1007/978-3-662-47070-1}
  {\path{doi:https://doi.org/10.1007/978-3-662-47070-1}}.

\bibitem{morAntSG01}
A.~C. Antoulas, D.~C. Sorensen, and S.~Gugercin.
\newblock A survey of model reduction methods for large-scale systems.
\newblock {\em Contemp. Math.}, 280:193--219, 2001.

\bibitem{supAumW22}
Q.~Aumann and S.~W.~R. Werner.
\newblock Code, data and results for numerical experiments in ``{S}tructured
  model order reduction for vibro-acoustic problems using interpolation and
  balancing methods'' (version 1.0), January 2022.
\newblock \href {https://doi.org/10.5281/zenodo.5836047}
  {\path{doi:10.5281/zenodo.5836047}}.

\bibitem{morAumW22}
Q.~Aumann and S.~W.~R. Werner.
\newblock Structured model order reduction for vibro-acoustic problems using
  interpolation and balancing methods.
\newblock {\em JSV}, 543:117363, 2023.
\newblock \href {https://doi.org/https://doi.org/10.1016/j.jsv.2022.117363}
  {\path{doi:https://doi.org/10.1016/j.jsv.2022.117363}}.

\bibitem{morBeaG05}
C.~A. Beattie and S.~Gugercin.
\newblock {K}rylov-based model reduction of second-order systems with
  proportional damping.
\newblock In {\em Proceedings of the 44th IEEE Conference on Decision and
  Control}, pages 2278--2283, December 2005.
\newblock \href {https://doi.org/10.1109/CDC.2005.1582501}
  {\path{doi:10.1109/CDC.2005.1582501}}.

\bibitem{morBetal}
P.~Benner, Y.~Filanova, D.~Karachalios, S.~M. Abdelhafez, J.~Przybilla, and
  S.~W.~R. Werner.
\newblock Mathematische {K}omplexit{\"a}tsreduktion: {M}odellreduktion
  dynamischer {S}ysteme.
\newblock {\em Mitt. DMV}, 29(4):198--204, 2021.
\newblock \href {https://doi.org/10.1515/dmvm-2021-0075}
  {\path{doi:10.1515/dmvm-2021-0075}}.

\bibitem{morBenGHetal22}
P.~Benner, P.~Goyal, J.~Heiland, and I.~Pontes~Duff.
\newblock Operator inference and physics-informed learning of low-dimensional
  models for incompressible flows.
\newblock {\em Electron. Trans. Numer. Anal.}, 56:28--51, 2022.
\newblock (Special Issue on Scientific Machine Learning).
\newblock URL: \url{https://epub.oeaw.ac.at/?arp=0x003d183f}.

\bibitem{morBenGKetal20}
P.~Benner, P.~Goyal, B.~Kramer, B.~Peherstorfer, and K.~Willcox.
\newblock Operator inference for non-intrusive model reduction of systems with
  non-polynomial nonlinear terms.
\newblock {\em Comp. Meth. Appl. Mech. Eng.}, 372:113433, 2020.
\newblock \href {https://doi.org/10.1016/j.cma.2020.113433}
  {\path{doi:10.1016/j.cma.2020.113433}}.

\bibitem{morBenGP19}
P.~Benner, P.~Goyal, and I.~Pontes~Duff.
\newblock Identification of dominant subspaces for linear structured parametric
  systems and model reduction.
\newblock e-prints 1910.13945, arXiv, 2019.
\newblock URL: \url{https://arxiv.org/abs/1910.13945}.

\bibitem{morPonGB20a}
P.~Benner, P.~Goyal, and I.~Pontes~Duff.
\newblock Data-driven identification of {R}ayleigh-damped second-order systems.
\newblock In C.~Beattie, P.~Benner, M.~Embree, S.~Gugercin, and S.~Lefteriu,
  editors, {\em Realization and Model Reduction of Dynamical Systems - A
  Festschrift in Honor of the 70th Birthday of {T}hanos {A}ntoulas}, pages
  255--272. Springer, Cham, 2022.
\newblock \href {https://doi.org/10.1007/978-3-030-95157-3_14}
  {\path{doi:10.1007/978-3-030-95157-3_14}}.

\bibitem{morBenMS05}
P.~Benner, V.~Mehrmann, and D.~C. Sorensen.
\newblock {\em Dimension Reduction of Large-Scale Systems}, volume~45 of {\em
  Lect. Notes Comput. Sci. Eng.}
\newblock Springer-Verlag, Berlin/Heidelberg, Germany, 2005.
\newblock \href {https://doi.org/10.1007/3-540-27909-1}
  {\path{doi:10.1007/3-540-27909-1}}.

\bibitem{morBerHL93}
G.~Berkooz, P.~Holmes, and J.~L. Lumley.
\newblock The proper orthogonal decomposition in the analysis of turbulent
  flows.
\newblock {\em Annual review of fluid mechanics}, 25:539--575, 1993.
\newblock \href {https://doi.org/10.1146/annurev.fl.25.010193.002543}
  {\path{doi:10.1146/annurev.fl.25.010193.002543}}.

\bibitem{BerG15}
A.~Bertram and R.~Gl{\"u}ge.
\newblock {\em Solid Mechanics Theory, Modeling, and Problems}.
\newblock Springer Cham, first edition, 2015.
\newblock \href {https://doi.org/https://doi.org/10.1007/978-3-319-19566-7}
  {\path{doi:https://doi.org/10.1007/978-3-319-19566-7}}.

\bibitem{morBil05}
D.~Billger.
\newblock The butterfly gyro.
\newblock In P.~Benner, D.~C. Sorensen, and V.~Mehrmann, editors, {\em
  Dimension Reduction of Large-Scale Systems}, volume~45 of {\em Lect. Notes
  Comput. Sci. Eng.}, pages 349--352. Springer-Verlag, Berlin/Heidelberg,
  Germany, 2005.
\newblock \href {https://doi.org/10.1007/3-540-27909-1_18}
  {\path{doi:10.1007/3-540-27909-1_18}}.

\bibitem{Boy94}
S.~Boyd, L.~El~Ghaoui, E.~Feron, and V.~Balakrishnan.
\newblock {\em Linear Matrix Inequalities in System and Control Theory}.
\newblock SIAM, 1994.
\newblock \href {https://doi.org/https://doi.org/10.1137/1.9781611970777}
  {\path{doi:https://doi.org/10.1137/1.9781611970777}}.

\bibitem{morChaLVetal06}
Y.~Chahlaoui, D.~Lemonnier, A.~Vandendorpe, and P.~Van~Dooren.
\newblock Second-order balanced truncation.
\newblock {\em Linear Algebra Appl.}, 415(2--3):373--384, 2006.
\newblock \href {https://doi.org/10.1016/j.laa.2004.03.032}
  {\path{doi:10.1016/j.laa.2004.03.032}}.

\bibitem{morCheTR12}
K.~K. Chen, J.~H. Tu, and R.~W. Rowley.
\newblock Variants of dynamic mode decomposition: Boundary condition,
  {K}oopman, and {F}ourier analyses.
\newblock {\em Nonlinear Science}, 22(6):887--915, 2012.
\newblock \href {https://doi.org/10.1007/s00332-012-9130-9}
  {\path{doi:10.1007/s00332-012-9130-9}}.

\bibitem{ChuH93}
J.~Chung and G.~M. Hulbert.
\newblock A time integration algorithm for structural dynamics with improved
  numerical dissipation: The generalized-$\alpha$ method.
\newblock {\em J. Appl. Mech.}, 60:371--375, 1993.
\newblock \href {https://doi.org/https://doi.org/10.1115/1.2900803}
  {\path{doi:https://doi.org/10.1115/1.2900803}}.

\bibitem{morCraBM68}
R.~R. Craig and M.~C.~C. Bampton.
\newblock Coupling of substructures for dynamic analyses.
\newblock {\em AIAA Journal}, 6(7):1313--1319, 1968.
\newblock \href {https://doi.org/https://doi.org/10.2514/3.4741}
  {\path{doi:https://doi.org/10.2514/3.4741}}.

\bibitem{morDav66}
E.~J. Davison.
\newblock A method for simplifying linear dynamic systems.
\newblock {\em {IEEE} Trans. Autom. Control}, AC--11:93--101, 1966.
\newblock \href {https://doi.org/10.1109/TAC.1966.1098264}
  {\path{doi:10.1109/TAC.1966.1098264}}.

\bibitem{BorCRetal12}
R.~de~Borst, M.~A. Crisfield, J.~J.~C. Remmers, and C.~V. Verhoosel.
\newblock {\em Non‐linear Finite Element Analysis of Solids and Structures}.
\newblock John Wiley and Sons, Ltd, 2012.
\newblock \href {https://doi.org/https://doi.org/10.1002/9781118375938.ch2}
  {\path{doi:https://doi.org/10.1002/9781118375938.ch2}}.

\bibitem{EicF98}
E.~Eich-Soellner and C.~F{\"{u}}hrer.
\newblock {\em Numerical Methods in Multibody Dynamics}.
\newblock Springer Fachmedien Wiesbaden GmbH, 1998.
\newblock \href {https://doi.org/https://doi.org/10.1007/978-3-663-09828-7}
  {\path{doi:https://doi.org/10.1007/978-3-663-09828-7}}.

\bibitem{FloPS14}
O.~Flod{\'{e}}n, K.~Persson, and G.~Sandberg.
\newblock Reduction methods for the dynamic analysis of substructure models of
  lightweight building structures.
\newblock {\em Comput. Struct.}, 138:49--61, 2014.
\newblock \href {https://doi.org/10.1016/j.compstruc.2014.02.011}
  {\path{doi:10.1016/j.compstruc.2014.02.011}}.

\bibitem{morFre03}
R.~W. Freund.
\newblock Model reduction methods based on {K}rylov subspaces.
\newblock {\em Acta Numer.}, 12:267--319, 2003.
\newblock \href {https://doi.org/10.1017/S0962492902000120}
  {\path{doi:10.1017/S0962492902000120}}.

\bibitem{GeR14}
M.~Geradin and D.~Rixen.
\newblock {\em Mechanical Vibrations: Theory and Application to Structural
  Dynamics}.
\newblock Hoboken, New Jersey : Wiley, 2014.
\newblock \href {https://doi.org/https://doi.org/10.1017/aer.2018.27}
  {\path{doi:https://doi.org/10.1017/aer.2018.27}}.

\bibitem{morGor94}
J.~H. Gordis.
\newblock Analysis of the improved reduced system ({IRS}) model reduction
  procedure.
\newblock {\em Int. J. Anal. Exp. Modal Anal.}, 9(4):269--285, 1994.

\bibitem{morGosG20}
I.~V. Gosea and S.~Gugercin.
\newblock The {AAA} framework for modeling linear dynamical systems with
  quadratic output.
\newblock {\em CoRR}, abs/2005.10316, 2020.
\newblock URL: \url{https://arxiv.org/abs/2005.10316}, \href
  {http://arxiv.org/abs/2005.10316} {\path{arXiv:2005.10316}}.

\bibitem{morGusS99}
B.~Gustavsen and A.~Semlyen.
\newblock Rational approximation of frequency domain responses by vector
  fitting.
\newblock {\em {IEEE} Trans. Power Del.}, 14(3):1052--1061, 1999.
\newblock \href {https://doi.org/10.1109/61.772353}
  {\path{doi:10.1109/61.772353}}.

\bibitem{morGuy65}
R.~J. Guyan.
\newblock Reduction of stiffness and mass matrices.
\newblock {\em AIAA J.}, 2:380, 1965.
\newblock \href {https://doi.org/10.2514/3.2874} {\path{doi:10.2514/3.2874}}.

\bibitem{morHaa17}
B.~Haasdonk.
\newblock Reduced basis methods for parametrized {PDE}s---a tutorial
  introduction for stationary and instationary problems.
\newblock In P.~Benner, A.~Cohen, M.~Ohlberger, and K.~Willcox, editors, {\em
  Model Reduction and Approximation: Theory and Algorithms}, pages 65--136.
  SIAM, 2017.
\newblock \href {https://doi.org/10.1137/1.9781611974829.ch2}
  {\path{doi:10.1137/1.9781611974829.ch2}}.

\bibitem{Hig88}
N.~J. Higham.
\newblock Computing a nearest symmetric positive semidefinite matrix.
\newblock {\em Linear Algebra Appl.}, 103(C):103--118, 1988.
\newblock \href {https://doi.org/10.1016/0024-3795(88)90223-6}
  {\path{doi:10.1016/0024-3795(88)90223-6}}.

\bibitem{HilHT77}
H.~M. Hilber, T.~J.~R. Hughes, and R.~L. Taylor.
\newblock Improved numerical dissipation for time integration algorithms in
  structural dynamics.
\newblock {\em Earthq. Eng. Struct. Dyn.}, 5(3):283--292, 1977.
\newblock \href {https://doi.org/https://doi.org/10.1002/eqe.4290050306}
  {\path{doi:https://doi.org/10.1002/eqe.4290050306}}.

\bibitem{morKorR05}
J.~G. Korvink and E.~B. Rudnyi.
\newblock Oberwolfach benchmark collection.
\newblock In P.~Benner, D.~C. Sorensen, and V.~Mehrmann, editors, {\em
  Dimension Reduction of Large-Scale Systems}, volume~45 of {\em Lect. Notes
  Comput. Sci. Eng.}, pages 311--315. Springer Berlin Heidelberg, 2005.
\newblock \href {https://doi.org/10.1007/3-540-27909-1_11}
  {\path{doi:10.1007/3-540-27909-1_11}}.

\bibitem{morKunV01}
K.~Kunisch and S.~Volkwein.
\newblock Galerkin proper orthogonal decomposition methods for parabolic
  systems.
\newblock {\em Numer. Math.}, 90:117--148, 2001.
\newblock \href {https://doi.org/https://doi.org/10.1007/s002110100282}
  {\path{doi:https://doi.org/10.1007/s002110100282}}.

\bibitem{morKunV08}
K.~Kunisch and S.~Volkwein.
\newblock Proper orthogonal decomposition for optimality systems.
\newblock {\em {ESAIM}: Math. Model. Numer. Anal.}, 42(1):1--23, 2008.
\newblock \href {https://doi.org/10.1051/m2an:2007054}
  {\path{doi:10.1051/m2an:2007054}}.

\bibitem{LaiRetal10}
W.~M. Lai, D.~Rubin, and E.~Krempl.
\newblock {\em Introduction to Continuum Mechanics}.
\newblock Butterworth-Heinemann, Amsterdam, 2010.
\newblock \href
  {https://doi.org/https://doi.org/10.1016/B978-0-7506-8560-3.X0001-1}
  {\path{doi:https://doi.org/10.1016/B978-0-7506-8560-3.X0001-1}}.

\bibitem{LiuMBetal13}
W.~K. Liu, B.~Moran, T.~Belytschko, and K.~Elkhodary.
\newblock {\em Nonlinear Finite Elements for Continua and Structures}.
\newblock Wiley, second edition, 2013.
\newblock URL: \url{https://books.google.de/books?id=e\_w8AgAAQBAJ}.

\bibitem{Lof04}
J.~L{\"o}fberg.
\newblock {YALMIP}: A toolbox for modeling and optimization in {MATLAB}.
\newblock In {\em IEEE International Conference on Robotics and Automation
  (IEEE Cat. No.04CH37508)}, pages 284--289, 2004.
\newblock \href {https://doi.org/10.1109/CACSD.2004.1393890}
  {\path{doi:10.1109/CACSD.2004.1393890}}.

\bibitem{morLuJCetal19}
K.~Lu, Y.~Jin, Y.~Chen, Y.~Yang, L.~Hou, Z.~Zhang, Z.~Li, and C.~Fu.
\newblock Review for order reduction based on proper orthogonal decomposition
  and outlooks of applications in mechanical systems.
\newblock {\em Mech. Syst. Signal Process.}, 123:264--297, 2019.
\newblock \href {https://doi.org/10.1016/j.ymssp.2019.01.018}
  {\path{doi:10.1016/j.ymssp.2019.01.018}}.

\bibitem{morMayA07}
A.~J. Mayo and A.~C. Antoulas.
\newblock A framework for the solution of the generalized realization problem.
\newblock {\em Linear Algebra Appl.}, 425(2--3):634--662, 2007.
\newblock \href {https://doi.org/10.1016/j.laa.2007.03.008}
  {\path{doi:10.1016/j.laa.2007.03.008}}.

\bibitem{morMoo79}
B.~C. Moore.
\newblock Principal component analysis in nonlinear systems: Preliminary
  results.
\newblock In {\em 18th IEEE Conference on Decision and Control including the
  Symposium on Adaptive Processes}, volume~2, pages 1057--1060, 1979.
\newblock \href {https://doi.org/10.1109/CDC.1979.270114}
  {\path{doi:10.1109/CDC.1979.270114}}.

\bibitem{Mue72}
P.~C. M{\"u}ller.
\newblock Stability of mechanical systems.
\newblock In {\em Special Problems of Gyrodynamics: Course Held at the
  Department of General Mechanics October 1970}, pages 34--49. Springer Vienna,
  1972.
\newblock \href {https://doi.org/10.1007/978-3-7091-2882-4}
  {\path{doi:10.1007/978-3-7091-2882-4}}.

\bibitem{NakST18}
Y.~Nakatsukasa, O.~S{\`e}te, and L.~N. Trefethen.
\newblock The {AAA} algorithm for rational approximation.
\newblock {\em SIAM J. Sci. Comput.}, 40(3):A1494--A1522, 2018.
\newblock \href {https://doi.org/10.1137/16M1106122}
  {\path{doi:10.1137/16M1106122}}.

\bibitem{New59}
N.~M. Newmark.
\newblock A method of computation for structural dynamics.
\newblock {\em ASCE J. of the Engrg. Mech. Division}, 85(EM 3, July):67--94,
  1959.
\newblock \href {https://doi.org/https://doi.org/10.1061/JMCEA3.0000098}
  {\path{doi:https://doi.org/10.1061/JMCEA3.0000098}}.

\bibitem{slicot_iss}
Niconet e.V.
\newblock {\em {SLICOT} - {S}ubroutine {L}ibrary in {S}ystems and {C}ontrol
  {T}heory}.
\newblock URL: \url{https://github.com/SLICOT}.

\bibitem{morwiki_gyro}
{Oberwolfach Benchmark Collection}.
\newblock Butterfly gyroscope.
\newblock {H}osted at {MORwiki} -- {M}odel {O}rder {R}eduction {W}iki, 2004.
\newblock URL: \url{http://modelreduction.org/index.php/Butterfly_Gyroscope}.

\bibitem{morPeh20}
B.~Peherstorfer.
\newblock Sampling low-dimensional {M}arkovian dynamics for preasymptotically
  recovering reduced models from data with operator inference.
\newblock {\em SIAM J. Sci. Comput.}, 42:A3489--A3515, 2020.
\newblock \href {https://doi.org/10.1137/19M1292448}
  {\path{doi:10.1137/19M1292448}}.

\bibitem{morPehW16}
B.~Peherstorfer and K.~Willcox.
\newblock Data-driven operator inference for nonintrusive projection-based
  model reduction.
\newblock {\em Comput. Methods Appl. Mech. Engrg.}, 306:196--215, 2016.
\newblock \href {https://doi.org/10.1016/j.cma.2016.03.025}
  {\path{doi:10.1016/j.cma.2016.03.025}}.

\bibitem{morQiaKPetal20}
E.~Qian, B.~Kr{\"a}mer, B.~Peherstorfer, and K.~Willcox.
\newblock Lift \& learn: Physics-informed machine learning for large-scale
  nonlinear dynamical systems.
\newblock {\em Physica D: Nonlinear Phenomena}, 406(1):art. 132401, 2020.
\newblock \href {https://doi.org/10.1016/j.physd.2020.132401}
  {\path{doi:10.1016/j.physd.2020.132401}}.

\bibitem{morSaaSW19}
J.~Saak, D.~Siebelts, and S.~W.~R. Werner.
\newblock A comparison of second-order model order reduction methods for an
  artificial fishtail.
\newblock {\em at-Auto\-mati\-sie\-rungs\-tech\-nik}, 67(8):648--667, 2019.
\newblock \href {https://doi.org/10.1515/auto-2019-0027}
  {\path{doi:10.1515/auto-2019-0027}}.

\bibitem{morMarLetal20}
L.~B. {Saint Martin}, R.~U. Mendes, and K.~L. Cavalca.
\newblock Model reduction and dynamic matrices extraction from state-space
  representation applied to rotating machines.
\newblock {\em Mech. Mach. Theory}, 149:103804, 2020.
\newblock \href {https://doi.org/10.1016/j.mechmachtheory.2020.103804}
  {\path{doi:10.1016/j.mechmachtheory.2020.103804}}.

\bibitem{morSalEL06}
B.~Salimbahrami, R.~Eid, and B.~Lohmann.
\newblock Model reduction by second order krylov subspaces: Extensions,
  stability and proportional damping.
\newblock In {\em Proc. IEEE Conf. Comput. Aided Control Syst. Design, Munich,
  Germany}, pages 2997--3002. IEEE, 2006.

\bibitem{morSch10}
P.~J. Schmid.
\newblock Dynamic mode decomposition of numerical and experimental data.
\newblock {\em J. Fluid Mech.}, 656:5--28, 2010.
\newblock \href {https://doi.org/10.1017/S0022112010001217}
  {\path{doi:10.1017/S0022112010001217}}.

\bibitem{morShaK22}
H.~Sharma and B.~Kramer.
\newblock Preserving lagrangian structure in data-driven reduced-order modeling
  of large-scale mechanical systems.
\newblock e-prints 2203.06361, arXiv, 2022.
\newblock URL: \url{https://arxiv.org/pdf/2203.06361.pdf}.

\bibitem{Tho21}
J.~J. Thomson.
\newblock {\em Vibrations and Stability: Advanced Theory, Analysis and Tools}.
\newblock Springer Cham, third edition, 2021.
\newblock \href {https://doi.org/https://doi.org/10.1007/978-3-030-68045-9}
  {\path{doi:https://doi.org/10.1007/978-3-030-68045-9}}.

\bibitem{TikGetal77}
A.~N. Tikhonov, A.~V. Goncharsky, V.~V. Stepanov, and A.~G. Yagola.
\newblock {\em Solutions of Ill-Posed Problems}.
\newblock Springer Dordrecht, 1995.
\newblock \href {https://doi.org/https://doi.org/10.1007/978-94-015-8480-7}
  {\path{doi:https://doi.org/10.1007/978-94-015-8480-7}}.

\bibitem{morTuRLetal14}
J.~H. Tu, C.~W. Rowley, D.~M. Luchtenburg, S.~L. Brunton, and J.~N. Kutz.
\newblock On dynamic mode decomposition: Theory and applications.
\newblock {\em J. Comput. Dynam.}, 1(2):391--421, 2014.
\newblock \href {https://doi.org/10.3934/jcd.2014.1.391}
  {\path{doi:10.3934/jcd.2014.1.391}}.

\bibitem{morVol12}
S.~Volkwein.
\newblock Model reduction using proper orthogonal decomposition.
\newblock Lecture notes, University of Konstanz, 2013.

\bibitem{morWerGG21}
S.~W.~R. Werner, I.~V. Gosea, and S.~Gugercin.
\newblock Structured vector fitting framework for mechanical systems.
\newblock e-print 2110.09220, arXiv, 2021.
\newblock URL: \url{https://arxiv.org/abs/2110.09220}.

\bibitem{morYilGBetal20}
S.~Y{\i}ld{\i}z, P.~Goyal, P.~Benner, and B.~Karas{\"o}zen.
\newblock Learning reduced-order dynamics for parametrized shallow water
  equations from data.
\newblock {\em Internat. J. Numer. Methods Fluids}, 93(8):2803--2821, 2021.
\newblock \href {https://doi.org/10.1002/fld.4998}
  {\path{doi:10.1002/fld.4998}}.

\bibitem{ZieTD14}
O.~C. Zienkiewicz, R.~L. Taylor, and D.~Fox.
\newblock {\em The Finite Element Method for Solid and Structural Mechanics}.
\newblock Butterworth-Heinemann, Oxford, seventh edition, 2014.
\newblock \href
  {https://doi.org/https://doi.org/10.1016/B978-1-85617-634-7.00017-X}
  {\path{doi:https://doi.org/10.1016/B978-1-85617-634-7.00017-X}}.

\end{thebibliography}
